\documentclass[a4paper,12pt,reqno]{amsart}

\usepackage{amsfonts}
\usepackage{amsmath}
\usepackage{stackengine}
\usepackage{amssymb}
\usepackage{cite}
\usepackage{mathrsfs}
\usepackage{enumitem}
\usepackage[colorlinks]{hyperref}
\usepackage{multicol}

\numberwithin{equation}{section}

\usepackage{graphicx}
\newcommand{\HH}{\mathbb{H}^n}

\newtheorem{theorem}{Theorem}[section]

\newtheorem{remark}[theorem]{Remark}

\newtheorem{lemma}[theorem]{Lemma}

	\allowdisplaybreaks
\begin{document}

\title[Fractional semilinear  damped wave equation on the Heisenberg group]
{Fractional semilinear  damped wave equation   on the Heisenberg group}
\author[A. Dasgupta]{Aparajita Dasgupta}
\address{
	Aparajita Dasgupta:
	\endgraf
	Department of Mathematics
	\endgraf
	Indian Institute of Technology, Delhi, Hauz Khas
	\endgraf
	New Delhi-110016, India 
	\endgraf
	{\it E-mail address:} {\rm adasgupta@maths.iitd.ac.in}
}

 \author[S. S. Mondal]{Shyam Swarup Mondal} \address{Shyam Swarup Mondal:    \endgraf Stat-Math unit \endgraf Indian Statistical Institute Kolkata \endgraf BT Road, Baranagar, Kolkata-700108, India 	\endgraf
 {\it E-mail address:} {\rm mondalshyam055@gmail.com}
}

\author[A. Tushir]{Abhilash Tushir}
\address{
	Abhilash Tushir:
	\endgraf
	Tata Institute of Fundamental Research, Centre For Applicable Mathematics
    \endgraf Bangalore, Karnataka-560065, India
	\endgraf
	{\it E-mail address:} {\rm abhilash2296@gmail.com}
}
\date{\today}
\begin{abstract}
This paper aims to investigate the Cauchy problem for the semilinear damped wave equation for the fractional  sub-Laplacian $(-\mathcal{L}_{\mathbb{H}})^{\alpha}$, $\alpha>0$  on the Heisenberg group $\HH$ with power type non-linearity. With the presence of a positive damping term and nonnegative mass term, we derive $L^2-L^2$ decay estimates for the solution of the homogeneous linear fractional damped wave equation on $\HH$, for its time derivative, and for its space derivatives. We also discuss how these estimates can be improved when we consider additional $L^1$-regularity for the Cauchy data in the absence of the mass term.  Also, in the absence of mass term, we prove the global well-posedness for $2\leq p\leq 1+\frac{2\alpha}{(\mathcal{Q}-2\alpha)}_{+}$ $(\text{or }1+\frac{4\alpha}{\mathcal{Q}}<p\leq 1+\frac{2\alpha}{(\mathcal{Q}-2\alpha)_{+}})$ in the case of $L^1\cap L^2$ $(\text{or } L^2)$ Cauchy data, respectively.   However,  in the presence of the mass term,    the global (in time) well-posedness for small data holds for $1<p \leq 1+ \frac{2\alpha}{(\mathcal{Q}-2\alpha)_{+}}$.  Finally, as an application of the linear decay estimates,  we investigate well-posedness for the Cauchy problem for a weakly coupled system with two semilinear fractional damped wave equations with positive mass term on  $\HH$.
 \end{abstract}
 
\maketitle
\tableofcontents
\section{Introduction}\label{intro}

In this paper, we study the following Cauchy problem for the semilinear damped wave equation for the fractional sub-Laplacian on the Heisenberg group $\mathbb{H}^n$ with power type non-linearity $|u|^p$, namely
\begin{align}\label{main:Heisenberg}
\begin{cases}
u_{tt}(t, \eta)+\left(-\mathcal{L}_{\mathbb{H}}\right)^{\alpha} u(t, \eta)+bu_{t}(t,\eta)+mu(t,\eta)
=|u|^{p}, &\eta\in \mathbb{H}^{n},t>0, \\
u(0, \eta)=u_{0}(\eta),\quad u_{t}(0, \eta)=u_{1}(\eta),  &\eta \in \mathbb{H}^{n},
\end{cases}
\end{align}
with the damping term determined by $b>0$ and with the mass term $m\geq 0$, where  $p>1$ and   $ (-\mathcal{L}_{\mathbb{H}})^{\alpha}$ is the fractional sub-Laplacian on $\mathbb{H}^{n}$. We refer to Section \ref{sec2} for a thorough exposition of the fractional sub-Laplacian on the Heisenberg group $\mathbb{H}^{n}$.

%More precisely, we derive decay estimates for the solution, for its time derivative, and for its space derivatives for the linear part of the damped wave equation (\ref{main:Heisenberg}), i.e.,
%\begin{align}\label{main:homoHeisenberg}
%\begin{cases}
%u_{tt}(t, \eta)+\left(-\mathcal{L}_{\mathbb{H}}\right)^{\alpha} u(t, \eta)+bu_{t}(t,\eta)+mu(t,\eta)
%=0, &\eta\in \mathbb{H}^{n},t>0, \\
%u(0, \eta)=u_{0}(\eta),\quad u_{t}(0, \eta)=u_{1}(\eta),& \eta \in \mathbb{H}^{n},
%\end{cases}
%\end{align} and as an application of the linear decay estimates for the vanishing mass term, we investigate well-posedness for the Cauchy problem for a weakly coupled
%system of semilinear fractional evolution equations on $\HH$.
%Further, we investigate the global well-posedness of the nonlinear Cauchy problem \eqref{main:Heisenberg} for the case of positive mass and in the absence of the mass term separately.

The study of linear or nonlinear PDEs in non-Euclidean structures has been the focus of several researchers in recent decades. For example, the semilinear wave equation has been investigated for the  Heisenberg group in \cite{Vla, Ruz18, MR660749}. In the case of graded groups, we refer to the recent works \cite{palmieri, gra1, Ruz18, gra3}. We also refer to \cite{garetto,27, 28, BKM22, DKM23} concerning the wave equation on compact Lie groups and \cite{AP,HWZ, AKR22,AZ} for Riemannian symmetric space of non-compact type. 

In particular, Georgiev and Palmieri in \cite{Vla} considered the following semilinear damped wave equation without mass term and with power type non-linearity $|u|^p$ on the Heisenberg group $\HH$:
\begin{align} \label{eq00100}
\begin{cases}
u_{tt}(t,\eta)-\mathcal{L}_{\mathbb{H}}u(t,\eta) +u_t(t,\eta) =|u|^p, & \eta\in \mathbb{H}^n,t>0,\\
u(0, \eta)=u_{0}(\eta),\quad u_{t}(0, \eta)=u_{1}(\eta), \quad &\eta \in \mathbb{H}^{n},
\end{cases}
\end{align}
 and proved that the critical exponent is given by
\begin{align}\label{eq36}
p_{\text{Fuji}}(\mathcal{Q}):=1+\frac{2}{\mathcal{Q}},
\end{align}
where $\mathcal{Q}:=2n+2$ represents the homogeneous dimension of $\HH$, which is also known as the so-called {\it Fujita exponent} for the semilinear heat equation; see \cite{Vla,palmieri,Nurgi}. Here the critical exponent signifies the threshold condition on the exponent $p $ for the global (in time) Sobolev solutions and the blow-up of local (in time) weak solutions with small data.
More precisely, a global existence for $p>p_{\text {Fuji}}(\mathcal{Q})$ and blow-up of the local (in time) solutions in the subcritical case $1<p<p_{\text {Fuji}}(\mathcal{Q})$. They also proved the blow-up of the weak solution for the case $p=p_{\text {Fuji}}(n).$ Its counterpart in the Euclidean framework, for the critical exponent for the semilinear damped wave equation with the power type non-linearity, one can see \cite{Matsumura, Todorova, Zhang} and references therein. Moreover,  in \cite{Palmieri2020}, Palmieri derived the $L^2$-decay estimates for solutions as well as for their time derivative and their space derivatives of solution to the homogeneous linear damped wave equation on the Heisenberg group. However, in the presence of the positive mass and damping term, the behavior of the decay for the linear damped wave equation changes drastically, as explained in \cite{Ruz18}. Furthermore, the estimates for the linear viscoelastic type damped wave equation on the Heisenberg group can be found in \cite{32}.

In this paper, we first aim to investigate the decay estimates of the linear fractional Cauchy problem for the damped wave equation \eqref{main:homoHeisenberg} without imposing any restrictions on the mass term. To the best of our knowledge, in the framework of the Heisenberg group, the fractional damped wave equation on the Heisenberg group has not been considered in the literature so far, even for the linear Cauchy problem. Therefore, an intriguing and viable problem is to investigate the behavior of the solution to the fractional damped wave equation (with a nonnegative mass term) on the Heisenberg group $\HH$. Note that the fractional power of the sub-Laplacian on the Heisenberg group $\HH$ has already been considered by several authors in literature; for example, one can see  \cite{shyam, Cowling,ahmad15} and references therein for works related to Hardy’s inequality and uncertainty inequalities for fractional powers of the sub-Laplacian on the Heisenberg group.  

 Before we provide our main results, let's review some essential terminology that will be used throughout the work.\\
\textbf{Notations:} Throughout the article,  we use the following notations: 
\begin{itemize}
	\item $f \lesssim g:$\,\, There exists a positive constant $C$ (whose value may change from line to line in this manuscript) such that $f \leq C g.$  
	\item  $f \simeq g$: Means that $f \lesssim g$ and $g \lesssim f$. 
	\item $\mathcal{Q}$:  The homogeneous dimension of $\HH$.
	\item $dg:$ The Haar measure on the Heisenberg  group $\HH.$
	\item $\mathcal{L}_{\mathbb{H}}:$ The sub-Laplacian  on $\HH.$
	\item  ${H}^{\alpha}$: The subelliptic Sobolev spaces of  order  $\alpha>0$ on $\HH$.
 \item $\|(u,v)\|_{\mathcal{D}}:$ For any normed linear space $W,X,Y,$ and $Z$, we introduce the spaces $\mathcal{D}:=(W\cap X)\times (Y\cap Z) $  with the norm
 \begin{equation*}
     \|(u,v)\|_{\mathcal{D}}:=\|u\|_{W\cap X}+\|v\|_{Y\cap Z}:=\|u\|_{W}+\|u\|_{X}+\|v\|_{Y}+\|v\|_{ Z}.
 \end{equation*}
	%\item   $(x)_{+}$: Finally, we define $(x)_{+}:=\max \{x, 0\}$ and $\frac{1}{(x)_{+}}=\infty$ when $x \leq0$.
\end{itemize}
Our first aim is to derive decay estimates for the solution, for its time derivative, and for its space derivatives for the linear part of the damped wave equation (\ref{main:Heisenberg}), i.e.,
\begin{align}\label{main:homoHeisenberg}
\begin{cases}
u_{tt}(t, \eta)+\left(-\mathcal{L}_{\mathbb{H}}\right)^{\alpha} u(t, \eta)+bu_{t}(t,\eta)+mu(t,\eta)
=0, &\eta\in \mathbb{H}^{n},t>0, \\
u(0, \eta)=u_{0}(\eta),\quad u_{t}(0, \eta)=u_{1}(\eta),& \eta \in \mathbb{H}^{n}.
\end{cases}
\end{align}
The decay estimates for the aforementioned Cauchy problem are given by the following theorem, which will be employed to further explore the global (in time) well-posedness of the Cauchy problem \eqref{main:Heisenberg} with mass and without mass.
\begin{theorem}\label{sob:well} Let $\alpha,b>0,$ and $m\geq 0$ such that $b^{2}> 4m$. Assume that the initial Cauchy data $(u_{0},u_{1})\in \mathcal{A}^{\alpha}:= H^{\alpha}(\mathbb{H}^{n})\times L^{2}(\mathbb{H}^{n})$ and let $u\in C(\mathbb{R}_{+};H^{\alpha}(\mathbb{H}^{n}))\cap C^{1}(\mathbb{R}_{+};L^{2}(\mathbb{H}^{n})) $  solve the Cauchy problem \eqref{main:homoHeisenberg}. Then the solution $u$ satisfies the following  decay estimates: 
	\begin{align}
		\|u(t, \cdot)\|_{L^2\left(\mathbb{H}^{n}\right)}&\lesssim e^{\left(-\frac{b}{2}+\sqrt{\frac{b^{2}}{4}-m}\right) t}\left\|(u_{0},u_{1})\right\|_{L^{2}\left(\mathbb{H}^{n}\right)},\label{l2l2final}\\
		\left\| u(t, \cdot)\right\|_{H^{\alpha}\left(\mathbb{H}^{n}\right)} &\lesssim (1+t)^{-\frac{1}{2}} e^{-\frac{m}{2b}t} \left\|(u_{0},u_{1})\right\|_{\mathcal{A}^{\alpha}},~~\text{and}\label{delul2l2final}\\
		\|\partial_{t}u(t, \cdot)\|_{L^2\left(\mathbb{H}^{n}\right)}&\lesssim\left[(1+t)^{-1}+m\right] e^{-\frac{m}{2b}t}\left\|(u_{0},u_{1})\right\|_{\mathcal{A}^{\alpha}},\label{deltl2l2final}
	\end{align}
	for all $t>0$. Additionally, if we suppose that the initial Cauchy data $u_{0},u_{1}\in L^{1}(\mathbb{H}^{n})$, then the following estimates hold
	\begin{align}
		\|u(t, \cdot)\|_{L^2\left(\mathbb{H}^{n}\right)}&\lesssim (1+t)^{-\frac{\mathcal{Q}}{4\alpha}}e^{-\frac{m}{2b}t}\left(\left\|(u_{0},u_{1})\right\|_{L^{1}\left(\mathbb{H}^{n}\right)}+\left\|(u_{0},u_{1})\right\|_{L^{2}\left(\mathbb{H}^{n}\right)}\right),\label{l1l2l2final}\\
		\left\| u(t, \cdot)\right\|_{H^{\alpha}\left(\mathbb{H}^{n}\right)} &\lesssim (1+t)^{-\frac{\mathcal{Q}}{4\alpha}-\frac{1}{2}} e^{-\frac{m}{2b}t} \left(\left\|(u_{0},u_{1})\right\|_{L^{1}\left(\mathbb{H}^{n}\right)}+\left\|(u_{0},u_{1})\right\|_{\mathcal{A}^{\alpha}}\right),\text{ and}\label{delul1l2final}\\
  \|\partial_{t}u(t, \cdot)\|_{L^2\left(\mathbb{H}^{n}\right)}&\lesssim\left[(1+t)^{-\frac{\mathcal{Q}}{4\alpha}-1}+m(1+t)^{-\frac{\mathcal{Q}}{4\alpha}}\right] e^{-\frac{m}{2b}t}\left(\left\|(u_{0},u_{1})\right\|_{L^{1}\left(\mathbb{H}^{n}\right)}+\left\|(u_{0},u_{1})\right\|_{\mathcal{A}^{\alpha}}\right),\label{deltl1l2l2final}
	\end{align}
	for all $t\geq0$. 
\end{theorem}
The estimates obtained for the linear fractional damped wave equation in the preceding theorem expand various previous results in the literature. The estimates for the damped wave equation with mass term $(m>0)$ and without mass term $(m=0)$  were developed separately in  \cite{Ruz18} and \cite{Palmieri2020},  respectively, for the Heisenberg group $\HH$. Our analysis used in the above theorem allows us to obtain the estimates unitedly for mass as well as massless scenarios on the Heisenberg group but for fractional cases. 
 \begin{remark}
      We can observe the following important fact from  Theorem \ref{sob:well}. 
 \begin{itemize}
 %\item  Given $\alpha=1$, $b>0$, and $m=0$, all the decay estimates in the above theorem coincide with the decay estimates proved in \cite{Palmieri2020}. On the other hand, for $\alpha=1$, $b>0$, and $m>0$,  the decay estimates  in \eqref{l2l2final}, \eqref{delul2l2final}, and \eqref{deltl2l2final},  coincide with the decay estimates obtained in \cite{Ruz18}.
 	\item  From \eqref{l2l2final}, \eqref{delul2l2final}, and \eqref{deltl2l2final}, one can observe the existence of an exponential decay rate for the $L^2$-norms of the solution of the homogeneous linear problem and its derivatives. This is due to the presence of a mass term in the linear damped wave equation \eqref{main:homoHeisenberg}.
 However, the situation changes drastically if we remove the mass term from the fractional linear equation \eqref{main:homoHeisenberg}. In this case, we cannot get exponential decay; in fact, there is no decay in \eqref{l2l2final} for the $L^2$-norms of the solution. In order to find still some decay, further we need to consider additional $L^1$-regularity for the initial Cauchy data, as can be seen in \eqref{l1l2l2final}.
\item Moreover, the above estimate provides an explicit dependence of the decay estimates over the parameters $b$ and $m$, which aids in extracting the polynomial decay from the exponential decay in the absence of a mass term.  In the light of relation $-\frac{b}{2}+\sqrt{\frac{b^{2}}{4}-m}=-\frac{m}{2b}-f(b,m)$, where $f(b,m)=\frac{b^{2}-m}{2b}-\sqrt{\frac{b^{2}}{4}-m}\geq 0$ (can be verify using contradiction). The estimate \eqref{l2l2final} can be further estimated as:
\begin{equation}\label{weak:l2l2}
 \|u(t, \cdot)\|_{L^2\left(\mathbb{H}^{n}\right)}\lesssim e^{-\frac{m}{2b}t}\left\|(u_{0},u_{1})\right\|_{L^{2}\left(\mathbb{H}^{n}\right)}.
\end{equation}
Despite being weaker than \eqref{l2l2final}, the aforementioned estimate will facilitate the computation in the proof of Theorem \ref{Main:m>0} and the coupled system.
%\item \textcolor{red}{In the massless scenario, we can split the relevant integral term from the Plancherel formula into two zones with regard to $\lambda$. For ``small'' $|\lambda|$, one can use the analogue of Riemann-Lebesgue inequality in order to get polynomial decay rates as in  \eqref{l1l2l2final}, \eqref{delul1l2final}, and  \eqref{deltl1l2l2final}. On the other hand, for $|\lambda|$ ``large" by Plancherel formula, we get exponential decay, provided the Cauchy data has adequate $L^2$-regularity.} 
\item The estimates in the above theorem correspond identically to the $L^2-L^2$ estimates with additional $L^1$-regularity for the data in the Euclidean situation was proven in \cite{Matsumura} with $\alpha=1$.  Furthermore, the fractional damped wave equation in Euclidean setting, as shown in \cite[Proposition 2.1]{ressing17} and \cite[Theorem 1.4]{radu}, also corresponds exactly to these estimates.
\end{itemize}
 \end{remark}
Coming back to the nonlinear case,  as we explained earlier, the critical exponent for the damped wave equation \eqref{eq00100} is given by $ p_{\text{Fuji}}(\mathcal{Q}):=1+\frac{2}{\mathcal{Q}}$, see \cite{Vla}. More precisely, they proved a global existence for $p>p_{\text {Fuji}}(\mathcal{Q})$ and blow-up of the local (in time) solutions in the subcritical case $1<p<p_{\text {Fuji}}(\mathcal{Q})$ for the equation \eqref{eq00100}. On the other hand,  in the presence of positive mass and damping term,  the authors in \cite{Ruz18} considered the Cauchy problem for the semilinear damped wave equation for the sub-Laplacian on the Heisenberg group (and the more general setting of graded Lie groups for the semilinear damped wave equation involving a Rockland operator)
with power type non-linearities and proved global (in time) well-posedness for small initial data.
 
One next purpose of this paper is to investigate the global (in time) well-posedness of the Cauchy problem for the fractional semilinear damped wave equation \eqref{main:Heisenberg}. In the presence of  damping term   $b>0$  and in the absence of the mass term $m$ from the equation \eqref{main:Heisenberg}, i.e.,   for the Cauchy problem  
\begin{align}\label{non:Heisenberg}
	\begin{cases}
		u_{tt}(t, \eta)+\left(-\mathcal{L}_{\mathbb{H}}\right)^{\alpha} u(t, \eta)+bu_{t}(t,\eta)
		=|u|^{p},  &\eta\in \mathbb{H}^{n},t>0, \\
		u(0, \eta)=u_{0}(\eta),\quad u_{t}(0, \eta)=u_{1}(\eta),  &\eta \in  \mathbb{H}^{n},
	\end{cases}
\end{align} 
for $p>1$, we have the following  global (in time) well-posedness result.
\begin{theorem}[Global Existence: $m=0$]\label{global:existence}
	Let $\alpha\geq 1$. Let $2\alpha<\mathcal{Q}\leq4\alpha $ and fix $p$ in Cauchy problem \eqref{non:Heisenberg} such that
	\begin{equation}\label{GNI}
		2\leq p\leq  1+\frac{2\alpha}{\mathcal{Q}-2\alpha}.
	\end{equation}
	Then, there exists a sufficiently small $\varepsilon>0$ such that for any initial data
	\begin{equation*}
		(u_{0},u_{1})\in \mathcal{B}^{\alpha}:= L^{1}\left(\mathbb{H}^{n}\right)\cap H^{\alpha}\left(\mathbb{H}^{n}\right)\times L^{1}\left(\mathbb{H}^{n}\right)\cap L^{2}\left(\mathbb{H}^{n}\right) ~~\text{ satisfying }~~  		\|(u_{0},u_{1})\|_{\mathcal{B}^{\alpha}}<\varepsilon,
	\end{equation*}
	there is a unique global solution $u \in C\left(\mathbb{R}_{+}; H^\alpha\left(\mathbb{H}^{n}\right)\right) \cap C^1\left(\mathbb{R}_{+}; L^2\left(\mathbb{H}^{n}\right)\right)$ to the Cauchy problem \eqref{non:Heisenberg}. Moreover, u satisfies the following estimates
	   \begin{equation*}
		\left\|\partial_{t}^{i}(-\mathcal{L}_{\mathbb{H}})^{j/2} u(t, \cdot)\right\|_{L^2\left(\mathbb{H}^{n}\right)} \lesssim (1+t)^{-\frac{\mathcal{Q}}{4\alpha}-\frac{j}{2\alpha}-i}  \|(u_{0},u_{1})\|_{\mathcal{B}^{\alpha}} ,
	\end{equation*}
for $(i,j)\in\{(0,0),(0,\alpha),(1,0)\}$ and for all $t\geq 0$.
\end{theorem}
\begin{remark}\label{l1:dropped}
The restriction $p\geq 2$ in \eqref{GNI} is due to the additional $L^1$ assumption on the initial data, together with the hypothesis of the Gagliardo-Nirenberg type inequality on $\HH$. The restriction $p\geq 2$ for $2\alpha<\mathcal{Q}\leq 4\alpha$ corresponds identically to the Euclidean case; see \cite[Section 18.1]{reissig:book}.
\end{remark}
As we observed in Theorem \ref{sob:well} that relaxation of initial Cauchy data from $L^{1}\cap L^{2}$-regularity to $ L^{2}$ leads to the loss of decay rate. Furthermore, it is known that in the Euclidean situation, the critical exponent for semilinear problems changes if the initial data's $L^1$ assumption is replaced by a $L^p$ assumption for $p\in(1,2]$; see \cite{ryo}. Analogous to the Euclidean case \cite{abbico2017}, the dropping of $L^1$-regularity also affects the range of $p$ for global existence for the Cauchy problem \eqref{non:Heisenberg}. 
More precisely, the following theorem gives global (in time) well-posedness result in the case of $L^2$-regularity only.
\begin{theorem}[Global Existence: $m=0$]\label{global:existence:l2}
	Let $\alpha> 1$. Let $2\alpha<\mathcal{Q}<4\alpha $ and fix $p$ in Cauchy problem \eqref{non:Heisenberg} such that
	\begin{equation}\label{GNII}
		1+\frac{4\alpha}{\mathcal{Q}}< p\leq  1+\frac{2\alpha}{\mathcal{Q}-2\alpha}.
	\end{equation}
	Then, there exists a sufficiently small $\varepsilon>0$ such that for any initial data
	\begin{equation*}
		(u_{0},u_{1})\in \mathcal{A}^{\alpha}=  H^{\alpha}\left(\mathbb{H}^{n}\right)\times  L^{2}\left(\mathbb{H}^{n}\right) ~~\text{ satisfying }  ~~		\|(u_{0},u_{1})\|_{\mathcal{A}^{\alpha}}<\varepsilon,
	\end{equation*}
	there is a unique global solution $u \in C\left(\mathbb{R}_{+}; H^\alpha\left(\mathbb{H}^{n}\right)\right) \cap C^1\left(\mathbb{R}_{+}; L^2\left(\mathbb{H}^{n}\right)\right)$ to the Cauchy problem \eqref{non:Heisenberg}. Moreover, u satisfies the following estimates
  \begin{equation*}
		\left\|\partial_{t}^{i}(-\mathcal{L}_{\mathbb{H}})^{j/2} u(t, \cdot)\right\|_{L^2\left(\mathbb{H}^{n}\right)} \lesssim (1+t)^{-\frac{j}{2\alpha}-i}  \|(u_{0},u_{1})\|_{\mathcal{A}^{\alpha}} ,
	\end{equation*}
for $(i,j)\in\{(0,0),(0,\alpha),(1,0)\}$ and for all $t\geq 0$.
\end{theorem}
\begin{remark} The following justifies the technical assumption in Theorem \ref{global:existence} and Theorem \ref{global:existence:l2}:
\begin{itemize}
    \item $\alpha\geq 1$ or $\alpha>1$  is required to guarantee the existence of $\mathcal{Q}$ satisfying
\begin{equation*}
    2\alpha<\mathcal{Q}\leq 4\alpha \quad \text{and}\quad     2\alpha<\mathcal{Q}< 4\alpha,
\end{equation*}
respectively, since $\mathcal{Q}=2n+2\geq 4$. 
\item  $\mathcal{Q}>2\alpha$ is due to the application of Gagliardo-Nirenberg inequality.
  \item   $\mathcal{Q}\leq 4\alpha$ or $\mathcal{Q}<4\alpha$  is required to make the range of $p$ in \eqref{GNI} and \eqref{GNII} nonempty, respectively.
\end{itemize}
\end{remark}
\begin{remark}
    In the massless case, we prove global (in time) well-posedness for $p>1+\frac{4\alpha}{\mathcal{Q}}$ with $L^2$-regularity on the initial data and $p\geq 2$ with $L^1\cap L^2$-regularity on the initial data. Whereas,    for $\alpha=1$, Georgiev and 
 Palmieri in \cite{Vla} prove the global (in time) existence for $p>p_{\text{Fuji}}(\mathcal{Q}):=1+\frac{2}{\mathcal{Q}}$ in the exponential weighted energy space, which is a stronger assumption than $L^1$-regularity.% Under additional restrictions on the initial Cauchy data or using $L^p$-estimates$p\in[1,2]$, we believe it is possible to prove the global existence for $p>1+\frac{2\alpha}{\mathcal{Q}}$. Moreover, it is also possible to prove the local well-posedness result for the missing range of $p$.  
\end{remark}
In \cite{Ruz18}, Ruzhansky and Tokmagambetov studied the semilinear damped wave equation \eqref{main:Heisenberg} with positive mass in a wider framework of graded Lie groups. They observed the global (in time) well-posedness for $p>1$ because of the exponential decay in the presence of positive mass and initial data having $L^2$-regularity only. No further lower bound for the non-linearity exponent $p > 1$ is necessary in this case. 
The decay estimates obtained in \cite[Theorem 4.6]{Ruz18} for the nonlinear Cauchy problem \eqref{main:Heisenberg} with positive mass can be improved by explicit dependency of decay estimates in Theorem \ref{sob:well} on the parameters $\alpha,b,$ and $m$.
\begin{theorem}[Global Existence: $m>0$]\label{Main:m>0}
	Let $\alpha>0$, $\mathcal{Q}>2\alpha$ and $1<p \leq 1+ \frac{2\alpha}{\mathcal{Q}-2\alpha}$.
 Then, there exists a sufficiently small $\varepsilon>0$ such that for any initial data 
 \begin{equation*}
 (u_{0},u_{1})\in \mathcal{A}^{\alpha}= H^{\alpha}(\mathbb{H}^n) \times L^{2}\left(\mathbb{H}^{n}\right)~~\text{ satisfying }~~   \|(u_{0},u_{1})\|_{\mathcal{A}^{\alpha}}<\varepsilon, 
	\end{equation*}
	there is a unique global solution $u \in C\left(\mathbb{R}_{+}; H^\alpha\left(\mathbb{H}^{n}\right)\right) \cap C^1\left(\mathbb{R}_{+}; L^2\left(\mathbb{H}^{n}\right)\right)$ to the Cauchy problem \eqref{main:Heisenberg}.  
	Moreover, the solution  $u$ satisfies the following estimates
 \begin{equation*}
		\left\|\partial_{t}^{i}(-\mathcal{L}_{\mathbb{H}})^{j/2} u(t, \cdot)\right\|_{L^2\left(\mathbb{H}^{n}\right)} \lesssim \left[(1+t)^{-\frac{j}{2\alpha}-i}+mi\right] e^{-\frac{m}{2b}t}  \|(u_{0},u_{1})\|_{\mathcal{A}^{\alpha}},
	\end{equation*}
	for $(i,j)\in\{(0,0),(0,\alpha),(1,0)\}$ and for all $t>0$.
\end{theorem} 
\begin{remark}
 The upper bound $\mathcal{Q}<4\alpha$ is now relaxed with the inclusion of the mass term, since the assumption $\mathcal{Q}>2\alpha$ is enough to guarantee the existence of $p$ satisfying $1<p \leq 1+ \frac{2\alpha}{\mathcal{Q}-2\alpha}$.
\end{remark} 
  
In addition to the paper's main findings, we will examine the well-posedness of the weakly coupled system on the Heisenberg group $\mathbb{H}^{n}$ as an application of our linear estimates found in Theorem \ref{sob:well}.

Apart from the introduction, let's provide a quick outline of the paper's structure for the convenience of reading:
 \begin{itemize}
 \item In Section \ref{sec2},  we recall some basics of the Fourier analysis on the Heisenberg groups $\HH$ to
make the manuscript self-contained.  More precisely,  in Subsection \ref{sec2.3}, we recall the Fourier analysis on the Heisenberg group. In Subsection \ref{sec2.1}, we recall the Hermite functions and Hermite operator, which act as a building block to define the fractional sub-Laplacian on the Heisenberg group. Also, in Subsection \ref{sec2.2}, we give a thorough exposition of the fractional sub-Laplacian $(-\mathcal{L}_{\mathbb{H}})^{\alpha}$ and the associated Sobolev spaces on the Heisenberg group.

  \item in Section \ref{sec3}, we first develop a preparatory lemma to prove Theorem \ref{sob:well}. Then we develop the estimates for the solution, its space derivative and time derivative in Subsections \ref{subsec3.1}, \ref{subsec3.2}, and \ref{subsec3.3}, respectively, to conclude the proof of Theorem \ref{sob:well}. 
  %   \item  In Section \ref{sec3}, by employing the group Fourier transform and using the Plancherel formula for low and high-frequency parts, we prove decay estimates for the  linear Cauchy problem 
 %\eqref{main:homoHeisenberg}. In the absence of the mass term, we also discuss the decay of the linear equation by considering additional $L^1$-regularity for the initial Cauchy data.  
     \item  in Section \ref{sec4}, we supply the  proof of Theorem \ref{global:existence} and  \ref{global:existence:l2}.
     \item in Section \ref{sec5}, as an application of the linear estimate, we investigate the well-posedness results for the   weakly coupled system and consequently we will shed some light regarding the proof of Theorem \ref{Main:m>0}.
     \item in Section \ref{sec6}, we conclude the manuscript with some final remarks.
 \end{itemize}

\section{Preliminaries}\label{sec2}
In this section, we recall some basics of the Fourier analysis on the Heisenberg  groups $\mathbb{H}^n$ to make the manuscript self-contained. A complete account of the representation theory on $\mathbb{H}^n$ can be found in 
\cite{Vla,Ruz18,Fischer,shyam,shyam1}. 
However, we mainly adopt the notation and terminology given in \cite{Fischer} for convenience.

\subsection{Fourier analysis on the Heisenberg group $\mathbb{H}^n$}\label{sec2.3} 
One of the simplest examples of a non-commutative and non-compact group is the famous Heisenberg group $\mathbb{H}^n$.  The theory of the Heisenberg group   plays a crucial  role in several branches of mathematics and physics. The Heisenberg group $\mathbb{H}^n$  is a nilpotent Lie group whose underlying manifold is $ \mathbb{R}^{2n+1} $ and the group operation is defined by
$$(x, y, t)\circ (x', y', t')=(x+x', y+y', t+t'+ \frac{1}{2}(xy'-x'y)),$$
where $(x, y, t)$, $ (x', y', t')$ are in $\mathbb{R}^n \times \mathbb{R}^n \times \mathbb{R}$ and  $xy'$ denotes the standard scalar product in $\mathbb{R}^n$. Moreover, $\mathbb{H}^n$ is a unimodular Lie group on which the left-invariant Haar measure $\mathrm{d}g$ is the usual Lebesgue measure $\, \mathrm{d}x \, \mathrm{d}y \,\mathrm{d}t.$ 

We start this subsection by  recalling the definition of the operator-valued group Fourier transform on $\mathbb{H}^n$.  By Stone-Von Neumann theorem, the only infinite-dimensional unitary irreducible representations
(up to unitary equivalence) are given by $\pi_{\lambda}$, $\lambda$ in $\mathbb{R}^*$, where the mapping  $\pi_{\lambda}$  is a strongly continuous unitary representation defined by 
$$\pi_{\lambda}(\eta)f(u)=e^{i \lambda(t+\frac{1}{2}xy)} e^{i \sqrt{\lambda}yu} {  f(u+\sqrt{|\lambda|}x)}, \quad    \eta=(x, y, t)\in \mathbb{H}^n,$$ for all  $ f  \in L^2(\mathbb{R}^n).$ We use the convention $$\sqrt{\lambda}:={\rm sgn}(\lambda)\sqrt{|\lambda|}
= \begin{cases}  \sqrt{\lambda},&\lambda>0,\\ -\sqrt{|\lambda|}, &\lambda<0.\end{cases}
$$
For each ${\lambda} \in \mathbb{R}^*$, the group Fourier transform of $f\in L^1(\mathbb{H}^n)$  is a bounded linear operator   on  $L^2(\mathbb{R}^n)$ defined by
\begin{equation}\label{ftheisen}
 \widehat{f}(\lambda)=\int_{\mathbb{H}^{n}}f(\eta)\pi_{\lambda}(\eta)^{*}\operatorname{d}\eta,\quad \lambda\in\mathbb{R}^{*}
\end{equation}
Let   $B(L^2(\mathbb{R}^n))$  be the set of all bounded operators on $L^2(\mathbb{R}^n)$. As the Schr\"odinger representations are unitary,  for any  {  $\lambda\in \mathbb{R}^*$}, we have  
\begin{align}\label{fhatnorm}
 \|\widehat{f}\|_{B(L^2(\mathbb{R}^n))}\leq \|f\|_{L^1(\HH)}.
\end{align}
If $f \in L^2(\mathbb{H}^n)$, then $\widehat{f}(\lambda)$ is a Hilbert-Schmidt operator on $L^2(\mathbb{R}^n)$ and satisfies the following Plancherel formula $$\|{f}\|_{L^2(\mathbb{H}^n)}^2=\int_{\mathbb{R}^*}\|\widehat{f}(\lambda)\|_{S_2}^{2} \, \mathrm{d}\mu( \lambda),$$
where $\| . \|_{S_2}$ stands for  the norm in the Hilbert space $S_2$,   the set of all Hilbert-Schmidt operators  on  $L^2(\mathbb{R}^n)$ and $\mathrm{d}\mu(\lambda)=c_n {|\lambda|}^n \, \mathrm{d}\lambda$ with $c_n$ being  a positive constant. Thus, with the help of the orthonormal basis $\left\{e_k\right\}_{k \in \mathbb{N}_0^n}$ for  $L^2\left(\mathbb{R}^n\right)$, we have  the following  Parseval's identity 
\begin{equation}\label{parseval}
	\left\|\widehat{f}(\lambda)e_{k}\right\|^{2}_{L^{2}(\mathbb{R}^{n})}=\sum\limits_{l\in\mathbb{N}^{n}}\left|(\widehat{f}(\lambda)e_{k},e_{\ell})_{L^{2}(\mathbb{R}^{n})}\right|^{2}
\end{equation}
and  using this,  the definition of the Hilbert-Schmidt norm can be formulate as 
\begin{equation*}
	\|\widehat{f}(\lambda)\|_{S_2}^2  \doteq \operatorname{Tr}\left[  (\pi_\lambda(f) )^* \pi_\lambda(f)\right]  =\sum_{k \in \mathbb{N}_0^n} \|\widehat{f}(\lambda) e_k \|_{L^2\left(\mathbb{R}^n\right)}^2=\sum_{k, \ell \in \mathbb{N}^n_0} | (\widehat{f}(\lambda) e_k, e_{\ell} )_{L^2\left(\mathbb{R}^n\right)} |^2.
\end{equation*}
The above expression allows us to write the Plancherel formula  in the following way:
\begin{align}\label{Plancherel}
	\|{f}\|_{L^2(\mathbb{H}^n)}^2= c_n\int_{\mathbb{R}^*}\sum_{k, \ell \in \mathbb{N}^n_0} | (\widehat{f}(\lambda) e_k, e_{\ell} )_{L^2\left(\mathbb{R}^n\right)} |^2 \, {|\lambda|}^n \, \mathrm{d}\lambda. 
\end{align} 
For  $f \in \mathcal{S}(\mathbb{H}^n)$, the space of all Schwartz class functions on $\HH,$ the Fourier inversion formula  takes the form
\begin{align*}
	f(\eta) &:=  \int_{\mathbb{R}^*}\operatorname{Tr}[\pi_{\lambda}(\eta)\widehat{f}(\lambda)]\,\mathrm{d}\mu(\lambda), \quad   \eta\in \mathbb{H}^n,
\end{align*}
where $\operatorname{Tr(A)}$ denotes the trace  of the operator $A$.
%Let $f^\lambda$ stand for the inverse Fourier transform of $f$ with respect to the   central variable $t$ given by 
%$$f^\lambda(x, y)=\int_{\mathbb{R}^*} f(x, y, t) e^{i\lambda t} dt.$$
%Now taking usual Euclidean Fourier transform of $f^\lambda$ in the variable $\lambda$, we get
%$$f(x, y, t)=\frac{1}{2\pi}\int_{\mathbb{R}^*} f^\lambda(x, y) e^{-i\lambda t} \mathrm{d}\lambda.$$
%Furthermore,  the convolution of two function $f$ and  $g$ on  the Heisenberg group $\HH$ is defined by
%$$f\ast g(\eta)=\int_{\HH} f(\eta \zeta^{-1})g(\zeta)\; d \zeta, \quad \eta, \zeta\in \HH.$$
%With $\HH=\mathbb{C}^n\times \mathbb{R}$ in view and writing  an element $\eta=(z, t)\in \mathbb{C}^n\times \mathbb{R}$,  a straightforward calculation yields
%\begin{align*} 
%    (f\ast g)^\lambda (z)&=\int_{\mathbb{C}^{n}}f^\lambda(z-z')g^\lambda (z') e^{\frac{i}{2} \operatorname{Im(z\cdot \bar{z'})}} \;dz', \quad z, z'\in \mathbb{C}^{n}\\
%    &=f^\lambda \ast_\lambda g^\lambda (z),
%\end{align*}
%where $$f^\lambda \ast_\lambda g^\lambda (z)=\int_{\mathbb{C}^{n}}f^\lambda(z-z')g^\lambda (z') e^{\frac{i}{2} \operatorname{Im(z\cdot \bar{z'})}} \;dz', \quad z, z'\in \mathbb{C}^{n}$$ is called the $\lambda$-twisted convolution.}

\subsection{Hermite operator and its properties}\label{sec2.1} In this  subsection, we recall some  definitions and properties of Hermite functions which we will use frequently in order to study Schr\"odinger representations and sub-Laplacian on  the Heisenberg group $\HH$.  We start with the definition of Hermite polynomials on $\mathbb{R}.$

Let $H_k$ denote the Hermite polynomial on $\mathbb{R}$, defined by
$$H_k(x)=(-1)^k \frac{d^k}{dx^k}(e^{-x^2} )e^{x^2}, \quad k=0, 1, 2, \dots   ,$$\vspace{.30mm}
and $h_k$ denote the normalized Hermite functions on $\mathbb{R}$ defined by
$$h_k(x)=(2^k\sqrt{\pi} k!)^{-\frac{1}{2}} H_k(x)e^{-\frac{1}{2}x^2}, \quad k=0, 1, 2, \dots,$$\vspace{.30mm}
The Hermite functions $\{h_k \}$ are the eigenfunctions of the Hermite operator (or the one-dimensional harmonic oscillator) $\mathrm{H}=-\frac{d^2}{dx^2}+x^2$ with eigenvalues $2k+1,  k=0, 1, 2, \dots$. These functions form an orthonormal basis for $L^2(\mathbb{R})$. The higher dimensional Hermite functions denoted by $e_{k}$ are then obtained by taking tensor products of one dimensional Hermite functions. Thus for any multi-index $k=(k_1, \dots, k_n) \in \mathbb{N}_0^n$ and $x =(x_1, \dots, x_n)\in \mathbb{R}^n$, we define
$e_{k}(x)=\prod_{j=1}^{n}h_{k_j}(x_j).$
The family $\{e_{k}\}_{k\in \mathbb{N}^n_0}$ is then   an orthonormal basis for $L^2(\mathbb{R}^n)$. They are eigenfunctions of the Hermite operator $\mathrm{H}=-\Delta+|x|^2$, namely, we have an ordered set of  {  natural  numbers} $\{\mu_k\}_{k \in \mathbb{N}_0^n}$ such that 
\begin{equation*}
  \mathrm{H} e_{k}(x)=\mu_ke_{k}(x), \quad x\in \mathbb{R}^n,  
\end{equation*}
for all $k\in \mathbb{N}_0^n$.  More precisely, $\mathrm{H}$  has eigenvalues 
$$\mu_k=\sum_{j=1}^{n}(2k_j+1)=2|k|+n,$$
corresponding to the eigenfunction $e_{k}$ for $k\in \mathbb{N}_0^n.$

  Given $f \in L^{2}(\mathbb{R}^{n})$, we have the  following Hermite expansion
	$$f=\sum_{k \in \mathbb{N}_0^{n}}\left(f, e_{k}\right) e_{k}=\sum_{m=0}^\infty \sum_{|k|=m}\left(f, e_{k}\right) e_{k}= \sum_{m=0}^\infty P_mf,$$ where $P_{m}$ denotes the orthogonal projection of $L^{2}(\mathbb{R}^{n})$ onto the eigenspace spanned by $\{e_{k}:|k|=m\}.$  
	%Observe the regularising effect of $P_m$, which takes temperate distributions into Schwartz functions: $P_m: \mathcal{S}^{\prime}\left(\mathbb{R}^n\right) \rightarrow \mathcal{S}\left(\mathbb{R}^n\right)$.  
	Then the spectral decomposition of $\mathrm{H}$ on $\mathbb{R}^n$ is given by
	$$
	\mathrm{H}f=\sum_{m=0}^{\infty}(2 m+n) P_mf.
	$$
	Since 0 is not in the spectrum of $H$,  for any $s \in \mathbb{R}$, we can define the fractional powers $H^s$  by means of the spectral theorem, namely
	\begin{align}\label{Hermite power}
		\mathrm{H}^s f=\sum_{m=0}^{\infty}(2 m+n)^s P_m f.
	\end{align}
    
%On the other hand,  for $\lambda\in \mathbb{R}^*$, one can  define the scaled Laguerre functions of type $(n-1)$ as$$\Phi_k^\lambda(z)=L_k^{n-1} \left( \frac{1}{2} |\lambda||z|^2\right)e^{-\frac{1}{4}|\lambda||z|^2}, \quad z\in \mathbb{C}^n.$$ Here $L_k^{n-1}$'s are the Laguerre polynomials of type $(n-1)$. For the definition and several important properties of the Laguerre polynomials, we refer to   \cite{shyam}. Moreover,  the set $\{\Phi_k^\lambda\}_{k=0}^\infty$ forms an orthogonal basis for the subspace consisting of radial functions in $L^2(\mathbb{C}^n).$

 \subsection{Sub-Laplacian and Sobolev spaces on the Heisenberg group $\HH$}\label{sec2.2}
 The canonical basis for the Lie algebra $\mathfrak{h}_n$ of $\mathbb{H}^n$ is given by the left-invariant vector fields:
\begin{align}\label{CH00VF}
	X_j&=\partial _{x_j}-\frac{y_j}{2} \partial _{t}, &Y_j=\partial _{y_j}+\frac{x_j}{2} \partial _{t}, \quad j=1, 2, \dotsc   n,   ~ \mathrm{and} ~ T=\partial _{t},
\end{align} \vspace{.30mm}
which 	satisfy the commutator relations
$[X_{i}, Y_{j}]=\delta_{ij}T, \quad  \text{for} ~i, j=1, 2, \dotsc n.$

Moreover, the canonical basis for   $\mathfrak{h}_n$  admits the  decomposition  $\mathfrak{h}_n=V\oplus W,$ where $V=\operatorname{span}\{X_j, Y_j\}_{j=1}^n$ and $W=\operatorname{span} \{T\}$.   Thus, the Heisenberg group $\HH$ is a step 2 stratified Lie group, and its homogeneous dimension is $Q:=2n +2$.
The sub-Laplacian $\mathcal{L}_{\mathbb{H}}$ on $\HH$ is   defined  as
\begin{align*}
	\mathcal{L}_{\mathbb{H}} &=\sum_{j=1}^{n}(X_j^2+Y_j^2)=\sum_{j=1}^{n}\bigg(\bigg(\partial _{x_j}-\frac{y_j}{2} \partial _{t}\bigg)^2+\bigg(\partial _{y_j}+\frac{x_j}{2} \partial _{t}\bigg)^2\bigg)\\& =\Delta_{\mathbb{R}^{2n}}+\frac{1}{4}(|x|^2+|y|^2) \partial_t^2+\sum_{j=1}^n\left(x_j \partial_{y_j }-y_j \partial_{x_j}\right)\partial_{t},
\end{align*} 
where $\Delta_{\mathbb{R}^{2n}}$ is the standard Laplacian on $\mathbb{R}^{2n},$ $|x|^{2}=x_{1}^{2}+\cdots+x_{n}^{2}$, and $|y|^{2}=y_{1}^{2}+\cdots+y_{n}^{2}$.
It is known that  the operator $\mathcal{L}_{\mathbb{H}}$ is a non-negative hypoelliptic differential operator, which is essentially self-adjoint
on the Schwartz class $\mathcal{S}(\mathbb{H}^n) = \mathcal{S}(\mathbb{R}^{2n+1})$ and therefore  admits a unique self-adjoint extension on $L^2(\HH)$ with   maximal domain $\mathcal{D}(\mathcal{L}_{\mathbb{H}})=\{f\in L^2(\HH):\mathcal{L}_{\mathbb{H}}f \in L^2(\HH) \}$.

%Then, using the spectral decomposition of the scaled Hermite operator $H(\lambda)=-\Delta+|\lambda|^2|x|^2$ and the $\lambda$-twisted convolution,  the spectral decomposition of the sub-Laplacian $\mathcal{L}_\mathbb{H}$ on $\HH$ is given by $$ (-\mathcal{L}_{\mathbb{H}}) f(z, t)=(2\pi)^{-n-1}\int_{\mathbb{R}^*} \left ( \sum_{k=0}^\infty (2k+n)|\lambda| f^\lambda \ast_\lambda \Phi_k^\lambda (z)\right)e^{-it\lambda} \;d\mu(\lambda), $$ where $\Phi_k^\lambda(z) $  is the  scaled Laguerre functions of type $(n-1)$.  Moreover, the fractional powers of the sub-Laplacian $\mathcal{L}_\mathbb{H}^s, s>0$ can be equivalently defined via  spectral decomposition\begin{align}\label{def} (-\mathcal{L}_{\mathbb{H}} )^sf(z, t)=(2\pi)^{-n-1}\int_{\mathbb{R}^*} \left ( \sum_{k=0}^\infty \left((2k+n)|\lambda|\right)^{s} f^\lambda \ast_\lambda \Phi_k^\lambda (z)\right)e^{-it\lambda} \;d\mu(\lambda). \end{align}
%Furthermore, one can notice that  $\widehat{\mathcal{L}_{\mathbb{H}} ^sf}(\lambda)=\widehat{f}(\lambda) H(\lambda)^s. $

On the other hand, the action of the infinitesimal representation $\mathrm{d}\pi_\lambda$ of $\pi_\lambda$ on the generators of the first layer of the Lie algebra $\mathfrak{h}_n$ is given by
\begin{equation*}
    \mathrm{d} \pi_\lambda\left(X_j\right)=\sqrt{|\lambda|} \partial_{x_j}\quad \text{and}\quad\mathrm{d} \pi_\lambda\left(Y_j\right)=i \operatorname{sign}(\lambda) \sqrt{|\lambda|} x_j \text { for } j=1, \dots, n .
\end{equation*}
Since the action of $\;\mathrm{d} \pi_\lambda$ can be extended to the universal enveloping algebra of $\mathfrak{h}_n$, combining  the above two expressions,  we obtain
$$
\;\mathrm{d}\pi_\lambda\left(\mathcal{L}_{\mathbb{H}} \right)=\;\mathrm{d} \pi_\lambda\bigg(\sum_{j=1}^n(X_j^2+Y_j^2)\bigg)=|\lambda| \sum_{j=1}^n\big(\partial_{x_j}^2-x_j^2\big)=-|\lambda| \mathrm{H},
$$
where $\mathrm{H} =-\Delta+|x|^2$ is the Hermite operator on $\mathbb{R}^n$. Thus the operator valued symbol $\sigma_{\mathcal{L}_{\mathbb{H}}}(\lambda)$ of $\mathcal{L}_{\mathbb{H}}$ acting on $L^2(\mathbb{R}^n)$  takes the form
$\sigma_{\mathcal{L}_{\mathbb{H}}}(\lambda)=-|\lambda| \mathrm{H}.$ {  Furthermore, for $s \in \mathbb{R}$, using  the functional calculus, the     symbol of $(-\mathcal{L}_{\mathbb{H}})^s $ is 
	$|\lambda|^s \mathrm{H}^s,$ where the notion of $\mathrm{H}^s$ is defined in (\ref{Hermite power}).
}

The Sobolev spaces $H^s(\HH), s \in \mathbb{R}$, associated to the sub-Laplacian $\mathcal{L}_{\mathbb{H}}$, are defined as
$$
H^s\left(\mathbb{H}^n\right):=\left\{f \in \mathcal{D}^{\prime}\left(\mathbb{H}^n\right):(I-\mathcal{L}_{\mathbb{H}} )^{s / 2} f \in L^2\left(\mathbb{H}^n\right)\right\},
$$
with the norm $$\|f\|_{H^s\left(\mathbb{H}^n\right)}:=\left\|(I-\mathcal{L}_{\mathbb{H}} )^{s / 2} f\right\|_{L^2\left(\mathbb{H}^n\right)}.$$ 
Similarly,  we denote by $ \dot{H}^{ s}(\mathbb{H}^n),$ the  homogeneous Sobolev  defined as the space of all $f\in \mathcal{D}'(\HH)$ such that $(-\mathcal{L}_{\mathbb{H}})^{{s}/{2}}f\in L^2(\HH)$. More generally, 	we define  $ \dot{H}^{p, s}(\mathbb{H}^n)$ as the homogeneous Sobolev space defined as the space of all  $f\in \mathcal{D}'(\HH)$ such that $(-\mathcal{L}_{\mathbb{H}})^{{s}/{2}}f\in L^p(\HH)$ for $s <\frac{Q}{p}$ and $1<p<\infty$. For a detailed definition and important properties related to the fractional sub-Laplacian on the Heisenberg group $\HH$ using the $\lambda$-twisted convolution and the spectral decomposition of the scaled Hermite operator on $\mathbb{R}^n$, we refer to \cite{roncal24,shyam}.

Now we recall the following important inequalities; see e.g. \cite{Ruz18,Fischer,GKR} for a more general graded Lie group framework. However, we will state those in the Heisenberg group setting.
\begin{theorem}[Hardy-Littlewood-Sobolev inequality]\label{eq177} Let  $\HH$ be the Heisenberg group with the homogeneous dimension $\mathcal{Q}:=2n+2$.  Let $a\geq 0$ and $1<p\leq q<\infty$ be such that 	$$	\frac{a}{\mathcal{Q}}= \frac{1}{p}-\frac{1}{q}.$$ Then we have   the following inequality $$ \|f\|_{\dot{H}^{q, -a}(\HH)} \lesssim \|f\|_{L^p(\HH)}. 	$$
\end{theorem}
%\begin{thm}[Fractional Gagliardo-Nirenbeng Inequality]
%	Let $	p, q, r \in(1, \infty)$ and  $0<a<b$. Then 
%	$$\|f\|_{\dot{H}_{\mathcal{L}}^{p, a}(\HH)} \lesssim \|f\|_{L^q(\mathbb{H}^n)}^\theta\| \|_{\dot{H}_\mathcal{L}^{r, b}(\mathbb{H}^n)}^{1-\theta},$$
%	for all $f\in L^q (\mathbb{H}^n)\cap   \dot{H}_\mathcal{L}^{r, b}(\mathbb{H}^n)$, 	where $\theta=1-\frac{a}{b}$ and $			\frac{1}{p}=\frac{\theta}{q}+\frac{1-\theta}{r}.$
%\end{thm}	 
We have the Gagliardo-Nirenberg inequality  on  $\mathbb{H}^n$ as follows:
\begin{theorem}[Gagliardo-Nirenberg inequality]\label{eq16}
	Let  $\mathcal{Q} $ be the homogeneous dimension on the Heisenberg group $\HH$.   Assume that
	$$
	s\in(0,1], 1<r<\frac{\mathcal{Q}}{s},\text { and }~  2 \leq q \leq \frac{r\mathcal{Q}}{\mathcal{Q}-sr} 		.$$
	Then  
 \begin{align}\label{Gagliardo}
\|u\|_{L^q(\mathbb{H}^n)} \lesssim\|u\|_{\dot{H}^{r,s}(\mathbb{H}^n)}^\theta\|u\|_{L^2(\mathbb{H}^n)}^{1-\theta},%\simeq\left\|\mathcal{L} u\right\|_{L^r(\mathbb{H}^n)}^\theta\|u\|_{L^p(\mathbb{H}^n)}^{1-\theta},    
\end{align}
	for $\theta=\left(\frac{1}{2}-\frac{1}{q}\right)/{\left(\frac{s}{\mathcal{Q}}+\frac{1}{2}-\frac{1}{r}\right)}\in[0,1]$, provided $\frac{s}{\mathcal{Q}}+\frac{1}{2}\neq \frac{1}{r}$.
\end{theorem}
%\begin{remark}    Note that the condition $ 1<r<\frac{Q}{s}$ in the above Gagliardo-Nirenberg inequality can be relaxed up to $r\geq \frac{Q}{s}$ also.  Indeed, for  $q\geq 2$ and $r\geq \frac{Q}{s}$, i.e., $\frac{s}{Q} -\frac{1}{r}\geq 0,$  we can see that   $$0\leq \theta=\frac{\left(\frac{1}{2}-\frac{1}{q}\right)}{\left(\frac{s}{Q}+\frac{1}{2}-\frac{1}{r}\right)} \leq \frac{ \frac{1}{2} }{\left(\frac{s}{Q}+\frac{1}{2}-\frac{1}{r}\right)} \leq \frac{ \frac{1}{2} }{  \frac{1}{2} }= 1.$$  Therefore, the modified Gagliardo-Nirenberg inequality can be seen as:    For $q\geq 2$ and $s>0$, the inequality  \eqref{Gagliardo} holds for any	  $\theta=\left(\frac{1}{2}-\frac{1}{q}\right)/{\left(\frac{s}{Q}+\frac{1}{2}-\frac{1}{r}\right)}\in[0,1]$, provided $\frac{s}{Q}+\frac{1}{2}\neq \frac{1}{r}$. \end{remark}

%	We refer to \cite{Ruz18, Fischer} for an extensive analysis of these spaces and their properties in the setting of general graded Lie groups.  We also refer to \cite{Fischer} and \cite{thanga}	for a detailed study on the Heisenberg group $\HH$.

\section{Linear damped wave equation}\label{sec3}
This section is devoted to derive the  decay estimates for the solution of the linear Cauchy problem 
 \eqref{main:homoHeisenberg}. Taking the group Fourier transform \eqref{ftheisen} to the Cauchy problem \eqref{main:homoHeisenberg} w.r.t. the space variable $\eta\in\mathbb{H}^{n}$, we obtain
     \begin{align}\label{hom:eq1}
	\begin{cases}
		\partial^2_t\widehat{u}(t,\lambda)+\sigma_{(-\mathcal{L}_{\mathbb{H}})^{\alpha}}(\lambda)\widehat{u}(t,\lambda)+b\partial_t \widehat{u}(t,\lambda)+ m\widehat{u}(t,\lambda) =0,& \lambda\in\mathbb{R}^{*},~t>0,\\ \widehat{u}(0,\lambda)=\widehat{u}_0(\lambda), &\lambda\in\mathbb{R}^{*},\\ \partial_t\widehat{u}(0,\lambda)=\widehat{u}_1(\lambda), &\lambda\in\mathbb{R}^{*},
	\end{cases} 
\end{align}
where  $\sigma_{(-\mathcal{L}_{\mathbb{H}})^{\alpha}}(\lambda)=|\lambda|^{\alpha}\mathrm{H}_{w}^{\alpha}$ is the symbol of  the fractional sub-Laplacian  $\left(-\mathcal{L}_{\mathbb{H}}\right)^{\alpha}$ on the Heisenberg group $\mathbb{H}^{n}$.  
Let us set the notation
\begin{equation*}
	\widehat{u}(t,\lambda)_{k,\ell}=(	\widehat{u}(t,\lambda)e_{k},e_{\ell})_{L^{2}(\mathbb{R}^{n})},
\end{equation*}
where $\{e_{k}\}_{k\in\mathbb{N}^{n}}$ is the family of Hermite functions. Utilizing the above notation and the relation $\mathrm{H}^{\alpha}_{w}e_{k}=\mu_{k}^{\alpha}e_{k}$, the Cauchy problem \eqref{hom:eq1} takes the form:
     \begin{align}\label{hom:eq2}
	\begin{cases}
		\partial^2_t\widehat{u}(t,\lambda)_{k,\ell}+b\partial_t \widehat{u}(t,\lambda)_{k,\ell}+|\lambda|^{\alpha}\mu_{k}^{\alpha}\widehat{u}(t,\lambda)_{k,\ell}+ m\widehat{u}(t,\lambda)_{k,\ell} =0,& \lambda\in\mathbb{R}^{*},~t>0,\\ \widehat{u}(0,\lambda)_{k,\ell}=\widehat{u}_{0}(\lambda)_{k,\ell}, &\lambda\in\mathbb{R}^{*},\\ \partial_t\widehat{u}(0,\lambda)_{k,\ell}=\widehat{u}_{1}(\lambda)_{k,\ell}, &\lambda\in\mathbb{R}^{*}.
	\end{cases} 
\end{align}
We now fix $\lambda\in\mathbb{R}^{*}$ and an ordered pair $(k,l)\in \mathbb{N}^{n}\times \mathbb{N}^{n}$, and then deal with each of the equations given by \eqref{hom:eq2} independently.
The characteristic equation of Cauchy problem \eqref{hom:eq2} is given by
\[\tau ^2+ b \tau  +|\lambda|^{\alpha}\mu_{k}^{\alpha}+m =0,\]
and the characteristic roots are given by
\begin{equation*}
	\tau_{\pm}=\begin{cases}
		-\frac{b}{2}\pm \sqrt{\frac{b^2}{4}-m-|\lambda|^{\alpha}\mu_{k}^{\alpha}}& \text{if } |\lambda|^{\alpha}\mu_{k}^{\alpha}<  \frac{b^2}{4}-m,\\
		-\frac{b}{2}&\text{if } |\lambda|^{\alpha}\mu_{k}^{\alpha}=  \frac{b^2}{4}-m,\\
		-\frac{b}{2}\pm i \sqrt{|\lambda|^{\alpha}\mu_{k}^{\alpha}-\frac{b^2}{4}+m}&\text{if } |\lambda|^{\alpha}\mu_{k}^{\alpha}>  \frac{b^2}{4}-m.
	\end{cases}
\end{equation*}
Using straightforward calculations, the solution $\widehat{u}(t,\lambda)_{k,\ell}$ of the Cauchy  problem \eqref{hom:eq2} can be written as
\begin{equation}\label{usoln}
		\widehat{u}(t,\lambda)_{k,\ell}=e^{-\frac{b}{2}t}A_0(t,\lambda;b,m,\alpha)_{k} \widehat{u}_0(\lambda)_{k,\ell}+e^{-\frac{b}{2}t}A_1(t,\lambda;b,m,\alpha)_{k} \left(\frac{b}{2} \widehat{u}_0(\lambda)_{k,\ell}+\widehat{u}_1(\lambda)_{k,\ell}\right),
\end{equation}
for all $\lambda \in \mathbb{R}^{*}$ and $t>0$, where
\begin{equation}\label{eq2}
	A_0(t,\lambda;b,m,\alpha)_{k}=\begin{cases}
		\cosh \left(\sqrt{\frac{b^2}{4} -m-|\lambda|^{\alpha}\mu_{k}^{\alpha}}~t    \right)& \text{if } |\lambda|^{\alpha}\mu_{k}^{\alpha}<  \frac{b^2}{4}-m,\\
		1&  \text{if } |\lambda|^{\alpha}\mu_{k}^{\alpha}= \frac{b^2}{4}-m,\\
		\cos \left(\sqrt{|\lambda|^{\alpha}\mu_{k}^{\alpha}-\frac{b^2}{4} +m}~t    \right)& \text{if } |\lambda|^{\alpha}\mu_{k}^{\alpha}> \frac{b^2}{4}-m,
	\end{cases} 
\end{equation}
and 
\begin{equation}\label{eq3}
	A_1(t,\lambda;b,m,\alpha)_{k}=\begin{cases}
		\frac{	\sinh \left(\sqrt{\frac{b^2}{4} -m-|\lambda|^{\alpha}\mu_{k}^{\alpha}}~t    \right)}{ \sqrt{\frac{b^2}{4} -m-|\lambda|^{\alpha}\mu_{k}^{\alpha}} }& \text{if } |\lambda|^{\alpha}\mu_{k}^{\alpha}<  \frac{b^2}{4}-m,\\
		t&  \text{if } |\lambda|^{\alpha}\mu_{k}^{\alpha}= \frac{b^2}{4}-m,\\
		\frac{	\sin \left(\sqrt{|\lambda|^{\alpha}\mu_{k}^{\alpha}-\frac{b^2}{4} +m}~t    \right)}{\sqrt{|\lambda|^{\alpha}\mu_{k}^{\alpha}-\frac{b^2}{4} +m}}& \text{if } |\lambda|^{\alpha}\mu_{k}^{\alpha}> \frac{b^2}{4}-m.
	\end{cases}  
\end{equation}
Furthermore, it is easy to verify that
\begin{eqnarray}
\partial_t A_{1}(t,\lambda;b,m,\alpha)_{k}&=&A_{0}(t,\lambda;b,m,\alpha)_{k}\label{eq1dt} \\
\partial_t A_{0}(t,\lambda;b,m,\alpha)_{k}&=&\left(\frac{b^2}{4} -m-|\lambda|^{\alpha}\mu_{k}^{\alpha}\right)A_{1}(t,\lambda;b,m,\alpha)_{k} \label{eq2dt}
\end{eqnarray}
for all $\lambda \in \mathbb{R}^{*}$ and $t>0$. Using the above relations, we deduce that
\begin{multline}\label{eq4}
		\partial_t	\widehat{u}(t,\lambda)_{k,\ell}=e^{-\frac{b}{2}t}A_{0}(t,\lambda;b,m,\alpha)_{k} \widehat{u}_1(\lambda)_{k,\ell}\\-e^{-\frac{bt}{2}}A_{1}(t,\lambda;b,m,\alpha)_{k} \left(\frac{b}{2} \widehat{u}_1(\lambda)_{k,\ell}+(|\lambda|^{\alpha}\mu_{k}^{\alpha}+m)\widehat{u}_0(\lambda)_{k,\ell}\right).
\end{multline}  
Keeping track for the explicit dependency of the decay estimates on parameters $b,m,\alpha$, and $\mathcal{Q}$ is very crucial in this study, and it will be made possible by the following lemma. We will split the frequencies $|\lambda|$ into two zones: small frequencies i.e.,  $|\lambda|<\frac{1}{\mu_{k}}\left[\frac{1}{2}\left(\frac{b^{2}}{4}-m\right)\right]^{\frac{1}{\alpha}}$ and large frequencies i.e., $|\lambda|>\frac{1}{\mu_{k}}\left[\frac{1}{2}\left(\frac{b^{2}}{4}-m\right)\right]^{\frac{1}{\alpha}}$. For ``small'' $|\lambda|$, one can use the Parseval's identity and 
 the relation \eqref{fhatnorm} in addition to Plancherel formula to get decay rates as in  \eqref{l1l2l2final}, \eqref{delul1l2final}, and  \eqref{deltl1l2l2final}. On the other hand, for  ``large" $|\lambda|$,  Plancherel formula is enough to obtain the required decay estimates,  provided the Cauchy data has adequate $L^2$-regularity. 
%In order to obtain better decay estimate for $\widehat{u}(t,\lambda)_{k,\ell}$, we will estimate  separately for small frequencies i.e.,  $|\lambda|<\frac{1}{\mu_{k}}\left[\frac{1}{2}\left(\frac{b^{2}}{4}-m\right)\right]^{\frac{1}{\alpha}}$ and large frequencies i.e., $|\lambda|>\frac{1}{\mu_{k}}\left[\frac{1}{2}\left(\frac{b^{2}}{4}-m\right)\right]^{\frac{1}{\alpha}}$.
\begin{lemma}\label{lemmaest} Let $b>0$ and $m\geq 0$ satisfies the condition $b^{2}>4m$. Let $\widehat{u}(t,\lambda)_{k,\ell}$ be a solution of the Cauchy problem \eqref{hom:eq2}. Then 
	\begin{enumerate}
	\item[$1$.] $\widehat{u}(t,\lambda)_{k,\ell}$ given by the formula \eqref{usoln} satisfies the following estimates:
	\begin{equation}\label{smallfreq1}
		|\widehat{u}(t,\lambda)_{k,\ell}|^{2}\lesssim 	e^{\left(	-b+ \sqrt{b^2-4m-4|\lambda|^{\alpha}u_{k}^{\alpha}}\right)t}\left(	|\widehat{u}_{0}(\lambda)_{k,\ell}|^{2}+	|\widehat{u}_{1}(\lambda)_{k,\ell}|^{2}\right),
	\end{equation}
	when $|\lambda|<\frac{1}{\mu_{k}}\left[\frac{1}{2}\left(\frac{b^{2}}{4}-m\right)\right]^{\frac{1}{\alpha}}$ and
	\begin{equation}\label{largefreq1}
		|\widehat{u}(t,\lambda)_{k,\ell}|^{2}\lesssim 		e^{\left(-b+\sqrt{\frac{1}{2}\left(b^{2}-4m\right)}\right) t}\left(	|\widehat{u}_{0}(\lambda)_{k,\ell}|^{2}+	|\widehat{u}_{1}(\lambda)_{k,\ell}|^{2}\right),
	\end{equation}
	when $|\lambda|>\frac{1}{\mu_{k}}\left[\frac{1}{2}\left(\frac{b^{2}}{4}-m\right)\right]^{\frac{1}{\alpha}}$.
	\item[$2$.] $\partial_{t}\widehat{u}(t,\lambda)_{k,\ell}$ given by the formula \eqref{eq4} satisfies the following estimates:
	\begin{align}\label{smallfreq2}
			|	\partial_{t}	\widehat{u}(t,\lambda)_{k,\ell}|^{2}&\lesssim e^{\left(	-b+ \sqrt{b^2-4m-4|\lambda|^{\alpha}u_{k}^{\alpha}}\right)t}(|\lambda|^{\alpha}\mu_{k}^{\alpha}+m)^{2}\left(|\widehat{u}_{0}(\lambda)_{k,\ell}|^{2}+|\widehat{u}_{1}(\lambda)_{k,\ell}|^{2}\right)\nonumber\\&+e^{\left(	-b- \sqrt{b^2-4m-4|\lambda|^{\alpha}u_{k}^{\alpha}}\right)t}|\widehat{u}_{1}(\lambda)_{k,\ell}|^{2},
	\end{align}
	when $|\lambda|<\frac{1}{\mu_{k}}\left[\frac{1}{2}\left(\frac{b^{2}}{4}-m\right)\right]^{\frac{1}{\alpha}}$ and
	\begin{equation}\label{largefreq2}
		|	\partial_{t}	\widehat{u}(t,\lambda)_{k,\ell}|^{2}\lesssim e^{\left(-b+\sqrt{\frac{1}{2}\left(b^{2}-4m\right)}\right) t}\left((|\lambda|^{\alpha}\mu_{k}^{\alpha}+m)|\widehat{u}_{0}(\lambda)_{k,\ell}|^{2}+|\widehat{u}_{1}(\lambda)_{k,\ell}|^{2}\right),
	\end{equation}
	when $|\lambda|>\frac{1}{\mu_{k}}\left[\frac{1}{2}\left(\frac{b^{2}}{4}-m\right)\right]^{\frac{1}{\alpha}}$.
		\item[$3$.] $(|\lambda|\mu_{k})^{\alpha/2}\widehat{u}(t,\lambda)_{k,\ell}$  satisfies the following estimates:
			\begin{equation}\label{smallfreq11}
			\left|(|\lambda|\mu_{k})^{\alpha/2}\widehat{u}(t,\lambda)_{k,\ell}\right|^{2}\lesssim 	e^{\left(	-b+ \sqrt{b^2-4m-4|\lambda|^{\alpha}u_{k}^{\alpha}}\right)t}(|\lambda|\mu_{k})^{\alpha}\left(	|\widehat{u}_{0}(\lambda)_{k,\ell}|^{2}+	|\widehat{u}_{1}(\lambda)_{k,\ell}|^{2}\right),
		\end{equation}
		when $|\lambda|<\frac{1}{\mu_{k}}\left[\frac{1}{2}\left(\frac{b^{2}}{4}-m\right)\right]^{\frac{1}{\alpha}}$ and
			\begin{equation}\label{largefreq22}
			\left|(|\lambda|\mu_{k})^{\alpha/2}\widehat{u}(t,\lambda)_{k,\ell}\right|^{2}\lesssim  e^{\left(-b+\sqrt{\frac{1}{2}\left(b^{2}-4m\right)}\right) t}\left(|\lambda|^{\alpha}\mu_{k}^{\alpha}|\widehat{u}_{0}(\lambda)_{k,\ell}|^{2}+|\widehat{u}_{1}(\lambda)_{k,\ell}|^{2}\right),
		\end{equation}
		when $|\lambda|>\frac{1}{\mu_{k}}\left[\frac{1}{2}\left(\frac{b^{2}}{4}-m\right)\right]^{\frac{1}{\alpha}}$.
	\end{enumerate}
\end{lemma}
\begin{remark} Combining the estimates for small frequencies and large frequencies, we can also obtain the following uniform estimates
\begin{align}
			|\widehat{u}(t,\lambda)_{k,\ell}|^{2}&\lesssim 		e^{\left(-b+\sqrt{b^{2}-m}\right) t}\left(	|\widehat{u}_{0}(\lambda)_{k,\ell}|^{2}+	|\widehat{u}_{1}(\lambda)_{k,\ell}|^{2}\right),\label{uniformestforu}\\
				|	\partial_{t}	\widehat{u}(t,\lambda)_{k,\ell}|^{2}&\lesssim e^{\left(-b+\sqrt{b^{2}-m}\right) t}\left((|\lambda|^{\alpha}\mu_{k}^{\alpha}+m)|\widehat{u}_{0}(\lambda)_{k,\ell}|^{2}+|\widehat{u}_{1}(\lambda)_{k,\ell}|^{2}\right),\label{uniformestfordeltu}\\
				\left|(|\lambda|\mu_{k})^{\alpha/2}\widehat{u}(t,\lambda)_{k,\ell}\right|^{2}&\lesssim  e^{\left(-b+\sqrt{b^{2}-4m}\right) t}\left(|\lambda|^{\alpha}\mu_{k}^{\alpha}|\widehat{u}_{0}(\lambda)_{k,\ell}|^{2}+|\widehat{u}_{1}(\lambda)_{k,\ell}|^{2}\right),\label{uniformestforlamu}
\end{align}
for all $t>0$.
\end{remark}
In order to prove the Lemma \ref{lemmaest} and decay estimates in Theorem \ref{sob:well}, we will be utilizing the following well-known fact.
\begin{remark}\label{estexp}
For any $\gamma,\delta>0$, the function $f(t)=t^{\gamma}e^{-\delta t}$ is bounded $(0,\infty)$. Furthermore,
for any $\gamma>0$ and $\beta>\delta>0$, we have 
    \begin{equation}\label{expestimate1}
        t^{\gamma}e^{-\beta t}\leq C(\gamma,\delta)e^{-(\beta-\delta)t} \quad \text{for all } t> 0,
    \end{equation}
    for some positive constant $C(\gamma,\delta)>0$. Moreover, if 
    $\beta,\beta_{1},\beta_{2},\gamma>0$ such that $\beta=\beta_{1}+\beta_{2}$, then  we have 
    \begin{equation}\label{expestimate2}
        e^{-\beta t}\leq C(\gamma,\beta_{1})t^{-\gamma}e^{-\beta_{2}t}\quad \text{for all } t> 0,
    \end{equation}
        for some positive constant $C(\gamma,\beta_{1})>0$.
\end{remark}
It is evident to verify the boundedness of $t^{\gamma}e^{-\delta t}$ on $(0,t)$ with bound depending on $\gamma$ and $\delta$. Consequently, the estimates \eqref{expestimate1} and \eqref{expestimate2} can be obtained as follows:
\begin{equation*}
    t^{\gamma}e^{-\beta t}=t^{\gamma}e^{-\delta t}e^{-(\beta-\delta)}\leq C(\gamma,\delta)e^{-(\beta-\delta)t}\quad \text{for all } t> 0,
\end{equation*}
and
\begin{equation*}
    e^{-\beta t}=t^{\gamma}e^{-\beta_{1}t}t^{-\gamma}e^{-\beta_{2}t}\leq C(\gamma,\beta_{1})t^{-\gamma}e^{-\beta_{2}t}\quad \text{for all } t> 0,
\end{equation*}
respectively.

We are now in position to supply the proof of Lemma \ref{lemmaest}.
\begin{proof}[Proof of Lemma \ref{lemmaest}] \sloppy In order to estimate $\widehat{u}(t,\lambda)_{k,\ell}$ (or $\partial_{t}\widehat{u}(t,\lambda)_{k,\ell}$), we first need to estimate $e^{-\frac{b}{2}t}|A_{0}(t,\lambda;b,m,\alpha)_{k}|$ and $e^{-\frac{b}{2}t}|A_{1}(t,\lambda;b,m,\alpha)_{k}|$. For  $|\lambda|<\frac{1}{\mu_{k}}\left[\frac{1}{2}\left(\frac{b^{2}}{4}-m\right)\right]^{\frac{1}{\alpha}}$, i.e., $|\lambda|^{\alpha}\mu_{k}^{\alpha}< \frac{1}{2}\left(\frac{b^2}{4}-m\right)$, it is easy to verify that
\begin{align}
	e^{-\frac{b}{2}t}|A_{0}(t,\lambda;b,m,\alpha)_{k}|&=e^{-\frac{b}{2}t}\cosh \left(\sqrt{\frac{b^2}{4} -m-|\lambda|^{\alpha}\mu_{k}^{\alpha}}~t    \right)\lesssim e^{\left(	-\frac{b}{2}+ \sqrt{\frac{b^2}{4}-m-|\lambda|^{\alpha}\mu_{k}^{\alpha}}\right)t},\label{bound1}\\
	e^{-\frac{b}{2}t}|A_{1}(t,\lambda;b,m,\alpha)_{k}|&=e^{-\frac{b}{2}t}\frac{	\sinh \left(\sqrt{\frac{b^2}{4} -m-|\lambda|^{\alpha}\mu_{k}^{\alpha}}~t    \right)}{ \sqrt{\frac{b^2}{4} -m-|\lambda|^{\alpha}\mu_{k}^{\alpha}} }\lesssim e^{\left(	-\frac{b}{2}+ \sqrt{\frac{b^2}{4}-m-|\lambda|^{\alpha}\mu_{k}^{\alpha}}\right)t}.\label{bound2}
\end{align}
Combining the above bounds with the formula \eqref{usoln}, we can obtain our required estimate \eqref{smallfreq1}.

For  $|\lambda|>\frac{1}{\mu_{k}}\left[\frac{1}{2}\left(\frac{b^{2}}{4}-m\right)\right]^{\frac{1}{\alpha}}$, i.e., $|\lambda|^{\alpha}\mu_{k}^{\alpha}> \frac{1}{2}\left(\frac{b^2}{4}-m\right)$, we have
\begin{align}
	e^{-\frac{b}{2}t}|A_{0}(t,\lambda;b,m,\alpha)_{k}|&\lesssim
\begin{cases}
	e^{-\frac{b}{2}t}\cosh \left(\sqrt{\frac{1}{2}\left(\frac{b^2}{4} -m\right)}~t    \right)&\text{if }\frac{1}{2}\left(\frac{b^2}{4}-m\right)<|\lambda|^{\alpha}\mu_{k}^{\alpha}< \frac{b^2}{4}-m,\\
	e^{-\frac{b}{2}t}&\text{if }|\lambda|^{\alpha}\mu_{k}^{\alpha}= \frac{b^2}{4}-m,\\
	e^{-\frac{b}{2}t}&\text{if }|\lambda|^{\alpha}\mu_{k}^{\alpha}> \frac{b^2}{4}-m,
\end{cases}\label{bound3a}\\
&\lesssim e^{\left(	-\frac{b}{2}+ \sqrt{\frac{1}{2}\left(\frac{b^2}{4}-m\right)}\right)t},\label{bound3}
\end{align}
 and 
\begin{align}
	e^{-\frac{b}{2}t}|A_{1}(t,\lambda;b,m,\alpha)_{k}|&\lesssim
	\begin{cases}
		e^{-\frac{b}{2}t}\frac{	\sinh \left(\sqrt{\frac{b^2}{4} -m-|\lambda|^{\alpha}\mu_{k}^{\alpha}}~t    \right)}{ \sqrt{\frac{b^2}{4} -m-|\lambda|^{\alpha}\mu_{k}^{\alpha}} }&\text{if }\frac{1}{2}\left(\frac{b^2}{4}-m\right)<|\lambda|^{\alpha}\mu_{k}^{\alpha}< \frac{b^2}{4}-m,\\
		te^{-\frac{b}{2}t}&\text{if }|\lambda|^{\alpha}\mu_{k}^{\alpha}= \frac{b^2}{4}-m,\\
		te^{-\frac{b}{2}t}&\text{if }|\lambda|^{\alpha}\mu_{k}^{\alpha}> \frac{b^2}{4}-m,
	\end{cases}\label{bound4a}\\
	&\lesssim e^{\left(	-\frac{b}{2}+ \sqrt{\frac{1}{2}\left(\frac{b^2}{4}-m\right)}\right)t}\label{bound4},
\end{align}
 where in the last steps we utilize the Remark \ref{estexp} for $0<\sqrt{\frac{1}{2}\left(\frac{b^{2}}{4}-m\right)}<\frac{b}{2}$. 
Now using the the bounds \eqref{bound3} and \eqref{bound4} in formula \eqref{usoln}, we can obtain our required estimate \eqref{largefreq1}.

Let us now estimate the derivative term $\partial_{t}\widehat{u}(t,\lambda)_{k,\ell}$.
For  $|\lambda|<\frac{1}{\mu_{k}}\left[\frac{1}{2}\left(\frac{b^{2}}{4}-m\right)\right]^{\frac{1}{\alpha}}$, the expression \eqref{eq4} gives
\begin{align}
	&\partial_{t}	\widehat{u}(t,\lambda)_{k,\ell}=-e^{-\frac{b}{2}t}	\frac{	\sinh \left(\sqrt{\frac{b^2}{4} -m-|\lambda|^{\alpha}\mu_{k}^{\alpha}}~t    \right)}{\sqrt{\frac{b^2}{4} -m-|\lambda|^{\alpha}\mu_{k}^{\alpha}}}(|\lambda|^{\alpha}\mu_{k}^{\alpha}+m)\widehat{u}_0(\lambda)_{k,\ell}\nonumber\\
	&+e^{-\frac{b}{2}t}\left(\cosh \left(\sqrt{\frac{b^2}{4} -m-|\lambda|^{\alpha}\mu_{k}^{\alpha}}~t    \right)-\frac{b}{2}\frac{	\sinh \left(\sqrt{\frac{b^2}{4} -m-|\lambda|^{\alpha}\mu_{k}^{\alpha}}~t    \right)}{\sqrt{\frac{b^2}{4} -m-|\lambda|^{\alpha}\mu_{k}^{\alpha}}}\right)\widehat{u}_{1}(\lambda)_{k,\ell}\label{deltsolution1}\\
	&=-e^{-\frac{b}{2}t}	\frac{	\sinh \left(\sqrt{\frac{b^2}{4} -m-|\lambda|^{\alpha}\mu_{k}^{\alpha}}~t    \right)}{\sqrt{\frac{b^2}{4} -m-|\lambda|^{\alpha}\mu_{k}^{\alpha}}}(|\lambda|^{\alpha}\mu_{k}^{\alpha}+m)\widehat{u}_0(\lambda)_{k,\ell}\nonumber\\
	&+e^{-\frac{b}{2}t}\left[\left(\frac{1}{2}-\frac{b}{4\sqrt{\frac{b^2}{4} -m-|\lambda|^{\alpha}\mu_{k}^{\alpha}}}\right)e^{\sqrt{\frac{b^2}{4} -m-|\lambda|^{\alpha}\mu_{k}^{\alpha}}}\right.\nonumber\\
	&\left.+\left(\frac{1}{2}+\frac{b}{4\sqrt{\frac{b^2}{4} -m-|\lambda|^{\alpha}\mu_{k}^{\alpha}}}\right)e^{-\sqrt{\frac{b^2}{4} -m-|\lambda|^{\alpha}\mu_{k}^{\alpha}}}\right]\widehat{u}_{1}(\lambda)_{k,\ell}\label{deltsolution}.
\end{align}
Using the identity
\begin{equation}\label{relat}
	-4x\leq -1+\sqrt{1-4x}\leq -2x\quad \text{for any } x\in \left[0,\frac{1}{4}\right],
\end{equation}
the quantity in last two terms can be estimated as follows:
\begin{equation}\label{bound6}
\left|\frac{1}{2}-\frac{b}{4\sqrt{\frac{b^2}{4} -m-|\lambda|^{\alpha}\mu_{k}^{\alpha}}}\right|\lesssim|\lambda|^{\alpha}\mu_{k}^{\alpha}+m\quad \text{and }\quad
\left|\frac{1}{2}+\frac{b}{4\sqrt{\frac{b^2}{4} -m-|\lambda|^{\alpha}\mu_{k}^{\alpha}}}\right|\lesssim 1.
\end{equation}
Combining the expression \eqref{deltsolution} with the estimates \eqref{bound2} for the first term and \eqref{bound6} for last two terms, we deduce that
\begin{multline}
	|	\partial_{t}	\widehat{u}(t,\lambda)_{k,\ell}|\lesssim e^{\left(-\frac{b}{2}+\sqrt{\frac{b^2}{4} -m-|\lambda|^{\alpha}\mu_{k}^{\alpha}}\right)t}(|\lambda|^{\alpha}\mu_{k}^{\alpha}+m)\left(|\widehat{u}_{0}(\lambda)_{k,\ell}|+|\widehat{u}_{1}(\lambda)_{k,\ell}|\right)\\+e^{\left(-\frac{b}{2}-\sqrt{\frac{b^2}{4} -m-|\lambda|^{\alpha}\mu_{k}^{\alpha}}\right)t}|\widehat{u}_{1}(\lambda)_{k,\ell}|,
\end{multline}
for $|\lambda|<\frac{1}{\mu_{k}}\left[\frac{1}{2}\left(\frac{b^{2}}{4}-m\right)\right]^{\frac{1}{\alpha}}$. Thus, this instantly provide our required estimate \eqref{smallfreq2}.

\medskip

For $\frac{1}{\mu_{k}}\left[\frac{1}{2}\left(\frac{b^{2}}{4}-m\right)\right]^{\frac{1}{\alpha}}<|\lambda|<\frac{1}{\mu_{k}}\left[\frac{b^{2}}{4}-m\right]^{\frac{1}{\alpha}}$, utilizing the bounds \eqref{bound3} and \eqref{bound4} for the representation \eqref{eq4} for $\partial_{t}\widehat{u}(t,\lambda)_{k,\ell}$, we deduce that
\begin{align}\label{deltfirst}
	|	\partial_{t}	\widehat{u}(t,\lambda)_{k,\ell}&|\lesssim e^{\left(	-\frac{b}{2}+ \sqrt{\frac{1}{2}\left(\frac{b^2}{4}-m\right)}\right)t}\left(|\widehat{u}_{0}(\lambda)_{k,\ell}|+|\widehat{u}_{1}(\lambda)_{k,\ell}|\right)\nonumber\\
	&\cong e^{\left(	-\frac{b}{2}+ \sqrt{\frac{1}{2}\left(\frac{b^2}{4}-m\right)}\right)t}\left((|\lambda|^{\alpha}\mu_{k}^{\alpha}+m)|\widehat{u}_{0}(\lambda)_{k,\ell}|+|\widehat{u}_{1}(\lambda)_{k,\ell}|\right)\nonumber\\
	&\cong e^{\left(	-\frac{b}{2}+ \sqrt{\frac{1}{2}\left(\frac{b^2}{4}-m\right)}\right)t}\left(\sqrt{|\lambda|^{\alpha}\mu_{k}^{\alpha}+m}|\widehat{u}_{0}(\lambda)_{k,\ell}|+|\widehat{u}_{1}(\lambda)_{k,\ell}|\right),
\end{align}
where we use the boundedness of $|\lambda|$ to obtain the equivalent representations. 

 For $|\lambda|>\frac{1}{\mu_{k}}\left[\frac{b^{2}}{4}-m\right]^{\frac{1}{\alpha}}$, we have
\begin{multline}\label{deltrep}
	\partial_{t}	\widehat{u}(t,\lambda)_{k,\ell}
	=-e^{-\frac{b}{2}t}	\frac{	\sin \left(\sqrt{|\lambda|^{\alpha}\mu_{k}^{\alpha}-\frac{b^2}{4}+m}~t    \right)}{\sqrt{|\lambda|^{\alpha}\mu_{k}^{\alpha}-\frac{b^2}{4}+m}}(|\lambda|^{\alpha}\mu_{k}^{\alpha}+m)\widehat{u}_0(\lambda)_{k,\ell}\\
	+e^{-\frac{b}{2}t}\left(\cos \left(\sqrt{|\lambda|^{\alpha}\mu_{k}^{\alpha}-\frac{b^2}{4}+m}~t    \right)-\frac{b}{2}\frac{	\sin \left(\sqrt{|\lambda|^{\alpha}\mu_{k}^{\alpha}-\frac{b^2}{4}+m}~t    \right)}{\sqrt{|\lambda|^{\alpha}\mu_{k}^{\alpha}-\frac{b^2}{4}+m}}\right)\widehat{u}_{1}(\lambda)_{k,\ell}.
\end{multline}
Now for $\frac{1}{\mu_{k}}\left[\frac{b^{2}}{4}-m\right]^{\frac{1}{\alpha}}<|\lambda|<\frac{1}{\mu_{k}}\left[2\left(\frac{b^{2}}{4}-m\right)\right]^{\frac{1}{\alpha}}$, if we utilize the bounds \eqref{bound3a} and \eqref{bound4a} for the above representaion, we can again immediately obtain the estimate \eqref{deltfirst}. 

Furthermore,  it is easy to verify that
\begin{equation*}
    \sqrt{\frac{|\lambda|^{\alpha}\mu_{k}^{\alpha}+m}{|\lambda|^{\alpha}\mu_{k}^{\alpha}-\frac{b^{2}}{4}+m}}\quad \text{is bounded for }|\lambda|> \frac{1}{\mu_{k}}\left[2\left(\frac{b^{2}}{4}-m\right)\right]^{\frac{1}{\alpha}}.
\end{equation*}
Combining the above fact for the first term of the equation \eqref{deltrep} and using the bounds \eqref{bound3a} and \eqref{bound4a} for the last two terms, we can obtain
\begin{align}\label{delt2}
	|	\partial_{t}	\widehat{u}(t,\lambda)_{k,\ell}|
&\lesssim e^{-\frac{b}{2}t}\left(\sqrt{|\lambda|^{\alpha}\mu_{k}^{\alpha}+m}|\widehat{u}_{0}(\lambda)_{k,\ell}|+|\widehat{u}_{1}(\lambda)_{k,\ell}|\right)\nonumber\\
&\lesssim e^{\left(	-\frac{b}{2}+ \sqrt{\frac{1}{2}\left(\frac{b^2}{4}-m\right)}\right)t}\left(\sqrt{|\lambda|^{\alpha}\mu_{k}^{\alpha}+m}|\widehat{u}_{0}(\lambda)_{k,\ell}|+|\widehat{u}_{1}(\lambda)_{k,\ell}|\right).
\end{align}
 Thus, we obtain our required estimate \eqref{largefreq2} for $|\lambda|>\frac{1}{\mu_{k}}\left[\frac{1}{2}\left(\frac{b^{2}}{4}-m\right)\right]^{\frac{1}{\alpha}}$.

Similarly, using the arguments used for $\partial_{t}	\widehat{u}(t,\lambda)_{k,\ell}$, we can also obtain the estimates \eqref{smallfreq11} and \eqref{largefreq22}. This completes the proof of the lemma.
\end{proof}

We are now in a position to establish the $L^{2}-L^{2}$ and $(L^{1}\cap L^{2})-L^{2}$ decay estimates for the solution of the linear Cauchy problem \eqref{main:homoHeisenberg}. We will start by estimating the $(L^{1}\cap L^{2})-L^{2}$ decay and at suitable stages, we will additionally extract the $L^{2}-L^{2}$ decay estimates. 
\subsection{Estimates for  $\|u(t,\cdot)\|_{L^{2}(\mathbb{H}^{n})}$}\label{subsec3.1}
Here, we estimate the $L^{2}$-norm of $u(t,\eta),$ i.e.,  $\|u(t,\cdot)\|_{L^{2}(\mathbb{H}^{n})}$. Using the Plancherel formula \eqref{Plancherel}, we have
\begin{equation}
	\begin{aligned}\label{udnormest}
		\|u(t,&\cdot)\|^{2}_{L^{2}(\mathbb{H}^{n})}
		=c_{n}\sum_{k,\ell\in\mathbb{N}^{n}}\int_{\mathbb{R}{*}}(\widehat{u}(t,\lambda)e_{k},e_{\ell})^{2}_{L^{2}(\mathbb{R}^{n})}|\lambda|^{n}\mathrm{d}\lambda\nonumber\\
		&=c_{n}\sum_{k,\ell\in\mathbb{N}^{n}}\left(\int_{0<|\lambda|<\frac{1}{\mu_{k}}\left[\frac{1}{2}\left(\frac{b^{2}}{4}-m\right)\right]^{\frac{1}{\alpha}}}+\int_{|\lambda|>\frac{1}{\mu_{k}}\left[\frac{1}{2}\left(\frac{b^{2}}{4}-m\right)\right]^{\frac{1}{\alpha}}}\right)(\widehat{u}(t,\lambda)e_{k},e_{\ell})^{2}_{L^{2}(\mathbb{R}^{n})}|\lambda|^{n}\mathrm{d}\lambda\nonumber\\
		&=\mathrm{I}_{\operatorname{low}}+\mathrm{I}_{\operatorname{high}}.
	\end{aligned}
\end{equation}
Let us now estimate each term $\mathrm{I}_{\operatorname{low}}$ and $\mathrm{I}_{\operatorname{high}}$ one by one, beginning with $\mathrm{I}_{\operatorname{low}}$: 
\begin{equation}\label{lowfirst}
	\mathrm{I}_{\operatorname{low}}=c_{n}\sum_{k,\ell\in\mathbb{N}^{n}}\int_{0<|\lambda|<\frac{1}{\mu_{k}}\left[\frac{1}{2}\left(\frac{b^{2}}{4}-m\right)\right]^{\frac{1}{\alpha}}}(\widehat{u}(t,\lambda)e_{k},e_{\ell})^{2}_{L^{2}(\mathbb{R}^{n})}|\lambda|^{n}\mathrm{d}\lambda.
\end{equation}
Using the identity \eqref{relat}, the estimate \eqref{smallfreq1} for solution $	\widehat{u}(t,\lambda)_{k,\ell}$ can be further estimated as 
\begin{eqnarray}\label{utsolest}
	|\widehat{u}(t,\lambda)_{k,\ell}|^{2}&\lesssim& e^{\left(	-b+ \sqrt{b^2-4m-4|\lambda|^{\alpha}u_{k}^{\alpha}}\right)t}\left(	|\widehat{u}_{0}(\lambda)_{k,\ell}|^{2}+	|\widehat{u}_{1}(\lambda)_{k,\ell}|^{2}\right)\nonumber\\
	&=& e^{b\left(	-1+ \sqrt{1-\frac{4}{b^{2}}\left(m+|\lambda|^{\alpha}u_{k}^{\alpha}\right)}\right)t}\left(	|\widehat{u}_{0}(\lambda)_{k,\ell}|^{2}+	|\widehat{u}_{1}(\lambda)_{k,\ell}|^{2}\right)\nonumber\\
	&\lesssim& e^{-\frac{2}{b}\left(m+|\lambda|^{\alpha}u_{k}^{\alpha}\right)t}\left(	|\widehat{u}_{0}(\lambda)_{k,\ell}|^{2}+	|\widehat{u}_{1}(\lambda)_{k,\ell}|^{2}\right).
\end{eqnarray}
Now using the above estimate and the Parseval's identity \eqref{parseval} in \eqref{lowfirst}, we deduce that
\begin{align}\label{ilowest}
	&\mathrm{I}_{\mathrm{low}} \lesssim\sum_{k\in\mathbb{N}^{n}}\int_{0<|\lambda|<\frac{1}{\mu_{k}}\left[\frac{1}{2}\left(\frac{b^{2}}{4}-m\right)\right]^{\frac{1}{\alpha}}}e^{-\frac{2}{b}\left(m+|\lambda|^{\alpha}u_{k}^{\alpha}\right)t}\sum_{\ell\in\mathbb{N}^{n}}	\left(	|\widehat{u}_{0}(\lambda)_{k,\ell}|^{2}+	|\widehat{u}_{1}(\lambda)_{k,\ell}|^{2}\right)|\lambda|^{n}\mathrm{d}\lambda\nonumber\\
	&=\sum_{k\in\mathbb{N}^{n}}\int_{0<|\lambda|<\frac{1}{\mu_{k}}\left[\frac{1}{2}\left(\frac{b^{2}}{4}-m\right)\right]^{\frac{1}{\alpha}}}e^{-\frac{2}{b}\left(m+|\lambda|^{\alpha}u_{k}^{\alpha}\right)t}\left(	\|\widehat{u}_{0}(\lambda)e_{k}\|_{L^{2}(\mathbb{R}^{n})}^{2}+	\|\widehat{u}_{1}(\lambda)e_{k}\|_{L^{2}(\mathbb{R}^{n})}^{2}\right)|\lambda|^{n}\mathrm{d}\lambda\nonumber\\
	&\lesssim e^{-\frac{2m}{b}t}\sum_{k\in\mathbb{N}^{n}}\int_{0<|\lambda|<\frac{1}{\mu_{k}}\left[\frac{1}{2}\left(\frac{b^{2}}{4}-m\right)\right]^{\frac{1}{\alpha}}}e^{-\frac{2}{b}|\lambda|^{\alpha}u_{k}^{\alpha}t}|\lambda|^{n}\mathrm{d}\lambda\left(\|u_{0}\|^{2}_{L^{1}(\mathbb{H}^{n})}+\|u_{1}\|^{2}_{L^{1}(\mathbb{H}^{n})}\right)\nonumber\\
	&\lesssim e^{-\frac{2m}{b}t} \sum_{k\in\mathbb{N}^{n}}\int_{0}^{\frac{1}{\mu_{k}}\left[\frac{1}{2}\left(\frac{b^{2}}{4}-m\right)\right]^{\frac{1}{\alpha}}}e^{-\frac{2}{b}\lambda^{\alpha}u_{k}^{\alpha}t}\lambda^{n}\mathrm{d}\lambda \left(\|u_{0}\|^{2}_{L^{1}(\mathbb{H}^{n})}+\|u_{1}\|^{2}_{L^{1}(\mathbb{H}^{n})}\right),
\end{align}
where 2nd last inequality follows from the relation \eqref{fhatnorm} and the fact that
$\|e_{k}\|_{L^{2}(\mathbb{R}^{n})}=1$. Now using the substitution $\theta=\frac{2}{b}\lambda^{\alpha}u_{k}^{\alpha}t$ in \eqref{ilowest}, we have
\begin{equation*}
	\lambda^{n}\mathrm{d}\lambda =\frac{1}{\alpha}\left(\frac{b}{2}\right)^{\frac{n+1}{\alpha}}t^{-\left(\frac{n+1}{\alpha}\right)}\mu_{k}^{-(n+1)}\theta^{\frac{n+1-\alpha}{\alpha}}\mathrm{d}\theta \quad \text{and }\quad \theta \text{ varies from } 0\text{ to } \frac{1}{b}\left(\frac{b^{2}}{4}-m\right)t.
\end{equation*}
Consequently, the above estimate takes the form
\begin{align}\label{ilowestimate}
	&\mathrm{I}_{\mathrm{low}} \lesssim t^{-\left(\frac{n+1}{\alpha}\right)}e^{-\frac{2m}{b}t}  \sum_{k\in\mathbb{N}^{n}}\mu_{k}^{-(n+1)}\int_{0}^{\frac{1}{b}\left(\frac{b^{2}}{4}-m\right)t}e^{-\theta}\theta^{\frac{n+1-\alpha}{\alpha}}\mathrm{d}\theta\left(\|u_{0}\|^{2}_{L^{1}(\mathbb{H}^{n})}+\|u_{1}\|^{2}_{L^{1}(\mathbb{H}^{n})}\right) \nonumber\\
	&\lesssim t^{-\left(\frac{n+1}{\alpha}\right)}e^{-\frac{2m}{b}t}\sum_{k\in\mathbb{N}^{n}}(2|k|+n)^{-(n+1)}  \int_{0}^{\infty}e^{-\theta}\theta^{\frac{n+1-\alpha}{\alpha}}\mathrm{d}\theta\left(\|u_{0}\|^{2}_{L^{1}(\mathbb{H}^{n})}+\|u_{1}\|^{2}_{L^{1}(\mathbb{H}^{n})}\right) \nonumber\\
	&=t^{-\left(\frac{n+1}{\alpha}\right)}e^{-\frac{2m}{b}t}\sum_{k\in\mathbb{N}^{n}}(2|k|+n)^{-(n+1)}\Gamma\left(\frac{n+1}{\alpha}\right)\left(\|u_{0}\|^{2}_{L^{1}(\mathbb{H}^{n})}+\|u_{1}\|^{2}_{L^{1}(\mathbb{H}^{n})}\right)\nonumber\\
	&\lesssim t^{-\frac{\mathcal{Q}}{2\alpha}}e^{-\frac{2m}{b}t}\left(\|u_{0}\|^{2}_{L^{1}(\mathbb{H}^{n})}+\|u_{1}\|^{2}_{L^{1}(\mathbb{H}^{n})}\right)\quad (\because \mathcal{Q}=2n+2),
\end{align}
\sloppy where the last inequality follows from the convergence of the infinite series $\sum_{k\in\mathbb{N}^{n}}(2|k|+n)^{-(n+1)}$. Now we will estimate the term $\mathrm{I}_{\mathrm{high}}$ given by
\begin{equation*}
\mathrm{I}_{\mathrm{high}}=c_{n}\sum_{k,\ell\in\mathbb{N}^{n}} \int_{|\lambda|>\frac{1}{\mu_{k}}\left[\frac{1}{2}\left(\frac{b^{2}}{4}-m\right)\right]^{\frac{1}{\alpha}}} (\widehat{u}(t,\lambda)e_{k},e_{\ell})^{2}_{L^{2}(\mathbb{R}^{n})}|\lambda|^{n}\mathrm{d}\lambda.
\end{equation*}
Recalling the estimate \eqref{largefreq1} for $\lambda$ satisfying $|\lambda|>\frac{1}{\mu_{k}}\left[\frac{1}{2}\left(\frac{b^{2}}{4}-m\right)\right]^{\frac{1}{\alpha}}$, the above term $\mathrm{I}_{\mathrm{high}}$ can be estimated as follows: 
\begin{align}\label{ihighestimate}
	\mathrm{I}_{\mathrm{high}}
	&\lesssim e^{\left(-b+\sqrt{\frac{1}{2}\left(b^{2}-4m\right)}\right) t}\sum_{k,\ell\in\mathbb{N}^{n}} \int_{|\lambda|>\frac{1}{\mu_{k}}\left[\frac{1}{2}\left(\frac{b^{2}}{4}-m\right)\right]^{\frac{1}{\alpha}}} \left(	|\widehat{u}_{0}(\lambda)_{k,\ell}|^{2}+	|\widehat{u}_{1}(\lambda)_{k,\ell}|^{2}\right)|\lambda|^{n}\mathrm{d}\lambda\nonumber\\
	&\lesssim e^{\left(-b+\sqrt{\frac{1}{2}\left(b^{2}-4m\right)}\right) t}\int_{\mathbb{R}^{*}}\sum_{k,\ell\in\mathbb{N}^{n}}  \left(	|\widehat{u}_{0}(\lambda)_{k,\ell}|^{2}+	|\widehat{u}_{1}(\lambda)_{k,\ell}|^{2}\right)|\lambda|^{n}\mathrm{d}\lambda\nonumber\\
	&= e^{\left(-b+\sqrt{\frac{1}{2}\left(b^{2}-4m\right)}\right) t}  \left(	\|u_{0}\|^{2}_{L^{2}(\mathbb{H}^{n})}+	\|u_{1}\|^{2}_{L^{2}(\mathbb{H}^{n})}\right).
\end{align}
As the rates of decay for $\mathrm{I}_{\mathrm{low}}$ and $\mathrm{I}_{\mathrm{high}}$ differ we must first identify a uniform decay before combining the two estimates. Note that the factor $-b+\sqrt{\frac{1}{2}\left(b^{2}-4m\right)}$ in \eqref{ihighestimate} can be expressed as
\begin{equation}\label{equalityrel}
	-b+\sqrt{\frac{1}{2}\left(b^{2}-4m\right)}=-\frac{m}{b}-f_{1}(b,m),\quad b^{2}> 4m,
\end{equation}
where
\begin{eqnarray}
	f_{1}(b,m)&=&-\frac{m}{b}+b-\sqrt{\frac{1}{2}\left(b^{2}-4m\right)}\nonumber\\
	&=&-\frac{m}{b}+\frac{\sqrt{2}-1}{\sqrt{2}}b+\frac{b}{\sqrt{2}}\left(1-\sqrt{1-4\frac{m}{b^{2}}}\right)\nonumber\\
	&\geq&-\frac{m}{b}+\frac{\sqrt{2}-1}{\sqrt{2}}b+\frac{b}{\sqrt{2}}\left(2\frac{m}{b^{2}}\right)\nonumber\\
	&=&\left(\sqrt{2}-1\right)\left(\frac{m}{b}+\frac{b}{\sqrt{2}}\right)>0,
\end{eqnarray}
where the last inequality follows from the relation \eqref{relat}. Now utilizing the Remark \ref{estexp} for the equality \eqref{equalityrel}, the exponential factor in estimate \eqref{ihighestimate} can be further estimated as
\begin{equation}\label{expest}
	e^{\left(-b+\sqrt{\frac{1}{2}\left(b^{2}-4m\right)}\right) t}=e^{-\frac{m}{b}t}e^{-f_{1}(b,m)t}\leq C_{1}(b,m,\mathcal{Q},\alpha)t^{-\frac{\mathcal{Q}}{2\alpha}}e^{-\frac{m}{b}t}.
\end{equation}
Combining the estimates \eqref{ilowestimate} and \eqref{ihighestimate} with \eqref{expest}, we get
\begin{equation}\label{decayest1}
	\|u(t, \cdot)\|^{2}_{L^2\left(\mathbb{H}^{n}\right)}\lesssim t^{-\frac{\mathcal{Q}}{2\alpha}}e^{-\frac{m}{b}t}\left(\left\|u_{0}\right\|^{2}_{L^{2}\left(\mathbb{H}^{n}\right)}+\left\|u_{1}\right\|^{2}_{L^{2}\left(\mathbb{H}^{n}\right)}+\left\|u_{0}\right\|^{2}_{L^{1}\left(\mathbb{H}^{n}\right)}+\left\|u_{1}\right\|^{2}_{L^{1}\left(\mathbb{H}^{n}\right)}\right).
\end{equation}
It is important to take into account that in this case, using the uniform estimate \eqref{uniformestforu} rather than \eqref{smallfreq1} and \eqref{largefreq1} allows us to deduce that
\begin{equation}\label{l2l2lest}
	\|u(t, \cdot)\|^{2}_{L^2\left(\mathbb{H}^{n}\right)}\lesssim e^{\left(-b+\sqrt{b^{2}-4m}\right) t}\left(\left\|u_{0}\right\|^{2}_{L^{2}\left(\mathbb{H}^{n}\right)}+\left\|u_{1}\right\|^{2}_{L^{2}\left(\mathbb{H}^{n}\right)}\right).
\end{equation} 
Note that the above estimate gives exponential decay for $m>0$ and no decay for $m=0$ which allows us to exclude the singular behavior of $	\|u(t, \cdot)\|^{2}_{L^2\left(\mathbb{H}^{n}\right)}$  in \eqref{decayest1} as $t\to 0^{+}$. Thus, the estimate \eqref{decayest1} can be replaced by
\begin{equation}\label{l1l2l2est}
	\|u(t, \cdot)\|^{2}_{L^2\left(\mathbb{H}^{n}\right)}\lesssim (1+t)^{-\frac{\mathcal{Q}}{2\alpha}}e^{-\frac{m}{b}t}\left(\left\|u_{0}\right\|^{2}_{L^{2}\left(\mathbb{H}^{n}\right)}+\left\|u_{1}\right\|^{2}_{L^{2}\left(\mathbb{H}^{n}\right)}+\left\|u_{0}\right\|^{2}_{L^{1}\left(\mathbb{H}^{n}\right)}+\left\|u_{1}\right\|^{2}_{L^{1}\left(\mathbb{H}^{n}\right)}\right),
\end{equation}
for all $t>0$ and hence the estimates \eqref{l2l2lest} and \eqref{l1l2l2est} will immediately give our required estimate \eqref{l2l2final} and \eqref{l1l2l2final}.
\subsection{Estimates for  $\left\|(-\mathcal{L}_{\mathbb{H}})^{\alpha/2}u(t,\cdot)\right\|_{L^{2}(\mathbb{H}^{n})}$}\label{subsec3.2}
 Now, we estimate the $L^{2}$-norm of $(-\mathcal{L}_{\mathbb{H}})^{\alpha/2}u(t,\eta),$ i.e.,  $\left\|(-\mathcal{L}_{\mathbb{H}})^{\alpha/2}u(t,\cdot)\right\|_{L^{2}(\mathbb{H}^{n})}$. Using the Plancherel formula \eqref{Plancherel}, we have 
\begin{align}\nonumber
	&\left\|(-\mathcal{L}_{\mathbb{H}})^{\alpha/2}u(t,\cdot)\right\|^{2}_{L^{2}(\mathbb{H}^{n})}\\\nonumber
	&=c_{n}\sum_{k,\ell\in\mathbb{N}^{n}}\left(\int_{0<|\lambda|<\frac{1}{\mu_{k}}\left[\frac{1}{2}\left(\frac{b^{2}}{4}-m\right)\right]^{\frac{1}{\alpha}}}+\int_{|\lambda|>\frac{1}{\mu_{k}}\left[\frac{1}{2}\left(\frac{b^{2}}{4}-m\right)\right]^{\frac{1}{\alpha}}}\right)(|\lambda|^{\alpha/2}\mu_{k}^{\alpha/2}\widehat{u}(t,\lambda)e_{k},e_{\ell})^{2}_{L^{2}(\mathbb{R}^{n})}|\lambda|^{n}\mathrm{d}\lambda\\
	&=\operatorname{J}_{\operatorname{low}}+\mathrm{J}_{\operatorname{high}}.
\end{align}
 Let us estimate $\operatorname{J}_{\operatorname{low}}$ first given by $$\mathrm{J}_{\mathrm{low}}=c_{n}\sum_{k,\ell\in\mathbb{N}^{n}}\int_{|\lambda|<\frac{1}{\mu_{k}}\left[\frac{1}{2}\left(\frac{b^{2}}{4}-m\right)\right]^{\frac{1}{\alpha}}}(|\lambda|^{\alpha/2}\mu_{k}^{\alpha/2}\widehat{u}(t,\lambda)e_{k},e_{\ell})^{2}_{L^{2}(\mathbb{R}^{n})}|\lambda|^{n}\mathrm{d}\lambda.$$ 
Using the estimate \eqref{utsolest} and the relation \eqref{fhatnorm}, analogous to the estimate \eqref{ilowest}, we can obtain
\begin{align}
	&\mathrm{J}_{\mathrm{low}} \lesssim\sum_{k\in\mathbb{N}^{n}}\mu_{k}^{\alpha}\int_{0<|\lambda|<\frac{1}{\mu_{k}}\left[\frac{1}{2}\left(\frac{b^{2}}{4}-m\right)\right]^{\frac{1}{\alpha}}}e^{-\frac{2}{b}\left(m+|\lambda|^{\alpha}u_{k}^{\alpha}\right)t}\sum_{\ell\in\mathbb{N}^{n}}\left(	|\widehat{u}_{0}(\lambda)_{k,\ell}|^{2}+	|\widehat{u}_{1}(\lambda)_{k,\ell}|^{2}\right)|\lambda|^{n+\alpha}\mathrm{d}\lambda\nonumber\\
	&=\sum_{k\in\mathbb{N}^{n}}\mu_{k}^{\alpha}\int_{0<|\lambda|<\frac{1}{\mu_{k}}\left[\frac{1}{2}\left(\frac{b^{2}}{4}-m\right)\right]^{\frac{1}{\alpha}}}e^{-\frac{2}{b}\left(m+|\lambda|^{\alpha}u_{k}^{\alpha}\right)t}\left(	\|\widehat{u}_{0}(\lambda)e_{k}\|_{L^{2}(\mathbb{R}^{n})}^{2}+	\|\widehat{u}_{1}(\lambda)e_{k}\|_{L^{2}(\mathbb{R}^{n})}^{2}\right)|\lambda|^{n+\alpha}\mathrm{d}\lambda\nonumber\\
	&\lesssim e^{-\frac{2m}{b}t} \sum_{k\in\mathbb{N}^{n}}\mu_{k}^{\alpha}\int_{0<|\lambda|<\frac{1}{\mu_{k}}\left[\frac{1}{2}\left(\frac{b^{2}}{4}-m\right)\right]^{\frac{1}{\alpha}}}e^{-\frac{2}{b}|\lambda|^{\alpha}u_{k}^{\alpha}t}|\lambda|^{n+\alpha}\mathrm{d}\lambda \left(\|u_{0}\|^{2}_{L^{1}(\mathbb{H}^{n})}+\|u_{1}\|^{2}_{L^{1}(\mathbb{H}^{n})}\right) \nonumber\\
	&=e^{-\frac{2m}{b}t}\sum_{k\in\mathbb{N}^{n}}\mu_{k}^{\alpha}\int_{0}^{\frac{1}{\mu_{k}}\left[\frac{1}{2}\left(\frac{b^{2}}{4}-m\right)\right]^{\frac{1}{\alpha}}}e^{-\frac{2}{b}\lambda^{\alpha}u_{k}^{\alpha}t}\lambda^{n+\alpha}\mathrm{d}\lambda\left(\|u_{0}\|^{2}_{L^{1}(\mathbb{H}^{n})}+\|u_{1}\|^{2}_{L^{1}(\mathbb{H}^{n})}\right)
.\end{align}
 Using the substitution $\theta=\frac{2}{b}\lambda^{\alpha}u_{k}^{\alpha}t$, we have $\lambda^{n+\alpha}\mathrm{d}\lambda =\frac{1}{\alpha}\left(\frac{b}{2}\right)^{\frac{n+1}{\alpha}+1}\theta^{\frac{n+1}{\alpha}}\mu_{k}^{-(n+1+\alpha)}t^{-\left(\frac{n+1}{\alpha}+1\right)}\mathrm{d}\theta$
 and consequently, similar to the estimate \eqref{ilowestimate}, using the convergence of infinite series and  gamma function, we deduce that
\begin{align}\label{jlowestimate}
	&\mathrm{J}_{\mathrm{low}} \lesssim t^{-\left(\frac{n+1}{\alpha}+1\right)}e^{-\frac{2m}{b}t} \sum_{k\in\mathbb{N}^{n}}\mu_{k}^{\alpha}\mu_{k}^{-(n+1+\alpha)} 
	\int_{0}^{\frac{1}{b}\left(\frac{b^{2}}{4}-m\right)t}e^{-\theta}\theta^{\frac{n+1}{\alpha}}\mathrm{d}\theta\left(\|u_{0}\|^{2}_{L^{1}(\mathbb{H}^{n})}+\|u_{1}\|^{2}_{L^{1}(\mathbb{H}^{n})}\right)\nonumber\\
	&\lesssim t^{-\left(\frac{n+1}{\alpha}+1\right)}e^{-\frac{2m}{b}t}\sum_{k\in\mathbb{N}^{n}}(2|k|+n)^{-(n+1)}\Gamma\left(\frac{n+1}{\alpha}+1\right)\left(\|u_{0}\|^{2}_{L^{1}(\mathbb{H}^{n})}+\|u_{1}\|^{2}_{L^{1}(\mathbb{H}^{n})}\right)\nonumber\\
	&\lesssim t^{-\frac{\mathcal{Q}}{2\alpha}-1}e^{-\frac{2m}{b}t}\left(\|u_{0}\|^{2}_{L^{1}(\mathbb{H}^{n})}+\|u_{1}\|^{2}_{L^{1}(\mathbb{H}^{n})}\right).
\end{align}
By making the sharper intermediate steps, we may also get some decay by using just $L^{2}$-regularity. Let us recalculate the term $\operatorname{J}_{\operatorname{low}}$ by using the Plancherel formula instead of Parseval's identity:
\begin{align}\label{jlowl20}
	&\mathrm{J}_{\mathrm{low}} \lesssim\sum_{k,\ell\in\mathbb{N}^{n}}\int_{0<|\lambda|<\frac{1}{\mu_{k}}\left[\frac{1}{2}\left(\frac{b^{2}}{4}-m\right)\right]^{\frac{1}{\alpha}}}|\lambda|^{\alpha}\mu_{k}^{\alpha}e^{-\frac{2}{b}\left(m+|\lambda|^{\alpha}u_{k}^{\alpha}\right)t}\left(	|\widehat{u}_{0}(\lambda)_{k,\ell}|^{2}+	|\widehat{u}_{1}(\lambda)_{k,\ell}|^{2}\right)|\lambda|^{n}\mathrm{d}\lambda\nonumber\\
	&\lesssim e^{-\frac{2m}{b}t}\sum_{k,\ell\in\mathbb{N}^{n}}\int_{0<|\lambda|<\frac{1}{\mu_{k}}\left[\frac{1}{2}\left(\frac{b^{2}}{4}-m\right)\right]^{\frac{1}{\alpha}}}t^{-1}\left(	|\widehat{u}_{0}(\lambda)_{k,\ell}|^{2}+	|\widehat{u}_{1}(\lambda)_{k,\ell}|^{2}\right)|\lambda|^{n}\mathrm{d}\lambda\nonumber\\
	&\lesssim t^{-1}e^{-\frac{2m}{b}t}\int_{\mathbb{R}^{*}}\sum_{k,\ell\in\mathbb{N}^{n}}\left(	|\widehat{u}_{0}(\lambda)_{k,\ell}|^{2}+	|\widehat{u}_{1}(\lambda)_{k,\ell}|^{2}\right)|\lambda|^{n}\mathrm{d}\lambda\nonumber\\
	&=t^{-1}e^{-\frac{2m}{b}t}\left(\|u_{0}\|^{2}_{L^{2}(\mathbb{H}^{n})}+\|u_{1}\|^{2}_{L^{2}(\mathbb{H}^{n})}\right) ,
\end{align}
where we utilized the Remark \ref{estexp} in the 2nd step. 

Now we will estimate the term $\mathrm{J}_{\mathrm{high}}$ given by
\begin{equation*}
\mathrm{J}_{\mathrm{high}}=c_{n}\sum_{k,\ell\in\mathbb{N}^{n}} \int_{|\lambda|>\frac{1}{\mu_{k}}\left[\frac{1}{2}\left(\frac{b^{2}}{4}-m\right)\right]^{\frac{1}{\alpha}}}(|\lambda|^{\alpha/2}\mu_{k}^{\alpha/2}\widehat{u}(t,\lambda)e_{k},e_{\ell})^{2}_{L^{2}(\mathbb{R}^{n})}|\lambda|^{n}\mathrm{d}\lambda.
\end{equation*}
Using the estimate \eqref{largefreq22} for  the above equation, we have
\begin{align}\label{jhighestimatedel}
	\mathrm{J}_{\mathrm{high}}
	&\lesssim e^{\left(-b+\sqrt{\frac{1}{2}\left(b^{2}-4m\right)}\right) t}\sum_{k,\ell\in\mathbb{N}^{n}} \int_{|\lambda|>\frac{1}{\mu_{k}}\left[\frac{1}{2}\left(\frac{b^{2}}{4}-m\right)\right]^{\frac{1}{\alpha}}} \left((|\lambda|\mu_{k})^{\alpha}	|\widehat{u}_{0}(\lambda)_{k,\ell}|^{2}+	|\widehat{u}_{1}(\lambda)_{k,\ell}|^{2}\right)|\lambda|^{n}\mathrm{d}\lambda\nonumber\\
	&\lesssim e^{\left(-b+\sqrt{\frac{1}{2}\left(b^{2}-4m\right)}\right) t}\int_{\mathbb{R}^{*}}\sum_{k,\ell\in\mathbb{N}^{n}}  \left(	(|\lambda|\mu_{k})^{\alpha}|\widehat{u}_{0}(\lambda)_{k,\ell}|^{2}+	|\widehat{u}_{1}(\lambda)_{k,\ell}|^{2}\right)|\lambda|^{n}\mathrm{d}\lambda\nonumber\\
	&= e^{\left(-b+\sqrt{\frac{1}{2}\left(b^{2}-4m\right)}\right) t}  \left(	\left\|(-\mathcal{L}_{\mathbb{H}})^{\alpha/2}u_{0}\right\|^{2}_{L^{2}(\mathbb{H}^{n})}+	\|u_{1}\|^{2}_{L^{2}(\mathbb{H}^{n})}\right).
\end{align}
Let us now combine the estimates $\mathrm{J}_{\operatorname{low}}$ and $\mathrm{J}_{\operatorname{high}}$. Similar to the estimate \eqref{expest}, using the Remark \ref{estexp} for the expression \eqref{equalityrel}, we deduce that
\begin{equation}\label{expest10}
	e^{\left(-b+\sqrt{\frac{1}{2}\left(b^{2}-4m\right)}\right) t}\leq C_{2}(b,m,\mathcal{Q},\alpha)t^{-\frac{\mathcal{Q}}{2\alpha}-1}e^{-\frac{m}{b}t},
\end{equation}
and 
\begin{equation}\label{expest20}
	e^{\left(-b+\sqrt{\frac{1}{2}\left(b^{2}-4m\right)}\right) t}\leq C_{3}(b,m)t^{-1}e^{-\frac{m}{b}t},
\end{equation}
for some positive constants $C_{2}(b,m,\mathcal{Q},\alpha)$ and $C_{3}(b,m)$. Now combining the estimates \eqref{jlowestimate} and \eqref{jhighestimatedel} with \eqref{expest10}, we get
\begin{multline}\label{lamunorml1l2}
	\left\|(-\mathcal{L}_{\mathbb{H}})^{\alpha/2}u(t, \cdot)\right\|^{2}_{L^2\left(\mathbb{H}^{n}\right)}\\\lesssim t^{-\frac{\mathcal{Q}}{2\alpha}-1}e^{-\frac{m}{b}t}\left(\left\|(-\mathcal{L}_{\mathbb{H}})^{\alpha/2}u_{0}\right\|^{2}_{L^{2}\left(\mathbb{H}^{n}\right)}+\left\|u_{1}\right\|^{2}_{L^{2}\left(\mathbb{H}^{n}\right)}+\left\|u_{0}\right\|^{2}_{L^{1}\left(\mathbb{H}^{n}\right)}+\left\|u_{1}\right\|^{2}_{L^{1}\left(\mathbb{H}^{n}\right)}\right),
\end{multline}
and combining the estimates \eqref{jlowl20} and \eqref{jhighestimatedel} with \eqref{expest20}, we get
\begin{equation}\label{lamunorml2}
	\left\|(-\mathcal{L}_{\mathbb{H}})^{\alpha/2}u(t, \cdot)\right\|^{2}_{L^2\left(\mathbb{H}^{n}\right)}\lesssim t^{-1}e^{-\frac{m}{b}t}\left(\left\|(-\mathcal{L}_{\mathbb{H}})^{\alpha/2}u_{0}\right\|^{2}_{L^{2}\left(\mathbb{H}^{n}\right)}+\left\|u_{1}\right\|^{2}_{L^{2}\left(\mathbb{H}^{n}\right)}\right).
\end{equation}
Similar to the case of $	\|u(t, \cdot)\|^{2}_{L^2\left(\mathbb{H}^{n}\right)}$, if we utilize the uniform estimate \eqref{uniformestforlamu} instead of \eqref{smallfreq11} and \eqref{largefreq22}, then we have
\begin{equation*}
	\left\|(-\mathcal{L}_{\mathbb{H}})^{\alpha/2}u(t, \cdot)\right\|^{2}_{L^2\left(\mathbb{H}^{n}\right)}\lesssim 	e^{\left(-b+\sqrt{b^{2}-4m}\right) t}\left(\left\|(-\mathcal{L}_{\mathbb{H}})^{\alpha/2}u_{0}\right\|^{2}_{L^{2}\left(\mathbb{H}^{n}\right)}+\left\|u_{1}\right\|^{2}_{L^{2}\left(\mathbb{H}^{n}\right)}\right).
\end{equation*}
Using the argument similar to the previous estimate $\|u(t,\cdot)\|_{L^{2}(\mathbb{H}^{n})}$, we can also rule out the singular behaviour of 	$\left\|(-\mathcal{L}_{\mathbb{H}})^{\alpha/2}u(t, \cdot)\right\|_{L^2\left(\mathbb{H}^{n}\right)}$ as $t\to 0^{+}$ and hence the estimate \eqref{lamunorml1l2} and \eqref{lamunorml2} can be replaced by
\begin{multline}\label{2final1}
	\left\|(-\mathcal{L}_{\mathbb{H}})^{\alpha/2}u(t, \cdot)\right\|^{2}_{L^2\left(\mathbb{H}^{n}\right)}\\\lesssim (1+t)^{-\frac{\mathcal{Q}}{2\alpha}-1}e^{-\frac{m}{b}t}\left(\left\|(-\mathcal{L}_{\mathbb{H}})^{\alpha/2}u_{0}\right\|^{2}_{L^{2}\left(\mathbb{H}^{n}\right)}+\left\|u_{1}\right\|^{2}_{L^{2}\left(\mathbb{H}^{n}\right)}+\left\|u_{0}\right\|^{2}_{L^{1}\left(\mathbb{H}^{n}\right)}+\left\|u_{1}\right\|^{2}_{L^{1}\left(\mathbb{H}^{n}\right)}\right).
\end{multline}
and 
\begin{equation}\label{2final2}
	\left\|(-\mathcal{L}_{\mathbb{H}})^{\alpha/2}u(t, \cdot)\right\|^{2}_{L^2\left(\mathbb{H}^{n}\right)}\lesssim (1+t)^{-1}e^{-\frac{m}{b}t}\left(\left\|(-\mathcal{L}_{\mathbb{H}})^{\alpha/2}u_{0}\right\|^{2}_{L^{2}\left(\mathbb{H}^{n}\right)}+\left\|u_{1}\right\|^{2}_{L^{2}\left(\mathbb{H}^{n}\right)}\right),
\end{equation}
respectively. The estimates \eqref{2final1} and \eqref{2final2} will instantly give our required estimates \eqref{delul1l2final} and \eqref{delul2l2final}, respectively.
\subsection{Estimates for  $\left\|\partial_{t}u(t,\cdot)\right\|_{L^{2}(\mathbb{H}^{n})}$}\label{subsec3.3}
 We now estimate the $L^{2}$-norm of $\partial_{t}u(t,\eta)$, i.e.,  $\left\|\partial_{t}u(t,\cdot)\right\|_{L^{2}(\mathbb{H}^{n})}$. Using the Plancherel formula \eqref{Plancherel}, we have 
\begin{align}\nonumber
	&\|\partial_{t}u(t,\cdot)\|^{2}_{L^{2}(\mathbb{H}^{n})}\\\nonumber
	&=c_{n}\sum_{k,\ell\in\mathbb{N}^{n}}\left(\int_{0<|\lambda|<\frac{1}{\mu_{k}}\left[\frac{1}{2}\left(\frac{b^{2}}{4}-m\right)\right]^{\frac{1}{\alpha}}}+\int_{|\lambda|>\frac{1}{mu_{k}}\left[\frac{1}{2}\left(\frac{b^{2}}{4}-m\right)\right]^{\frac{1}{\alpha}}}\right)(\partial_{t}\widehat{u}(t,\lambda)e_{k},e_{\ell})^{2}_{L^{2}(\mathbb{R}^{n})}|\lambda|^{n}\mathrm{d}\lambda\\
	&=\mathrm{K}_{\operatorname{low}}+\mathrm{K}_{\operatorname{high}}.
\end{align}
 Let us begin to estimate $\mathrm{K}_{\operatorname{low}}$ given by
$$\mathrm{K}_{\operatorname{low}}=c_{n}\sum_{k,\ell\in\mathbb{N}^{n}} \int_{0<|\lambda|<\frac{1}{\mu_{k}}\left[\frac{1}{2}\left(\frac{b^{2}}{4}-m\right)\right]^{\frac{1}{\alpha}}} (\partial_{t}\widehat{u}(t,\lambda)e_{k},e_{\ell})^{2}_{L^{2}(\mathbb{R}^{n})}|\lambda|^{n}\mathrm{d}\lambda.$$
Utilizing the relation \eqref{relat}, the estimate \eqref{smallfreq2} can be further estimated as
\begin{equation*}
	|	\partial_{t}	\widehat{u}(t,\lambda)_{k,\ell}|\lesssim e^{-\frac{1}{b}\left(m+|\lambda|^{\alpha}u_{k}^{\alpha}\right)t}(|\lambda|^{\alpha}\mu_{k}^{\alpha}+m)\left(|\widehat{u}_{0}(\lambda)_{k,\ell}|+|\widehat{u}_{1}(\lambda)_{k,\ell}|\right)+e^{-\frac{b}{2}t}|\widehat{u}_{1}(\lambda)_{k,\ell}|.
\end{equation*}
Using the above estimate and the relation \eqref{fhatnorm}, analogous to the estimate \eqref{ilowest}, we can obtain
\begin{align}
	&\mathrm{K}_{\mathrm{low}} \lesssim\sum_{k\in\mathbb{N}^{n}}\mu_{k}^{2\alpha}\int_{0<|\lambda|<\frac{1}{\mu_{k}}\left[\frac{1}{2}\left(\frac{b^{2}}{4}-m\right)\right]^{\frac{1}{\alpha}}}e^{-\frac{2}{b}\left(m+|\lambda|^{\alpha}u_{k}^{\alpha}\right)t}\sum_{\ell\in\mathbb{N}^{n}}\left(	|\widehat{u}_{0}(\lambda)_{k,\ell}|^{2}+	|\widehat{u}_{1}(\lambda)_{k,
		\ell}|^{2}\right)|\lambda|^{n+2\alpha}\mathrm{d}\lambda\nonumber\\
	&+m^{2}\sum_{k\in\mathbb{N}^{n}}\int_{0<|\lambda|<\frac{1}{\mu_{k}}\left[\frac{1}{2}\left(\frac{b^{2}}{4}-m\right)\right]^{\frac{1}{\alpha}}}e^{-\frac{2}{b}\left(m+|\lambda|^{\alpha}u_{k}^{\alpha}\right)t}\sum_{\ell\in\mathbb{N}^{n}}\left(	|\widehat{u}_{0}(\lambda)_{k,\ell}|^{2}+	|\widehat{u}_{1}(\lambda)_{k,
		\ell}|^{2}\right)|\lambda|^{n}\mathrm{d}\lambda\nonumber\\
	&+e^{-bt}\sum_{k\in\mathbb{N}^{n}}\int_{0<|\lambda|<\frac{1}{\mu_{k}}\left[\frac{1}{2}\left(\frac{b^{2}}{4}-m\right)\right]^{\frac{1}{\alpha}}}\sum_{\ell\in\mathbb{N}^{n}}	|\widehat{u}_{1}(\lambda)_{k,\ell}|^{2}|\lambda|^{n}\mathrm{d}\lambda\nonumber\\
	&=\sum_{k\in\mathbb{N}^{n}}\mu_{k}^{2\alpha}\int_{0<|\lambda|<\frac{1}{\mu_{k}}\left[\frac{1}{2}\left(\frac{b^{2}}{4}-m\right)\right]^{\frac{1}{\alpha}}}e^{-\frac{2}{b}\left(m+|\lambda|^{\alpha}u_{k}^{\alpha}\right)t}\left(	\|\widehat{u}_{0}(\lambda)e_{k}\|_{L^{2}(\mathbb{R}^{n})}^{2}+	\|\widehat{u}_{1}(\lambda)e_{k}\|_{L^{2}(\mathbb{R}^{n})}^{2}\right)|\lambda|^{n+2\alpha}\mathrm{d}\lambda\nonumber\\
	&+m^{2}\sum_{k\in\mathbb{N}^{n}}\int_{0<|\lambda|<\frac{1}{\mu_{k}}\left[\frac{1}{2}\left(\frac{b^{2}}{4}-m\right)\right]^{\frac{1}{\alpha}}}e^{-\frac{2}{b}\left(m+|\lambda|^{\alpha}u_{k}^{\alpha}\right)t}\left(	\|\widehat{u}_{0}(\lambda)e_{k}\|_{L^{2}(\mathbb{R}^{n})}^{2}+	\|\widehat{u}_{1}(\lambda)e_{k}\|_{L^{2}(\mathbb{R}^{n})}^{2}\right)|\lambda|^{n}\mathrm{d}\lambda\nonumber\\
	&+e^{-bt}\sum_{k\in\mathbb{N}^{n}}\int_{0<|\lambda|<\frac{1}{\mu_{k}}\left[\frac{1}{2}\left(\frac{b^{2}}{4}-m\right)\right]^{\frac{1}{\alpha}}}	\|\widehat{u}_{1}(\lambda)e_{k}\|_{L^{2}(\mathbb{R}^{n})}^{2}|\lambda|^{n}\mathrm{d}\lambda\nonumber\\
	&\lesssim e^{-\frac{2m}{b}t}\sum_{k\in\mathbb{N}^{n}}\mu_{k}^{2\alpha}\int_{0<|\lambda|<\frac{1}{\mu_{k}}\left[\frac{1}{2}\left(\frac{b^{2}}{4}-m\right)\right]^{\frac{1}{\alpha}}}e^{-\frac{2}{b}|\lambda|^{\alpha}u_{k}^{\alpha}t}|\lambda|^{n+2\alpha}\mathrm{d}\lambda\left(\|u_{0}\|^{2}_{L^{1}(\mathbb{H}^{n})}+\|u_{1}\|^{2}_{L^{1}(\mathbb{H}^{n})}\right) \nonumber\\
	&+m^{2}e^{-\frac{2m}{b}t}\sum_{k\in\mathbb{N}^{n}}\int_{0<|\lambda|<\frac{1}{\mu_{k}}\left[\frac{1}{2}\left(\frac{b^{2}}{4}-m\right)\right]^{\frac{1}{\alpha}}}e^{-\frac{2}{b}|\lambda|^{\alpha}u_{k}^{\alpha}t}|\lambda|^{n}\mathrm{d}\lambda\left(\|u_{0}\|^{2}_{L^{1}(\mathbb{H}^{n})}+\|u_{1}\|^{2}_{L^{1}(\mathbb{H}^{n})}\right) \nonumber\\
	&+e^{-bt}  \sum_{k\in\mathbb{N}^{n}}\int_{0<|\lambda|<\frac{1}{\mu_{k}}\left[\frac{1}{2}\left(\frac{b^{2}}{4}-m\right)\right]^{\frac{1}{\alpha}}}|\lambda|^{n}\mathrm{d}\lambda \|u_{1}\|^{2}_{L^{1}(\mathbb{H}^{n})} \nonumber\\
	&=e^{-\frac{2m}{b}t}\sum_{k\in\mathbb{N}^{n}}\mu_{k}^{2\alpha}\int_{0}^{\frac{1}{\mu_{k}}\left[\frac{1}{2}\left(\frac{b^{2}}{4}-m\right)\right]^{\frac{1}{\alpha}}}e^{-\frac{2}{b}\lambda^{\alpha}u_{k}^{\alpha}t}\lambda^{n+2\alpha}\mathrm{d}\lambda\left(\|u_{0}\|^{2}_{L^{1}(\mathbb{H}^{n})}+\|u_{1}\|^{2}_{L^{1}(\mathbb{H}^{n})}\right) \nonumber\\
	&+m^{2}e^{-\frac{2m}{b}t}\sum_{k\in\mathbb{N}^{n}}\int_{0}^{\frac{1}{\mu_{k}}\left[\frac{1}{2}\left(\frac{b^{2}}{4}-m\right)\right]^{\frac{1}{\alpha}}}e^{-\frac{2}{b}\lambda^{\alpha}u_{k}^{\alpha}t}\lambda^{n}\mathrm{d}\lambda\left(\|u_{0}\|^{2}_{L^{1}(\mathbb{H}^{n})}+\|u_{1}\|^{2}_{L^{1}(\mathbb{H}^{n})}\right) \nonumber\\
	&+e^{-bt}\sum_{k\in\mathbb{N}^{n}}\int_{0}^{\frac{1}{\mu_{k}}\left[\frac{1}{2}\left(\frac{b^{2}}{4}-m\right)\right]^{\frac{1}{\alpha}}}\lambda^{n}\mathrm{d}\lambda \|u_{1}\|^{2}_{L^{1}(\mathbb{H}^{n})}.
\end{align}
 Using the substitute $\theta=\frac{2}{b}\lambda^{\alpha}u_{k}^{\alpha}t$ in the  first integral, we have
  \begin{equation*}
 \lambda^{n+2\alpha}\mathrm{d}\lambda =\frac{1}{\alpha}\left(\frac{b}{2}\right)^{\frac{n+1}{\alpha}+2}\theta^{\frac{n+1+\alpha}{\alpha}}\mu_{k}^{-(n+1+2\alpha)}t^{-\left(\frac{n+1}{\alpha}+2\right)}\mathrm{d}\theta,
 \end{equation*}
  and using estimate \eqref{ilowestimate} for the 2nd integral, we get
\begin{align}\label{klowest}
	&\mathrm{K}_{\mathrm{low}} \lesssim t^{-\left(\frac{n+1}{\alpha}+2\right)}  e^{-\frac{2m}{b}t} \sum_{k\in\mathbb{N}^{n}}\mu_{k}^{2\alpha}\mu_{k}^{-(n+1+2\alpha)}  \int_{0}^{\frac{1}{b}\left(\frac{b^{2}}{4}-m\right)t}e^{-\theta}\theta^{\frac{n+1+\alpha}{\alpha}}\mathrm{d}\theta\left(\|u_{0}\|^{2}_{L^{1}(\mathbb{H}^{n})}+\|u_{1}\|^{2}_{L^{1}(\mathbb{H}^{n})}\right)\nonumber\\
	&+m^{2}t^{-\frac{\mathcal{Q}}{2\alpha}}e^{-\frac{2m}{b}t}\left(\|u_{0}\|^{2}_{L^{1}(\mathbb{H}^{n})}+\|u_{1}\|^{2}_{L^{1}(\mathbb{H}^{n})}\right)+e^{-bt} \sum_{k\in\mathbb{N}^{n}}\mu_{k}^{-(n+1)}\|u_{1}\|^{2}_{L^{1}(\mathbb{H}^{n})} \nonumber\\
	&\lesssim t^{-\left(\frac{n+1}{\alpha}+2\right)}e^{-\frac{2m}{b}t}\sum_{k\in\mathbb{N}^{n}}(2|k|+n)^{-(n+1)}\Gamma\left(\frac{n+1}{\alpha}+2\right) \left(\|u_{0}\|^{2}_{L^{1}(\mathbb{H}^{n})}+\|u_{1}\|^{2}_{L^{1}(\mathbb{H}^{n})}\right) \nonumber\\
	&+m^{2}t^{-\frac{\mathcal{Q}}{2\alpha}}e^{-\frac{2m}{b}t}\left(\|u_{0}\|^{2}_{L^{1}(\mathbb{H}^{n})}+\|u_{1}\|^{2}_{L^{1}(\mathbb{H}^{n})}\right)+e^{-bt}\sum_{k\in\mathbb{N}^{n}}(2|k|+n)^{-(n+1)}\|u_{1}\|^{2}_{L^{1}(\mathbb{H}^{n})}\nonumber\\
	&\lesssim t^{-\left(\frac{n+1}{\alpha}+2\right)}e^{-\frac{2m}{b}t}\left(\|u_{0}\|^{2}_{L^{1}(\mathbb{H}^{n})}+\|u_{1}\|^{2}_{L^{1}(\mathbb{H}^{n})}\right)+m^{2}t^{-\frac{\mathcal{Q}}{2\alpha}}e^{-\frac{2m}{b}t}\left(\|u_{0}\|^{2}_{L^{1}(\mathbb{H}^{n})}+\|u_{1}\|^{2}_{L^{1}(\mathbb{H}^{n})}\right)\nonumber\\
	&+e^{-bt}\|u_{1}\|^{2}_{L^{1}(\mathbb{H}^{n})}\nonumber\\
	&\lesssim \left[t^{-\frac{\mathcal{Q}}{2\alpha}-2}+m^{2}t^{-\frac{\mathcal{Q}}{2\alpha}}\right]e^{-\frac{2m}{b}t}\left(\|u_{0}\|^{2}_{L^{1}(\mathbb{H}^{n})}+\|u_{1}\|^{2}_{L^{1}(\mathbb{H}^{n})}\right)+e^{-bt}\|u_{1}\|^{2}_{L^{1}(\mathbb{H}^{n})}.
\end{align}
As in the preceding instance, we can also get some decay by utilizing just $L^{2}$-regularity by making the sharper intermediate steps. Let us recalculate $\mathrm{K}_{\mathrm{low}}$ with $L^{2}$-regularity only
\begin{align}\label{klowl2}
	&\mathrm{K}_{\mathrm{low}} \lesssim\sum_{k,\ell\in\mathbb{N}^{n}}\int_{0<|\lambda|<\frac{1}{\mu_{k}}\left[\frac{1}{2}\left(\frac{b^{2}}{4}-m\right)\right]^{\frac{1}{\alpha}}}(|\lambda|
	\mu_{k})^{2\alpha}e^{-\frac{2}{b}\left(m+|\lambda|^{\alpha}u_{k}^{\alpha}\right)t}\left(	|\widehat{u}_{0}(\lambda)_{k,\ell}|^{2}+	|\widehat{u}_{1}(\lambda)_{k,
		\ell}|^{2}\right)|\lambda|^{n}\mathrm{d}\lambda\nonumber\\
	&+m^{2}\sum_{k,\ell\in\mathbb{N}^{n}}\int_{0<|\lambda|<\frac{1}{\mu_{k}}\left[\frac{1}{2}\left(\frac{b^{2}}{4}-m\right)\right]^{\frac{1}{\alpha}}}e^{-\frac{2}{b}\left(m+|\lambda|^{\alpha}u_{k}^{\alpha}\right)t}\left(	|\widehat{u}_{0}(\lambda)_{k,\ell}|^{2}+	|\widehat{u}_{1}(\lambda)_{k,
		\ell}|^{2}\right)|\lambda|^{n}\mathrm{d}\lambda\nonumber\\
	&+e^{-bt}\sum_{k,\ell\in\mathbb{N}^{n}}\int_{0<|\lambda|<\frac{1}{\mu_{k}}\left[\frac{1}{2}\left(\frac{b^{2}}{4}-m\right)\right]^{\frac{1}{\alpha}}}	|\widehat{u}_{1}(\lambda)_{k,\ell}|^{2}|\lambda|^{n}\mathrm{d}\lambda\nonumber\\
	&\lesssim t^{-2}e^{-\frac{2m}{b}t}\sum_{k,\ell\in\mathbb{N}^{n}}\int_{0<|\lambda|<\frac{1}{\mu_{k}}\left[\frac{1}{2}\left(\frac{b^{2}}{4}-m\right)\right]^{\frac{1}{\alpha}}}\left(	|\widehat{u}_{0}(\lambda)_{k,\ell}|^{2}+	|\widehat{u}_{1}(\lambda)_{k,
		\ell}|^{2}\right)|\lambda|^{n}\mathrm{d}\lambda\nonumber\\
	&+m^{2}e^{-\frac{2m}{b}t}\sum_{k,\ell\in\mathbb{N}^{n}}\int_{0<|\lambda|<\frac{1}{\mu_{k}}\left[\frac{1}{2}\left(\frac{b^{2}}{4}-m\right)\right]^{\frac{1}{\alpha}}}\left(	|\widehat{u}_{0}(\lambda)_{k,\ell}|^{2}+	|\widehat{u}_{1}(\lambda)_{k,
		\ell}|^{2}\right)|\lambda|^{n}\mathrm{d}\lambda\nonumber\\
	&+e^{-bt}\sum_{k,\ell\in\mathbb{N}^{n}}\int_{0<|\lambda|<\frac{1}{\mu_{k}}\left[\frac{1}{2}\left(\frac{b^{2}}{4}-m\right)\right]^{\frac{1}{\alpha}}}	|\widehat{u}_{1}(\lambda)_{k,\ell}|^{2}|\lambda|^{n}\mathrm{d}\lambda\nonumber\\
	&\lesssim \left[t^{-2}+m^{2}\right]e^{-\frac{2m}{b}t}\int_{\mathbb{R}^{*}}\sum_{k,\ell\in\mathbb{N}^{n}}\left(	|\widehat{u}_{0}(\lambda)_{k,\ell}|^{2}+	|\widehat{u}_{1}(\lambda)_{k,
		\ell}|^{2}\right)|\lambda|^{n}\mathrm{d}\lambda\nonumber\\
		&+e^{-bt}\sum_{k,\ell\in\mathbb{N}^{n}}\int_{\mathbb{R}^{*}}	|\widehat{u}_{1}(\lambda)_{k,\ell}|^{2}|\lambda|^{n}\mathrm{d}\lambda\nonumber\\
	&=\left[t^{-2}+m^{2}\right]e^{-\frac{2m}{b}t} \left(	\left\|u_{0}\right\|^{2}_{L^{2}(\mathbb{H}^{n})}+	\|u_{1}\|^{2}_{L^{2}(\mathbb{H}^{n})}\right)+e^{-bt}\|u_{1}\|^{2}_{L^{2}(\mathbb{H}^{n})},
\end{align}
where we utilized Remark \ref{estexp} to bound the exponential factor in the 1st and 2nd integrals in the 3rd step.

Now we will estimate the term $\mathrm{K}_{\mathrm{high}}$ given by 
\begin{equation*}
\mathrm{K}_{\mathrm{high}}=c_{n}\int_{|\lambda|>\frac{1}{\mu_{k}}\left[\frac{1}{2}\left(\frac{b^{2}}{4}-m\right)\right]^{\frac{1}{\alpha}}}(\partial_{t}\widehat{u}(t,\lambda)e_{k},e_{\ell})^{2}_{L^{2}(\mathbb{R}^{n})}|\lambda|^{n}\mathrm{d}\lambda.
\end{equation*}
Using the estimate \eqref{largefreq2}, we get
\begin{align}\label{khighestimate}
	&\mathrm{K}_{\mathrm{high}} \lesssim e^{\left(-b+\sqrt{\frac{1}{2}\left(b^{2}-4m\right)}\right) t}\sum_{k,\ell\in\mathbb{N}^{n}} \int_{|\lambda|>\frac{1}{\mu_{k}}\left[\frac{1}{2}\left(\frac{b^{2}}{4}-m\right)\right]^{\frac{1}{\alpha}}} \left((|\lambda|^{\alpha}\mu_{k}^{\alpha}+m)	|\widehat{u}_{0}(\lambda)_{k,\ell}|^{2}+	|\widehat{u}_{1}(\lambda)_{k,\ell}|^{2}\right)|\lambda|^{n}\mathrm{d}\lambda\nonumber\\
	&\lesssim e^{\left(-b+\sqrt{\frac{1}{2}\left(b^{2}-4m\right)}\right) t}\int_{\mathbb{R}^{*}}\sum_{k,\ell\in\mathbb{N}^{n}}  \left(	(|\lambda|^{\alpha}\mu_{k}^{\alpha}+m)|\widehat{u}_{0}(\lambda)_{k,\ell}|^{2}+	|\widehat{u}_{1}(\lambda)_{k,\ell}|^{2}\right)|\lambda|^{n}\mathrm{d}\lambda\nonumber\\
	&=e^{\left(-b+\sqrt{\frac{1}{2}\left(b^{2}-4m\right)}\right) t}  \left(	\left\|(-\mathcal{L}_{\mathbb{H}})^{\alpha/2}u_{0}\right\|^{2}_{L^{2}(\mathbb{H}^{n})}+	m\left\|u_{0}\right\|^{2}_{L^{2}(\mathbb{H}^{n})}+\|u_{1}\|^{2}_{L^{2}(\mathbb{H}^{n})}\right).
\end{align}
 Note that
\begin{equation*}
	-b=-\frac{2m}{b}-f_{2}(b,m),\quad \text{where }	f_{2}(b,m)=b-\frac{2m}{b}>\begin{cases}
		\frac{2m}{b}&\text{if } m>0,\\
		b &\text{if }m=0.
	\end{cases}
\end{equation*}
Now from Remark \ref{estexp} it follows that there exist constants $C_{4}(b,m,\mathcal{Q},\alpha)$ and $C_{5}(b,m)$ such that,
\begin{equation}\label{expest3}
	e^{-b t}=e^{-\frac{2m}{b}t}e^{-f_{2}(b,m)t}\leq C_{4}(b,m,\mathcal{Q},\alpha)t^{-\frac{\mathcal{Q}}{2\alpha}-2}e^{-\frac{2m}{b}t}
\end{equation}
and
\begin{equation}\label{expest4}
	e^{-b t}=e^{-\frac{2m}{b}t}e^{-f_{2}(b,m)t}\leq C_{5}(b,m)t^{-2}e^{-\frac{2m}{b}t}.
\end{equation}
Similar to the estimate \eqref{expest}, using the expression \eqref{equalityrel}, we deduce that
\begin{equation}\label{expest100}
	e^{\left(-b+\sqrt{\frac{1}{2}\left(b^{2}-4m\right)}\right) t}\leq C_{6}(b,m,\mathcal{Q},\alpha)t^{-\frac{\mathcal{Q}}{2\alpha}-2}e^{-\frac{m}{b}t}
\end{equation}
and 
\begin{equation}\label{expest200}
	e^{\left(-b+\sqrt{\frac{1}{2}\left(b^{2}-4m\right)}\right) t}\leq C_{7}(b,m)t^{-2}e^{-\frac{m}{b}t},
\end{equation}
for some positive constants $C_{6}(b,m,\mathcal{Q},\alpha)$ and $C_{7}(b,m)$. Now combining the estimates \eqref{klowest}, \eqref{khighestimate}, \eqref{expest3}, and \eqref{expest100}, we deduce that
\begin{multline}\label{deltfinal1}
	\|\partial_{t}u(t, \cdot)\|^{2}_{L^2\left(\mathbb{H}^{n}\right)}\lesssim\left[t^{-\frac{\mathcal{Q}}{2\alpha}-2}+m^{2}t^{-\frac{\mathcal{Q}}{2\alpha}}\right] e^{-\frac{m}{b}t}\\\left(\left\|(-\mathcal{L}_{\mathbb{H}})^{\alpha/2}u_{0}\right\|^{2}_{L^{2}\left(\mathbb{H}^{n}\right)}+\left\|u_{1}\right\|^{2}_{L^{2}\left(\mathbb{H}^{n}\right)}+\left\|u_{0}\right\|^{2}_{L^{1}\left(\mathbb{H}^{n}\right)}+\left\|u_{1}\right\|^{2}_{L^{1}\left(\mathbb{H}^{n}\right)}\right).
\end{multline}
and combining the estimates \eqref{klowest}, \eqref{khighestimate}, \eqref{expest4}, and \eqref{expest200}, we get
\begin{equation}\label{deltfinal2}
	\|\partial_{t}u(t, \cdot)\|^{2}_{L^2\left(\mathbb{H}^{n}\right)}\lesssim\\\left[t^{-2}+m^{2}\right] e^{-\frac{m}{b}t}\left(\left\|(-\mathcal{L}_{\mathbb{H}})^{\alpha/2}u_{0}\right\|^{2}_{L^{2}\left(\mathbb{H}^{n}\right)}+\left\|u_{1}\right\|^{2}_{L^{2}\left(\mathbb{H}^{n}\right)}\right).
\end{equation}
Similar to the previous cases, if we utilize the uniform estimate \eqref{uniformestfordeltu} instead of \eqref{smallfreq2} and \eqref{largefreq2}, then we have
\begin{equation}
	\left\|\partial_{t}u(t, \cdot)\right\|^{2}_{L^2\left(\mathbb{H}^{n}\right)}\lesssim 	e^{\left(-b+\sqrt{b^{2}-4m}\right) t}\left(\left\|(-\mathcal{L}_{\mathbb{H}})^{\alpha/2}u_{0}\right\|^{2}_{L^{2}\left(\mathbb{H}^{n}\right)}+\left\|u_{1}\right\|^{2}_{L^{2}\left(\mathbb{H}^{n}\right)}\right).
\end{equation}
In the light of above estimate, similar to the previous cases, we can neglect the singular behaviour of $	\left\|\partial_{t}u(t, \cdot)\right\|_{L^2\left(\mathbb{H}^{n}\right)}$ as $t\to 0^{+}$ and consequently the estimates \eqref{deltfinal1} and \eqref{deltfinal2} takes the form
\begin{multline}\label{nonsingestdelt}
	\|\partial_{t}u(t, \cdot)\|^{2}_{L^2\left(\mathbb{H}^{n}\right)}\lesssim\left[(1+t)^{-\frac{\mathcal{Q}}{2\alpha}-2}+m^{2}(1+t)^{-\frac{\mathcal{Q}}{2\alpha}}\right] e^{-\frac{m}{b}t}\\\left(\left\|(-\mathcal{L}_{\mathbb{H}})^{\alpha/2}u_{0}\right\|^{2}_{L^{2}\left(\mathbb{H}^{n}\right)}+\left\|u_{1}\right\|^{2}_{L^{2}\left(\mathbb{H}^{n}\right)}+\left\|u_{0}\right\|^{2}_{L^{1}\left(\mathbb{H}^{n}\right)}+\left\|u_{1}\right\|^{2}_{L^{1}\left(\mathbb{H}^{n}\right)}\right).
\end{multline}
and
\begin{equation*}
	\|\partial_{t}u(t, \cdot)\|^{2}_{L^2\left(\mathbb{H}^{n}\right)}\lesssim\left[(1+t)^{-2}+m^{2}\right] e^{-\frac{m}{b}t}\left(\left\|(-\mathcal{L}_{\mathbb{H}})^{\alpha/2}u_{0}\right\|^{2}_{L^{2}\left(\mathbb{H}^{n}\right)}+\left\|u_{1}\right\|^{2}_{L^{2}\left(\mathbb{H}^{n}\right)}\right),
\end{equation*}
respectively. Hence the above estimates will instantly give our required estimate
\eqref{deltl1l2l2final} and \eqref{deltl2l2final}. This completes the proof of Theorem 
\ref{sob:well}.
\section{Global existence}\label{sec4}
In this section, we will prove the global existence of a solution for the fractional damped wave equation with initial Cauchy data having $L^{1}\cap L^{2} $-regularity. We will also  highlight the proof of global existence with data having $L^2$-regularity only. 

 Let us recall some well-known integrals, which will be utilized to approximate the integrals over time. 
\begin{lemma}\label{intlemma0}\cite{SHANGBIN}
    Suppose that $\theta \in[0,1), a \geq 0$ and $b \geq 0$. Then, there exists $a$ constant $C=C(a, b, \theta)>0$, such that for all $t>0$ the following estimate holds:

\begin{align*}
\int_0^t(t  -\tau)^{-\theta}(1+t-\tau)^{-a}(1+\tau)^{-b} \mathrm{d} \tau 
 \leq\begin{cases}
C(1+t)^{-\min \{a+\theta, b\}} & \text {if }  \max \{a+\theta, b\}>1, \\
C(1+t)^{-\min \{a+\theta, b\}} \ln (2+t) & \text {if }  \max \{a+\theta, b\}=1, \\
C(1+t)^{1-a-\theta-b} & \text {if }  \max \{a+\theta, b\}<1.
\end{cases}
\end{align*}
\end{lemma}
 The above result also has another variation that goes back to \cite{segal} which have been utilized for plenty of linear as well as nonlinear equations. The following lemma will provide ease in our calculation, especially when different regularity is utilized for different intervals.
\begin{lemma}\label{intlemma} Let $\sigma,\beta$ be constants. 
\begin{enumerate}
\item Let $\sigma>0$ and $\beta>1$. Then
        \begin{equation*}
	\int_{0}^{t/2}(1+t-s)^{-\sigma}(1+s)^{-\beta}ds\lesssim (1+t)^{-\sigma}.
\end{equation*}
\item Let $\sigma<1$ and $\beta\in\mathbb{R}$. Then
        \begin{equation*}
	\int_{t/2}^{t}(1+t-s)^{-\sigma}(1+s)^{-\beta}ds\lesssim (1+t)^{1-\sigma-\beta}.
\end{equation*}
\end{enumerate}
\end{lemma}
    The first item is straightforward  following from the fact $(1+t-\tau)^{-\sigma}\lesssim (1+t)^{-\sigma }$ on $[0,t/2]$ and the convergence of $\int_{0}^{\infty}(1+\tau)^{-\beta}\mathrm{d}\tau$ for $\beta>1$. The second item can be found in \cite[Lemma 1.1]{takashi}. 
 
    In the following lemma, we extend the last item of the above lemma  for $\sigma=1$, which will be a key tool to estimate the time derivative $\partial_{t}u$:
    \begin{lemma}\label{intlemma2}
        Let $\sigma \in \mathbb{R}$ and $\beta>1$. Then
        \begin{equation*}
	\int_{t/2}^{t}(1+t-s)^{-1}(1+s)^{-\sigma-\beta}ds\lesssim (1+t)^{-\sigma-1}.
\end{equation*}
    \end{lemma}
    \begin{proof} Using the bound 
    \begin{equation*}
        (1+s)^{-\sigma-\beta}\leq\begin{cases}
        2^{\sigma+\beta}(1+t)^{-\sigma-\beta}&\text{if } \sigma+\beta>0,\\
            (1+t)^{-\sigma-\beta}&\text{if } \sigma+\beta\leq 0,
        \end{cases}  
    \end{equation*}
 for all $s\in[t/2,t]$, and the fact that $\frac{1}{2}(\beta-1)>0$, we have
    \begin{align*}
        \int_{t/2}^{t}(1+t-s)^{-1}(1+s)^{-\sigma-\beta}ds&\lesssim (1+t)^{-\sigma-\beta}\int_{t/2}^{t}(1+t-s)^{-1+\frac{1}{2}(\beta-1)}ds\nonumber\\&\lesssim (1+t)^{-\sigma-\beta+\frac{1}{2}(\beta-1)}\lesssim (1+t)^{-\sigma-1},
    \end{align*}
    where the last step follows from the fact that $\beta-\frac{1}{2}(\beta-1)>1$, since $\beta>1$.
    \end{proof}

\subsection{Global existence for initial data in $L^{1}\cap L^{2}$: The case $\mathbf{m=0}$}\label{sec4.1} In order to establish the global existence result stated in Theorem \ref{global:existence}, we will utilize the decay estimates developed in Theorem \ref{sob:well} for $m=0$ along with the Gagliardo-Nirenberg inequality, i.e., Theorem \ref{eq16} and the Banach fixed point theorem. 

 We can now supply the proof of the global existence of a solution for the Cauchy problem \eqref{non:Heisenberg} with initial data having $L^{1}\cap L^{2}$-regularity. 
\begin{proof}[Proof of Theorem \ref{global:existence}]
	The solution of the following homogeneous Cauchy problem 
 \begin{align}\label{hommzero:Heisenberg}	\begin{cases}	u_{tt}(t, \eta)+\left(-\mathcal{L}_{\mathbb{H}}\right)^{\alpha} u(t, \eta)+bu_{t}(t,\eta)=0,  &\eta\in \mathbb{H}^{n},t>0, \\		u(0, \eta)=u_{0}(\eta),\quad u_{t}(0, \eta)=u_{1}(\eta),   &\eta \in  \mathbb{H}^{n},	\end{cases} \end{align}
 can be expressed as
	\begin{equation*}
		u^{\operatorname{lin}}(t,\eta)=u_{0}(\eta)*_{(g)}E_{0}(t,\eta)+u_{1}(\eta)*_{(g)}E_{1}(t,\eta),
	\end{equation*}
	where $*_{(\eta)}$ denotes the convolution on $\mathbb{H}^{n}$;   $E_{0}(t,\eta)$ and $E_{1}(t,\eta)$ denote the fundamental solution to the above Cauchy problem with initial data $(u_{0},u_{1})=(\delta_{0},0)$ and $(u_{0},u_{1})=(0,\delta_{0})$, respectively. Utilizing  Duhamel's principle, the solution to the Cauchy problem \eqref{non:Heisenberg}, can be written as
	\begin{equation*}
		u(t,\eta)=u^{\operatorname{lin}}(t,\eta)+\int_{0}^{t}E_{1}(t-\tau,\eta)*_{(\eta)}|u(\tau,\eta)|^{p}\mathrm{d}\tau=u^{\operatorname{lin}}(t,\eta)+u^{\operatorname{non}}(t,\eta).
	\end{equation*}
Let's introduce a solution space $X$ defined as
$$X:=C(\mathbb{R}_{+};L^{2}(\mathbb{H}^{n})\cap H^{\alpha}(\mathbb{H}^{n}))\times C^{1}(\mathbb{R}_{+};L^{2}(\mathbb{H}^{n}))$$
associated with the norm
 \begin{equation*}
\left\|u\right\|_{X}:=\sup\limits_{t\in[0,\infty)}\left\{f_{1}(t)^{-1}\left\|u(t,\cdot)\right\|_{L^{2}(\mathbb{H}^{n})}+f_{2}(t)^{-1}\left\|u(t,\cdot)\right\|_{H^{\alpha}(\mathbb{H}^{n})}+f_{3}(t)^{-1}\left\|\partial_{t}u(t,\cdot)\right\|_{L^{2}(\mathbb{H}^{n})}\right\},
\end{equation*}
where
$$f_{1}(t)=(1+t)^{-\frac{\mathcal{Q}}{4\alpha}},\quad f_{2}(t)=(1+t)^{-\frac{\mathcal{Q}}{4\alpha}-\frac{1}{2}},\quad\text{and}\quad f_{3}(t)=(1+t)^{-\frac{\mathcal{Q}}{4\alpha}-1}.$$
 We now define the operator $\mathcal{N}:X\to X$  as follows:
\begin{equation}\label{nudef}
\mathcal{N} u(t,\eta):=u^{\operatorname{lin}}(t,\eta)+u^{\operatorname{non}}(t,\eta),\quad \text{for all }t>0,\eta\in \mathbb{H}^{n}.
\end{equation}
We will prove that the operator $\mathcal{N}$ satisfies the following inequalities:
\begin{eqnarray}
 \|\mathcal{N} u\|_{X}&\leq& L\label{ineq1}, \quad\text{for some suitable } L>0,\\
 \|\mathcal{N} u-\mathcal{N} v\|_{X} &\leq& \frac{1}{r}\|u-v\|_{X},\quad\text{for some suitable } r>1,\label{inequali2}
\end{eqnarray}
for all $u, v \in X$, where $u-v$ is the solution of the following Cauchy problem:
\begin{align}\label{uveqn}
	\begin{cases}
		(u-v)_{tt}(t, \eta)+\left(-\mathcal{L}_{\mathbb{H}}\right)^{\alpha} (u-v)(t, \eta)+b(u-v)_{t}(t,\eta)
		=|u|^{p}-|v|^{p}, & \eta \in  \mathbb{H}^{n},t>0\\
		u(0, \eta)=0,\quad u_{t}(0, \eta)=0,  & \eta \in  \mathbb{H}^{n}.
	\end{cases}
\end{align}
The inequalities \eqref{ineq1} and \eqref{inequali2}
 will allows us to prove that $\mathcal{N}$ is a contraction mapping on Banach space $X$, and consequently, using the Banach fixed point theorem, the operator $\mathcal{N}$ has a unique fixed point in $X$, whence we can obtain a unique global solution $u$ of the equation
\begin{equation*}
\mathcal{N}u=u,\quad \text{in } X.
\end{equation*}
Let us now proceed to establish the inequalities \eqref{ineq1} and \eqref{inequali2}.

\medskip

\noindent\textbf{Estimate for $\|\mathcal{N}u\|_{X}$:}
Utilizing the definition of norm $\|\cdot\|_{X}$ and the estimates in Theorem \ref{sob:well} for $m=0$, we get
\begin{equation}\label{fineqn0}
    \|u^{\operatorname{lin}}(t,\cdot)\|_{X}\leq \|(u_{0},u_{1})\|_{L^{1}(\mathbb{H}^{n})}+\|(-\mathcal{L}_{\mathbb{H}})^{\alpha/2}u_{0},u_{1})\|_{L^{2}(\mathbb{H}^{n})}=\|(u_{0},u_{1})\|_{\mathcal{B}^{\alpha}},
\end{equation}
and hence $u^{\operatorname{lin}}\in X$. Now
in order to compute the norm $\left\|u^{\operatorname{non}}\right\|_{X}$,  we need to estimate the norms 
\begin{equation*}
	\|\partial_{t}^{i}(-\mathcal{L}_{\mathbb{H}})^{j/2}u^{\operatorname{non}}(t)\|_{L^{2}(\mathbb{H}^{n})},\quad (i,j)\in\{(0,0),(0,\alpha),(1,0)\}.
\end{equation*}

\medskip

\noindent\textit{Estimate for $\|u^{\operatorname{non}}(t,\cdot)\|_{L^{2}(\mathbb{H}^{n})}$:}
Recalling the $(L^{1}\cap L^{2})-L^{2}$ estimate \eqref{l1l2l2final} for $u_{0}=0$ and $u_{1}=|u(\sigma,\cdot)|^{p}$, this gives
\begin{align}\label{unonestN}
	\|u^{\operatorname{non}}(t,\cdot)\|_{L^{2}(\mathbb{H}^{n})}&=\left\|\int_{0}^{t}E_{1}(t-\tau,\eta)*_{(\eta)}|u(\tau,\eta)|^{p}\mathrm{d}\tau\right\|_{L^{2}(\mathbb{H}^{n})}\nonumber\\
	&\leq \int_{0}^{t}(1+t-\tau)^{-\frac{\mathcal{Q}}{4\alpha}}\left(\left\||u(\tau,\cdot)|^{p}\right\|_{L^{1}(\mathbb{H}^{n})}+\left\||u(\tau,\cdot)|^{p}\right\|_{L^{2}(\mathbb{H}^{n})}\right)\mathrm{d}\tau.
\end{align}
Now utilizing the fractional Gagliardo-Nirenberg inequality along with the definition of $\|\cdot\|_{X}$, we deduce that for $\mathcal{Q}>2\alpha$ and $2\leq p\leq \frac{\mathcal{Q}}{\mathcal{Q}-2\alpha}=1+\frac{2\alpha}{\mathcal{Q}-2\alpha}$  the following inequalities hold:
\begin{align}
    \left\||u(\tau,\cdot)|^{p}\right\|_{L^{1}(\mathbb{H}^{n})}&=	\|u(\tau,\cdot)\|_{L^{p}}^{p}
\leq(1+\tau)^{-\frac{\mathcal{Q}}{2\alpha}(p-1)}\|u\|_{X}^{p}, \quad \text{and}\label{vpL1N}\\
\left\||u(\tau,\cdot)|^{p}\right\|_{L^{2}(\mathbb{H}^{n})}&=	\|u(\tau,\cdot)\|_{L^{2p}}^{p}\leq(1+\tau)^{-\frac{\mathcal{Q}}{2\alpha}\left(p-\frac{1}{2}\right)}\|u\|_{X}^{p}.\label{vpL2N}
\end{align}

Combining the above estimates with \eqref{unonestN},  we deduce that
\begin{eqnarray}\label{unonint}
	\|u^{\operatorname{non}}(t,\cdot)\|_{L^{2}(\mathbb{H}^{n})}&\lesssim& \int_{0}^{t}(1+t-\tau)^{-\frac{\mathcal{Q}}{4\alpha}}\left[(1+\tau)^{-\frac{\mathcal{Q}}{2\alpha}(p-1)}+(1+\tau)^{-\frac{\mathcal{Q}}{2\alpha}\left(p-\frac{1}{2}\right)}\right]\mathrm{d}\tau\|u\|_{X}^{p}\nonumber\\
	&\leq&2\int_{0}^{t}(1+t-\tau)^{-\frac{\mathcal{Q}}{4\alpha}}(1+\tau)^{-\frac{\mathcal{Q}}{2\alpha}(p-1)}\mathrm{d}\tau\|u\|_{X}^{p}\nonumber\\
 &\lesssim&(1+t)^{-\frac{\mathcal{Q}}{4\alpha}}\|u\|_{X}^{p},
\end{eqnarray}
where last step follows from Lemma \ref{intlemma0} for $\theta=0$, since 
\begin{equation*}
    \min\left\{\frac{\mathcal{Q}}{4\alpha},\frac{\mathcal{Q}}{2\alpha}(p-1)\right\}=\frac{\mathcal{Q}}{4\alpha} ~~\text{ and } ~~\max\left\{\frac{\mathcal{Q}}{4\alpha},\frac{\mathcal{Q}}{2\alpha}(p-1)\right\}=\frac{\mathcal{Q}}{2\alpha}(p-1)>1,
\end{equation*}
 for $p\geq2$ and $\mathcal{Q}>2\alpha$. Thus, we have
\begin{equation}\label{fineqn1}
	\|u^{\operatorname{non}}(t,\cdot)\|_{L^{2}(\mathbb{H}^{n})}\lesssim(1+t)^{-\frac{\mathcal{Q}}{4\alpha}}\|u\|_{X}^{p}.
\end{equation}

\medskip

\noindent\textit{Estimate for $\|(-\mathcal{L}_{\mathbb{H}})^{\alpha/2}u^{\operatorname{non}}(t,\cdot)\|_{L^{2}(\mathbb{H}^{n})}$:}
\sloppy In order to compute the estimate for $\|(-\mathcal{L}_{\mathbb{H}})^{\alpha/2}u^{\operatorname{non}}(t,\cdot)\|_{L^{2}(\mathbb{H}^{n})}$, we utilize the $(L^{1}\cap L^{2})-L^{2}$  estimate \eqref{delul1l2final} if $\tau\in[0,t/2]$ and the $L^{2}-L^{2}$ estimate \eqref{delul2l2final} for $\tau\in[t/2,t]$. Similar to the previous instance, using the estimates \eqref{vpL1N} and \eqref{vpL2N},  we obtain
\begin{align}\label{delualpha}
		\left\|(-\mathcal{L}_{\mathbb{H}})^{\alpha/2}u^{\operatorname{non}}(t,\cdot)\right\|_{L^{2}(\mathbb{H}^{n})}&\lesssim \int_{0}^{t/2}(1+t-\tau)^{-\frac{\mathcal{Q}}{4\alpha}-\frac{1}{2}}\left(\left\||u(\tau,\cdot)|^{p}\right\|_{L^{1}(\mathbb{H}^{n})}+\left\||u(\tau,\cdot)|^{p}\right\|_{L^{2}(\mathbb{H}^{n})}\right)\mathrm{d}\tau\nonumber\\
	&+\int_{t/2}^{t}(1+t-\tau)^{-\frac{1}{2}}\left\||u(\tau,\cdot)|^{p}\right\|_{L^{2}(\mathbb{H}^{n})}\mathrm{d}\tau\nonumber\\
	&\lesssim \int_{0}^{t/2}(1+t-\tau)^{-\frac{\mathcal{Q}}{4\alpha}-\frac{1}{2}}(1+\tau)^{-\frac{\mathcal{Q}}{2\alpha}(p-1)}\mathrm{d}\tau\nonumber\|u\|_{X}^{p}\\
	&+\int_{t/2}^{t}(1+t-\tau)^{-\frac{1}{2}}(1+\tau)^{-\frac{\mathcal{Q}}{2\alpha}\left(p-\frac{1}{2}\right)}\mathrm{d}\tau\|u\|_{X}^{p}.
\end{align}
Since $\frac{\mathcal{Q}}{2\alpha}(p-1)>1$, therefore utilizing the Lemma \ref{intlemma}, we get
\begin{equation*}
	\left\|(-\mathcal{L}_{\mathbb{H}})^{\alpha/2}u^{\operatorname{non}}(t,\cdot)\right\|_{L^{2}(\mathbb{H}^{n})}
		\lesssim\left[(1+t)^{-\frac{\mathcal{Q}}{4\alpha}-\frac{1}{2}}+(1+t)^{1-\frac{1}{2}-\frac{\mathcal{Q}}{2\alpha}\left(p-\frac{1}{2}\right)}\right]\|u\|_{X}^{p}
	\lesssim 2(1+t)^{-\frac{\mathcal{Q}}{4\alpha}-\frac{1}{2}}\|u\|_{X}^{p}.
\end{equation*}
 Thus, we have
\begin{equation}\label{fineqn2}
	\left\|u^{\operatorname{non}}(t,\cdot)\right\|_{H^{\alpha}(\mathbb{H}^{n})}\lesssim (1+t)^{-\frac{\mathcal{Q}}{4\alpha}-\frac{1}{2}}\|u\|_{X}^{p}.
\end{equation}

\medskip

\noindent\textit{Estimate for $\|\partial_{t}u^{\operatorname{non}}(t,\cdot)\|_{L^{2}(\mathbb{H}^{n})}$:}
In this instance, we utilize the $(L^{1}\cap L^{2})-L^{2}$ estimate \eqref{deltl1l2l2final} if $\tau\in[0,t/2]$ and the $ L^{2}-L^{2}$ estimate \eqref{deltl2l2final} for $\tau\in[t/2,t]$. Similar to the previous case, we get
\begin{align}\label{deltu}
		\left\|\partial_{t}u^{\operatorname{non}}(t,\cdot)\right\|_{L^{2}(\mathbb{H}^{n})}&\lesssim \int_{0}^{t/2}(1+t-\tau)^{-\frac{\mathcal{Q}}{4\alpha}-1}\left(\left\||u(\tau,\cdot)|^{p}\right\|_{L^{1}(\mathbb{H}^{n})}+\left\||u(\tau,\cdot)|^{p}\right\|_{L^{2}(\mathbb{H}^{n})}\right)\mathrm{d}\tau\nonumber\\
	&+\int_{t/2}^{t}(1+t-\tau)^{-1}\left\||u(\tau,\cdot)|^{p}\right\|_{L^{2}(\mathbb{H}^{n})}\mathrm{d}\tau\nonumber\\
	&\lesssim \int_{0}^{t/2}(1+t-\tau)^{-\frac{\mathcal{Q}}{4\alpha}-1}(1+\tau)^{-\frac{\mathcal{Q}}{2\alpha}(p-1)}\mathrm{d}\tau\nonumber\|u\|_{X}^{p}\\
	&+\int_{t/2}^{t}(1+t-\tau)^{-1}(1+\tau)^{-\frac{\mathcal{Q}}{2\alpha}\left(p-\frac{1}{2}\right)}\mathrm{d}\tau\|u\|_{X}^{p}.
\end{align}
Utilizing the Lemma \ref{intlemma} for the first integral and Lemma \ref{intlemma2} for the second integral, since $\frac{\mathcal{Q}}{2\alpha}(p-1)>1$ and $\frac{\mathcal{Q}}{2\alpha}\left(p-\frac{1}{2}\right)=-\frac{\mathcal{Q}}{4\alpha}-\frac{\mathcal{Q}}{2\alpha}\left(p-1\right)$,    we deduce that
\begin{equation}\label{fineqn3}
\left\|\partial_{t}u^{\operatorname{non}}(t,\cdot)\right\|_{L^{2}(\mathbb{H}^{n})}\lesssim(1+t)^{-\frac{\mathcal{Q}}{4\alpha}-1}\|u\|_{X}^{p}.
\end{equation}
Combining the estimates \eqref{fineqn0}, \eqref{fineqn1}, \eqref{fineqn2}, and \eqref{fineqn3}, we get
\begin{equation}\label{nufinal}
    \|\mathcal{N}u\|_{X}\leq A\|(u_{0},u_{1})\|_{\mathcal{B}^{\alpha}}+A\|u\|^{p}_{X}
\end{equation}
for some positive constant $A$.

\medskip

\noindent\textbf{Estimate for $\|\mathcal{N}u-\mathcal{N}v\|_{X}$:} For 
$u, v\in X$, the definition \eqref{nudef} gives
\begin{equation*}
	\mathcal{N}u(t,\eta)-\mathcal{N}v(t,\eta)=u^{\mathrm{lin}}(t,\eta)-v^{\mathrm{lin}}(t,\eta)+ u^{\mathrm{non}}(t,\eta)-v^{\mathrm{non}}(t,\eta).
\end{equation*}
It is easy to verify that $\|u^{\mathrm{lin}}-v^{\mathrm{lin}}\|_{X}=0$.
Now
in order to compute the norm $\left\|u^{\operatorname{non}}-v^{\operatorname{non}}\right\|_{X}$,  we need to estimate the norms 
\begin{equation*}
	\|\partial_{t}^{i}(-\mathcal{L}_{\mathbb{H}})^{j/2}(u^{\operatorname{non}}-v^{\operatorname{non}})(t)\|_{L^{2}(\mathbb{H}^{n})},\quad (i,j)\in\{(0,0),(0,\alpha),(1,0)\}.
\end{equation*}
Utilizing the estimate exactly similar to the case $\left\|\partial_{t}^{i}(-\mathcal{L}_{\mathbb{H}})^{j/2}u^{\operatorname{non}}(t,\cdot)\right\|_{L^{2}(\mathbb{H}^{n})}$, we get
\begin{equation*}
	\left\|\left(u^{\operatorname{non}}-v^{\operatorname{non}}\right)(t)\right\|_{L^{2}(\mathbb{H}^{n})} \lesssim \int_0^t(1+t-\tau)^{-\frac{\mathcal{Q}}{4\alpha}}\left\||u(\tau)|^p-|v(\tau)|^p\right\|_{L^{1}(\mathbb{H}^{n}) \cap L^{2}(\mathbb{H}^{n})} \mathrm{d} \tau,
\end{equation*}
\begin{multline*}
	\left\|(-\mathcal{L}_{\mathbb{H}})^{\alpha/2}\left(u^{\operatorname{non}}-v^{\operatorname{non}}\right)(t)\right\|_{L^{2}(\mathbb{H}^{n})} \lesssim \int_0^{t / 2}(1+t-\tau)^{-\frac{\mathcal{Q}}{4\alpha}-\frac{1}{2}}\left\||u(\tau)|^p-|v(\tau)|^p\right\|_{L^1(\mathbb{H}^{n}) \cap L^2(\mathbb{H}^{n})}  \mathrm{d} \tau \\
	\quad+\int_{t / 2}^t(1+t-\tau)^{-\frac{1}{2}}\left\||u(\tau)|^p-|v(\tau)|^p\right\|_{L^{2}(\mathbb{H}^{n})} \mathrm{d} \tau,\quad \text{and}
\end{multline*}
\begin{multline*}
	\left\|\partial_{t}\left(u^{\operatorname{non}}-v^{\operatorname{non}}\right)(t)\right\|_{L^{2}(\mathbb{H}^{n})} \lesssim \int_0^{t / 2}(1+t-\tau)^{-\frac{\mathcal{Q}}{4\alpha}-1}\left\||u(\tau)|^p-|v(\tau)|^p\right\|_{L^1(\mathbb{H}^{n}) \cap L^2(\mathbb{H}^{n})}  \mathrm{d} \tau \\
	\quad+\int_{t / 2}^t(1+t-\tau)^{-1}\left\||u(\tau)|^p-|v(\tau)|^p\right\|_{L^{2}(\mathbb{H}^{n})} \mathrm{d} \tau.
\end{multline*}
 Using H\"{o}lder's inequality, we deduce that
\begin{align}
	\left\||u(\tau)|^p-|v(\tau)|^p\right\|_{L^{1}(\mathbb{H}^{n})} &\lesssim\|u(\tau)-v(\tau)\|_{L^p(\mathbb{H}^{n})}\left(\|u(\tau)\|_{L^p(\mathbb{H}^{n})}^{p-1}+\|v(\tau)\|_{L^p(\mathbb{H}^{n})}^{p-1}\right)\text{ and } \\
	\left\||u(\tau)|^p-|v(\tau)|^p\right\|_{L^{2}(\mathbb{H}^{n})} &\lesssim\|u(\tau)-v(\tau)\|_{L^{2 p}(\mathbb{H}^{n})}\left(\|u(\tau)\|_{L^{2 p}(\mathbb{H}^{n})}^{p-1}+\|v(\tau)\|_{L^{2 p}(\mathbb{H}^{n})}^{p-1}\right).
\end{align}
Now utilizing the Gagliardo-Nirenberg inequality and the estimates from Theorem \ref{sob:well}, we deduce that
\begin{align*}
	&\left\||u(\tau)|^p-|v(\tau)|^p\right\|_{L^{1}(\mathbb{H}^{n})}\lesssim (1+\tau)^{-\frac{\mathcal{Q}}{2\alpha }(p-1)}\|u-v\|_{X}\left(\|u\|_{X}^{p-1}+\|v\|^{p-1}_{X}\right)\text{ and }\\
	&\left\||u(\tau)|^p-|v(\tau)|^p\right\|_{L^{2}(\mathbb{H}^{n})}\lesssim (1+\tau)^{-\frac{\mathcal{Q}}{2\alpha }\left(p-\frac{1}{2}\right)}\|u-v\|_{X}\left(\|u\|_{X}^{p-1}+\|v\|^{p-1}_{X}\right).
\end{align*}
Now similar to the estimates \eqref{fineqn1}, \eqref{fineqn2}, and \eqref{fineqn3}, we get
\begin{align*}
    \left\|\left(u^{\operatorname{non}}-v^{\operatorname{non}}\right)(t)\right\|_{L^{2}(\mathbb{H}^{n})} &\lesssim (1+t)^{-\frac{\mathcal{Q}}{4\alpha}}\|u-v\|_{X}\left(\|u\|_{X}^{p-1}+\|v\|^{p-1}_{X}\right),\\
        \left\|(-\mathcal{L}_{\mathbb{H}})^{\alpha/2}\left(u^{\operatorname{non}}-v^{\operatorname{non}}\right)(t)\right\|_{L^{2}(\mathbb{H}^{n})} &\lesssim (1+t)^{-\frac{\mathcal{Q}}{4\alpha}-\frac{1}{2}}\|u-v\|_{X}\left(\|u\|_{X}^{p-1}+\|v\|^{p-1}_{X}\right),\quad \text{and}\\
            \left\|\partial_{t}\left(u^{\operatorname{non}}-v^{\operatorname{non}}\right)(t)\right\|_{L^{2}(\mathbb{H}^{n})} &\lesssim (1+t)^{-\frac{\mathcal{Q}}{4\alpha}-1}\|u-v\|_{X}\left(\|u\|_{X}^{p-1}+\|v\|^{p-1}_{X}\right).
\end{align*}
Thus, combining the above estimate with $\|u^{\mathrm{lin}}-v^{\mathrm{lin}}\|_{X}=0$, we have
\begin{equation}\label{numinusnv}
    \|\mathcal{N}u-\mathcal{N}v\|_{X}\leq B\|u-v\|_{X}\left(\|u\|_{X}^{p-1}+\|v\|^{p-1}_{X}\right),
\end{equation}
for some positive constant $B>0$.

If we set $R:=rA\|(u_{0},u_{1})\|_{\mathcal{B}^{\alpha}}$ with sufficiently small $\|(u_{0},u_{1})\|_{\mathcal{B}^{\alpha}}\leq \varepsilon$ and for some $r>1$ such that
\begin{equation*}
    AR^p<\frac{R}{r} \quad \text{and} \quad 2BR^{p-1}<\frac{1}{r}.
\end{equation*}
Then     \eqref{nufinal}   and \eqref{numinusnv} reduce to  
\begin{align}\label{Final banach1}
	\|\mathcal{N} u\|_{X} \leq  \frac{2R}{r} \quad\text{and } \quad \|\mathcal{N}u-\mathcal{N}v\|_{X} \leq  \frac{1}{r}\|u-v\|_{X},	
\end{align}  
respectively for all $u, v\in \mathscr{B}_{R}:=\{u\in X: \|u\|_{X}\leq R\}$. Since $\mathscr{B}_{R}$ is a closed ball of Banach space $X$, hence itself a Banach space. Therefore using the Banach fixed theorem, we can conclude the proof of Theorem \ref{global:existence}.
\end{proof}
\subsection{Global existence for initial data in $ L^{2}$ only: The case $\mathbf{m=0}$}\label{sec4.2}
The proof of global existence for initial data in $ L^{2}$ only can be obtained from above proof with the following minor modifications:
\begin{itemize}
    \item Using the same solution space $X$ with the following modified norm
    \begin{equation*}
\left\|u\right\|_{X}:=\sup\limits_{t\in[0,\infty)}\left\{\left\|u(t,\cdot)\right\|_{L^{2}(\mathbb{H}^{n})}+(1+t)^{\frac{1}{2}}\left\|(-\mathcal{L}_{\mathbb{H}})^{\alpha/2}u(t,\cdot)\right\|_{L^{2}(\mathbb{H}^{n})}+(1+t)\left\|\partial_{t}u(t,\cdot)\right\|_{L^{2}(\mathbb{H}^{n})}\right\}.
\end{equation*}
\item The Gagliardo Nirenberg inequality gives
    \begin{equation*}
\left\||u(\tau,\cdot)|^{p}\right\|_{L^{2}(\mathbb{H}^{n})}\lesssim(1+\tau)^{-\frac{\mathcal{Q}}{4\alpha}\left(p-1\right)}\|u\|_{X}^{p},\quad 1\leq p<1+\frac{2\alpha}{\mathcal{Q}-2\alpha}.
\end{equation*}
\item Using $L^{2}-L^{2}$ estimates from Theorem \ref{sob:well}, we obtain
\begin{equation}\label{m=0:int}
    \|\partial_{t}^{i}(-\mathcal{L}_{\mathbb{H}})^{j/2}u^{\operatorname{non}}(t)\|_{L^{2}(\mathbb{H}^{n})}\lesssim \int_{0}^{t}(1+t-\tau)^{-\frac{j}{2\alpha}-i}(1+\tau)^{-\frac{\mathcal{Q}}{4\alpha}\left(p-1\right)}\mathrm{d}\tau\|u\|_{X}^{p},
\end{equation}
for $(i,j)\in\{(0,0),(0,\alpha),(1,0)\}$.
\item Under the hypothesis $p>1+\frac{4\alpha}{\mathcal{Q}}$, we have
\begin{equation*}
    \min\left\{\frac{j}{2\alpha}+i,\frac{\mathcal{Q}}{4\alpha}(p-1)\right\}=\frac{j}{2\alpha}+i ~~\& ~~\left\{\frac{j}{2\alpha}+i,\frac{\mathcal{Q}}{2\alpha}(p-1)\right\}=\frac{\mathcal{Q}}{2\alpha}(p-1)>1.
\end{equation*}
\item Utilizing  Lemma \ref{intlemma0} for $\theta=0$, the estimate \eqref{m=0:int} gives
\begin{equation}\label{m=0:2nd}
    \|\partial_{t}^{i}(-\mathcal{L}_{\mathbb{H}})^{j/2}u^{\operatorname{non}}(t)\|_{L^{2}(\mathbb{H}^{n})}\lesssim (1+t)^{-\frac{j}{2\alpha}-i}\|u\|_{X}^{p}.
\end{equation}
\end{itemize}

\section{Application}\label{sec5}
 This section is devoted to investigate the well-posedness of the weakly coupled system on the Heisenberg group $\HH$ as an application of our linear estimates addressed in Theorem \ref{sob:well}. 
Let us  consider the Cauchy problem for a weakly coupled system with  two  
  semi-linear  fractional  damped wave equations with positive mass term on $\HH$, namely 
\begin{align}\label{main:Coupled:Heisenberg}
	\begin{cases}
		u_{tt}(t, \eta)+(-\mathcal{L}_{\mathbb{H}})^{\alpha} u(t, \eta)+bu_{t}(t,\eta)+mu(t,\eta)
		=|v|^{p},  &\eta\in \mathbb{H}^{n},t>0, \\
		v_{tt}(t, \eta)+(-\mathcal{L}_{\mathbb{H}})^{\alpha}v(t, \eta)+bv_{t}(t,\eta)+mv(t,\eta)
		=|u|^{q}, &\eta\in \mathbb{H}^{n},t>0, \\
		u(0, \eta)=u_{0}(\eta),\quad u_{t}(0, \eta)=u_{1}(\eta),  &\eta\in \mathbb{H}^{n},\\
		v(0,\eta)=v_{0}(\eta)	,\quad v_{t}(0,\eta)=v_{1}(\eta),&\eta\in \mathbb{H}^{n},
	\end{cases}
\end{align}
for $p, q>1$.  
In the Eucledian framework, such type of problem is well-known as a part of the so-called Nakao’s problem proposed by Professor Mitsuhiro Nakao to determine the critical curve between the exponents $p$ and $q$ of the power nonlinearities,  and explored by several researchers for different system with different type of nonlinearities,  see  \cite{AKT2000,GTY2006, NW2014, abbico2017, CR2021, PT2023}  and references therein.  However, to the best of our knowledge, a weakly coupled system in the non-Eucledian framework has not been considered in the literature so far, even for the
power-type nonlinearities. Here, in the following result, we examine the well-posedness result for the weakly coupled system \eqref{main:Coupled:Heisenberg} on the Heisenberg group $\HH$ as an application of our linear estimates considered in Theorem  \ref{sob:well}. 
\begin{theorem}\label{global:Application}Let $\alpha>0$ and $\mathcal{Q}>2\alpha$. Let $1<p,q \leq 1+ \frac{2\alpha}{\mathcal{Q}-2\alpha}$,
 then, there exists a sufficiently small $\varepsilon>0$ such that for any initial data $$ 
		((u_{0},u_{1}),(v_{0},v_{1}))\in \mathcal{D}^{\alpha}\times \mathcal{D}^{\alpha},\quad \text{where }\mathcal{D}^{\alpha}:= H^{\alpha}(\mathbb{H}^n) \times L^{2}\left(\mathbb{H}^{n}\right),$$ 	satisfying
	\begin{equation*}
		\|(u_{0},u_{1})\|_{\mathcal{D}^{\alpha}}+\|(v_{0},v_{1})\|_{\mathcal{D}^{\alpha}}<\varepsilon, 
	\end{equation*}
	there is a unique global solution $$(u,v) \in C\left(\mathbb{R}_{+};H^\alpha\left(\mathbb{H}^{n}\right)\right) \cap C^1\left(\mathbb{R}_{+}; L^2\left(\mathbb{H}^{n}\right)\right)\times C\left(\mathbb{R}_{+};H^\alpha\left(\mathbb{H}^{n}\right)\right) \cap C^1\left(\mathbb{R}_{+}; L^2\left(\mathbb{H}^{n}\right)\right)$$ to the Cauchy problem \eqref{main:Coupled:Heisenberg} with $b>0,m>0$ satisfying $b^{2}>4m$.  
	Moreover, the following estimates hold:
 \begin{equation*}
     \|\partial_{t}^{i}(-\mathcal{L}_{\mathbb{H}})^{j/2}(u,v)(t)\|_{L^2\left(\mathbb{H}^{n}\right)}\lesssim \left[(1+\tau)^{-\frac{j}{2\alpha}-i}+mi\right]e^{-\frac{m}{2b}t}\left(\|(u_{0},u_{1})\|_{\mathcal{D}^{\alpha}}+\|(v_{0},v_{1})\|_{\mathcal{D}^{\alpha}}\right),
 \end{equation*}
 for $(i,j)\in\{(0,0),(0,\alpha),(1,0)\}$ and for all $t\geq 0$.
\end{theorem}

 The following lemma will be useful in the proof of the  global  well-posedness of the weakly coupled system. 
\begin{lemma}\cite{abbico2017}\label{intlemma3}
Let $c>0$ and $\sigma\in\mathbb{R}$. Then
        \begin{equation*}
	\int_{0}^{t}e^{-c(t-s)}(1+s)^{-\sigma}ds\lesssim (1+t)^{-\sigma}.
\end{equation*}    
\end{lemma}
One must be very careful while handling such integral  because slight modifications may lose the decay; e.g., the above integral has sufficient polynomial decay, but  $\int_{0}^{t}e^{-c(t-s)}(1+t-s)^{-\sigma}ds$ does not have any decay. We are now in position to supply the proof of Theorem  \ref{global:Application}.
\begin{proof}[Proof of Theorem \ref{global:Application}]
    Consider the Cauchy problems:
\begin{align}\label{couple:Eq1}
	\begin{cases}
		u_{tt}(t, \eta)+(-\mathcal{L}_{\mathbb{H}})^{\alpha} u(t, \eta)+bu_{t}(t,\eta)+mu(t,\eta)
		=|v|^{p},  &\eta\in \mathbb{H}^{n},t>0, \\
		u(0, \eta)=u_{0}(\eta),\quad u_{t}(0, \eta)=u_{1}(\eta),  &\eta\in \mathbb{H}^{n},
	\end{cases}
\end{align}
and 
\begin{align}\label{couple:Eq2}
	\begin{cases}
		v_{tt}(t, \eta)+(-\mathcal{L}_{\mathbb{H}})^{\alpha}v(t, \eta)+bv_{t}(t,\eta)+mv(t,\eta)
		=|u|^{q}, &\eta\in \mathbb{H}^{n},t>0, \\
		v(0,\eta)=v_{0}(\eta)	,\quad v_{t}(0,\eta)=v_{1}(\eta),&\eta\in \mathbb{H}^{n}.
	\end{cases}
\end{align}
The solution of the corresponding homogeneous Cauchy problems 
can be expressed as
\begin{eqnarray}
    u^{\operatorname{lin}}(t,\eta)&=&u_{0}(\eta)*_{(\eta)}K_{0}(t,\eta)+u_{1}(\eta)*_{(\eta)}K_{1}(t,\eta),\\
    	v^{\operatorname{lin}}(t,\eta)&=&v_{0}(\eta)*_{(\eta)}K_{0}(t,\eta)+v_{1}(\eta)*_{(\eta)}K_{1}(t,\eta),
\end{eqnarray} 
where  $K_{0}(t,\eta)$ and $K_{1}(t,\eta)$  denotes the fundamental solution to the homogeneous Cauchy problem  with Cauchy initial data $(u_{0},u_{1})=(\delta_{0},0)$ and $(u_{0},u_{1})=(0,\delta_{0})$, respectively. 

 Utilizing the Duhamel's principle, the solution to the Cauchy problem \eqref{couple:Eq1} and \eqref{couple:Eq2} can be written as
\begin{align}
	u(t,\eta)&=u^{\operatorname{lin}}(t,\eta)+\int_{0}^{t}K_{1}(t-\tau,\eta)*_{(\eta)}|v(\tau,\eta)|^{p}\mathrm{d}\tau=u^{\operatorname{lin}}(t,\eta)+u^{\operatorname{non}}(t,\eta),\nonumber\\
    v(t,\eta)&=v^{\operatorname{lin}}(t,\eta)+\int_{0}^{t}K_{1}(t-\tau,\eta)*_{(\eta)}|u(\tau,\eta)|^{q}\mathrm{d}\tau=v^{\operatorname{lin}}(t,\eta)+v^{\operatorname{non}}(t,\eta).\nonumber
\end{align}
Let us introduce a solution space $Z$ defined as
\begin{equation*}
	Z:=C\left(\mathbb{R}_{+}; H^{\alpha}\left(\mathbb{H}_n\right)\right) \cap C^1\left(\mathbb{R}_{+}; L^2\left(\mathbb{H}_n\right)\right)\times C\left(\mathbb{R}_{+}; H^\alpha\left(\mathbb{H}_n\right)\right) \cap C^1\left(\mathbb{R}_{+}; L^2\left(\mathbb{H}_n\right)\right),
\end{equation*}
associated with the norm
\begin{multline}\label{normdef:couple}
\|(u, v)\|_{Z}:=\\\sup _{t>0}\left\{f_1(\tau)^{-1}\|u(\tau, \cdot)\|_{L^{2}(\mathbb{H}^{n})}+f_2(\tau)^{-1}\left\|(-\mathcal{L})^{\alpha/2} u(\tau, \cdot)\right\|_{L^{2}(\mathbb{H}^{n})}+f_3(\tau)^{-1}\left\|u_t(\tau, \cdot)\right\|_{L^{2}(\mathbb{H}^{n})}\right. \\
\left.+f_{1}(\tau)^{-1}\|v(\tau, \cdot)\|_{L^{2}(\mathbb{H}^{n})}+f_2(\tau)^{-1}\left\|(-\mathcal{L})^{\alpha/2} v(\tau, \cdot)\right\|_{L^{2}(\mathbb{H}^{n})}+f_3(\tau)^{-1}\left\|v_t(\tau, \cdot)\right\|_{L^{2}(\mathbb{H}^{n})}\right\}
\end{multline}
where
\begin{equation*}
    f_{1}(\tau)=e^{-\frac{m}{2b}t},\quad f_{2}(\tau)=(1+\tau)^{-\frac{1}{2}}e^{-\frac{m}{2b}t},~~\text{and }~f_{3}(\tau)=\left[(1+\tau)^{-1}+m\right]e^{-\frac{m}{2b}t}.
\end{equation*}
We now introduce a operator $\Psi:Z\to Z $ defined as:
\begin{equation*}
	\Psi[u,v](t,\eta):=(u^{\operatorname{lin}}(t,\eta),v^{\operatorname{lin}}(t,\eta))+(u^{\operatorname{non}}(t,\eta),v^{\operatorname{non}}(t,\eta)).
\end{equation*}
From the spirit of proof of Theorem \ref{global:existence}, it is enough to prove that the operator $\Psi$ satisfies the following inequalities
\begin{equation}\label{ineq1:couple}
\|\Psi[u,v]\|_{Z}\lesssim \|(u_{0},u_{1})\|_{L^{2}(\mathbb{H}^{n})}+\|(v_{0},v_{1})\|_{L^{2}(\mathbb{H}^{n})}+\|(u,v)\|^{p}_{Z}+\|(u,v)\|^{q}_{Z},~~\text{and}
\end{equation}
\begin{multline}\label{ineq2:couple}
	\|\Psi[u,v]-\Psi[\overline{u},\overline{v}]\|_{Z}\lesssim \|(u,v)-(\overline{u},\overline{v})\|_{Z}\left(\|(u,v)\|_{Z}^{p-1}+\|(\overline{u},\overline{v})\|^{p-1}_{Z}\right.\\\left.+\|(u,v)\|_{Z}^{q-1}+\|(\overline{u},\overline{v})\|^{q-1}_{Z}\right).
\end{multline}
Afterward, similar to the proof of Theorem \ref{global:existence}, we can utilize the Banach fixed point theorem to determine the global (in time) existence results for small data.\\

\medskip

\noindent\textbf{Estimate for $\|\Psi[u,v]\|_{Z}$:}
Utilizing the definition of norm $\|(\cdot,\cdot)\|_{Z}$ and Theorem \ref{sob:well}, we get
\begin{equation}\label{ulinnorm:couple}
	\left\|(u^{\operatorname{lin}},v^{\operatorname{lin}})\right\|_{Z}\lesssim \|(u_{0},u_{1})\|_{L^{2}(\mathbb{H}^{n})}+\|(v_{0},v_{1})\|_{L^{2}(\mathbb{H}^{n})},
\end{equation}
this gives $(u^{\operatorname{lin}},v^{\operatorname{lin}})\in Z$.
 In order to compute the norm $\left\|(u^{\operatorname{non}},v^{\operatorname{non}})\right\|_{Z}$,  we need to estimate the norms 
\begin{equation*}
	\|\partial_{t}^{i}(-\mathcal{L}_{\mathbb{H}})^{j/2}u^{\operatorname{non}}(t)\|_{L^{2}(\mathbb{H}^{n})} \quad \text{and }	\|\partial_{t}^{i}(-\mathcal{L}_{\mathbb{H}})^{j/2}v^{\operatorname{non}}(t)\|_{L^{2}(\mathbb{H}^{n})},\quad (i,j)\in\{(0,0),(0,\alpha),(1,0)\}.
\end{equation*}
Recalling the $L^{2}-L^{2}$ estimate from Theorem \ref{sob:well}  for $(u_{0},u_{1})=(0,|v(\sigma,\cdot)|^{p})$ and $(v_{0},v_{1})=(0,|u(\sigma,\cdot)|^{q})$, similar to the estimate \eqref{unonestN}, we get
\begin{align}
	\|\partial_{t}^{i}(-\mathcal{L}_{\mathbb{H}})^{j/2}u^{\operatorname{non}}(t,\cdot)\|_{L^{2}(\mathbb{H}^{n})}
	&\lesssim \int_{0}^{t}\left[(1+t-\tau)^{-\frac{j}{2\alpha}-i}+mi\right]e^{-\frac{m}{2b}(t-\tau)}\left\||v(\tau,\cdot)|^{p}\right\|_{L^{2}(\mathbb{H}^{n})}\mathrm{d}\tau\label{unon:couple}\\
 	\|\partial_{t}^{i}(-\mathcal{L}_{\mathbb{H}})^{j/2}v^{\operatorname{non}}(t,\cdot)\|_{L^{2}(\mathbb{H}^{n})}
	&\lesssim \int_{0}^{t}\left[(1+t-\tau)^{-\frac{j}{2\alpha}-i}+mi\right]e^{-\frac{m}{2b}(t-\tau)}\left\||u(\tau,\cdot)|^{q}\right\|_{L^{2}(\mathbb{H}^{n})}\mathrm{d}\tau\label{vnon:couple}.
\end{align}
Now utilizing the fractional Gagliardo-Nirenberg inequality along with the definition of $\|\cdot\|_{X}$, we get
\begin{align*}
\left\||v(\tau,\cdot)|^{p}\right\|_{L^{2}(\mathbb{H}^{n})}&\lesssim(1+\tau)^{-\frac{\mathcal{Q}}{4\alpha}\left(p-1\right)}e^{-\frac{m}{2b}p\tau}\|(u,v)\|_{Z}^{p}, \quad\text{and}\\
\left\||u(\tau,\cdot)|^{q}\right\|_{L^{2}(\mathbb{H}^{n})}&\lesssim(1+\tau)^{-\frac{\mathcal{Q}}{4\alpha}\left(q-1\right)}e^{-\frac{m}{2b}q\tau}\|(u,v)\|_{Z}^{q},
\end{align*}
for $p>1$, where $\mathcal{Q}>2\alpha,\alpha>0$, and $1\leq p,q\leq \frac{\mathcal{Q}}{\mathcal{Q}-2\alpha}=1+\frac{2\alpha}{\mathcal{Q}-2\alpha},$ 
$(i,j)\in\{(0,0),(0,\alpha),(1,0)\}$.
The estimates \eqref{unon:couple} and \eqref{vnon:couple} along with the above estimates takes the form
\begin{align}\label{u1equation:couple}
	&\|\partial_{t}^{i}(-\mathcal{L}_{\mathbb{H}})^{j/2}u^{\operatorname{non}}(t,\cdot)\|_{L^{2}(\mathbb{H}^{n})}\nonumber\\&\lesssim \int_{0}^{t}\left[(1+t-\tau)^{-\frac{j}{2\alpha}-i}+mi\right]e^{-\frac{m}{2b}(t-\tau)}(1+\tau)^{-\frac{\mathcal{Q}}{4\alpha}(p-1)}e^{-\frac{m}{2b}p\tau}\mathrm{d}\tau\|(u,v)\|_{Z}^{p}\nonumber\\
 &\lesssim e^{-\frac{m}{2b}t}\int_{0}^{t}e^{-\frac{m}{2b}(p-1)\tau}\left[(1+t-\tau)^{-\frac{j}{2\alpha}-i}+mi\right](1+\tau)^{-\frac{\mathcal{Q}}{4\alpha}(p-1)}\mathrm{d}\tau\|(u,v)\|_{Z}^{p}\nonumber\\
  &\lesssim e^{-\frac{m}{2b}t}\int_{0}^{t}e^{-\frac{m}{2b}(p-1)(t-\tau)}\left[(1+\tau)^{-\frac{j}{2\alpha}-i}+mi\right](1+t-\tau)^{-\frac{\mathcal{Q}}{4\alpha}(p-1)}\mathrm{d}\tau\|(u,v)\|_{Z}^{p}\nonumber\\
  &\lesssim e^{-\frac{m}{2b}t}\int_{0}^{t}e^{-\frac{m}{2b}(p-1)(t-\tau)}\left[(1+\tau)^{-\frac{j}{2\alpha}-i}+mi\right]\mathrm{d}\tau\|(u,v)\|_{Z}^{p}\quad (\because \frac{\mathcal{Q}}{4\alpha}(p-1)\geq 0)\nonumber\\
 &\lesssim e^{-\frac{m}{2b}t}\left[(1+\tau)^{-\frac{j}{2\alpha}-i}+mi\right]\|(u,v)\|_{Z}^{p},
\end{align}
where last step follows from Lemma \ref{intlemma3}
for $p>1$. Similarly, we can obtain
\begin{equation}\label{v1equation:couple}
  \|\partial_{t}^{i}(-\mathcal{L}_{\mathbb{H}})^{j/2}v^{\operatorname{non}}(t,\cdot)\|_{L^{2}(\mathbb{H}^{n})}
 \lesssim \left[(1+\tau)^{-\frac{j}{2\alpha}-i}+mi\right]e^{-\frac{m}{2b}t}\|(u,v)\|_{Z}^{q}, \quad q>1,
\end{equation}
for $(i,j)\in\{(0,0),(0,\alpha),(1,0)\}.$

Combining the estimates \eqref{ulinnorm:couple}, \eqref{u1equation:couple}, and \eqref{v1equation:couple}, we 
 can obtain our required inequality \eqref{ineq1:couple}. 

Furthermore, using similar technique to estimate  $\|\partial_{t}^{i}(-\mathcal{L}_{\mathbb{H}})^{j/2}u^{\operatorname{non}}(t)\|_{L^{2}(\mathbb{H}^{n})}$ and $ \|\partial_{t}^{i}(-\mathcal{L}_{\mathbb{H}})^{j/2}v^{\operatorname{non}}(t)\|_{L^{2}(\mathbb{H}^{n})}$, we can estimate 
\begin{equation*}
\|\partial_{t}^{i}(-\mathcal{L}_{\mathbb{H}})^{j/2}(u^{\operatorname{non}}-\overline{u}^{\operatorname{non}})(t)\|_{L^{2}(\mathbb{H}^{n})}\quad \text{and  } \quad\|\partial_{t}^{i}(-\mathcal{L}_{\mathbb{H}})^{j/2}(v^{\operatorname{non}}-\overline{v}^{\operatorname{non}})(t)\|_{L^{2}(\mathbb{H}^{n})},
\end{equation*}
and subsequently following the argument used for \eqref{numinusnv}, we can obtain our required inequality \eqref{ineq2:couple}. This completes the proof of Theorem \ref{global:Application}.
\end{proof}

\begin{remark}
    If we set $v=u,q=p,v_{0}=u_{0}$, and $v_{1}=u_{1}$ in the coupled system \eqref{main:Coupled:Heisenberg}, then proof of Theorem \ref{Main:m>0} will follow from  Theorem \ref{global:Application}.
\end{remark}

\section{Final remarks}\label{sec6}
In this paper, we studied the following fractional Cauchy problem   
\begin{align}\label{New}
	\begin{cases}
		u_{tt}(t, \eta)+\left(-\mathcal{L} \right)^{\alpha} u(t, \eta)+bu_{t}(t,\eta)+mu(t,\eta)
_{\mathbb{H}}		=|u|^p,  &\eta\in \mathbb{H}^{n},t>0, \\
		u(0, \eta)=u_{0}(\eta),\quad u_{t}(0, \eta)=u_{1}(\eta),  &\eta \in  \mathbb{H}^{n},
	\end{cases}
\end{align}
with the damping term $b>0$, the mass term $m\geq 0$,   $\alpha\geq 1$, and  $p>1.$ Our linear decay  estimates for the above system generalize a number of works from the literature available on the  Heisenberg group $\mathbb{H}^{n}$. For example, 
 \begin{itemize}
     \item for $b>0$  and   $m=0$, our findings extends estimates considered by  \cite{Palmieri2020} and   coincidence for the  particular case $\alpha=1$;
          \item for $b>0$  and   $m>0$, our findings can be seen as a generalization of the estimates considered by  \cite{Ruz18} and, in particular, coincidence for   $\alpha=1$. 
\end{itemize}
We also investigated the global (in time) well-posedness of the   Cauchy problem  for the fractional damped wave equation \eqref{New} with or without mass  scenario. More precisely, we obtain the global well-posedness for:
\begin{itemize}
    \item for $m=0$, the Cauchy problem \eqref{New} with  sufficiently small $L^1\cap L^2$ Cauchy data is globally (in time) well-posed for
    \begin{equation*}
\alpha\geq 1,\quad    2\alpha<\mathcal{Q}\leq 4\alpha, \quad    2\leq p\leq 1+\frac{2\alpha}{\mathcal{Q}-2\alpha};
    \end{equation*}
    \item for $m=0$, the Cauchy problem \eqref{New} with  sufficiently small $L^2$ Cauchy data is globally (in time) well-posed for
    \begin{equation*}
\alpha> 1,\quad    2\alpha<\mathcal{Q}< 4\alpha, \quad    1+\frac{4\alpha}{\mathcal{Q}}< p\leq 1+\frac{2\alpha}{\mathcal{Q}-2\alpha};
    \end{equation*}
        \item for $m>0$, the Cauchy problem \eqref{New} with  sufficiently small $L^2$ Cauchy data is globally (in time) well-posed for
    \begin{equation*}
\alpha> 0,\quad    2\alpha<\mathcal{Q}, \quad    1< p\leq 1+\frac{2\alpha}{\mathcal{Q}-2\alpha}.
    \end{equation*}
\end{itemize}

%On the other hand, for the linear viscoelastic damped wave equation on the Heisenberg group  
%\begin{align}\label{main:Heisenberg1}
%	\begin{cases}
%		u_{tt}(t, \eta) -\mathcal{L}_{\mathbb{H}} u(t, \eta)-\mathcal{L}_{\mathbb{H}} u_{t}(t,\eta) 	=0,  &\eta\in \mathbb{H}^{n}, t>0,\\
%		u(0, \eta)=u_{0}(\eta),\quad u_{t}(0, \eta)=u_{1}(\eta),   &\eta \in  \mathbb{H}^{n},
%	\end{cases}
%\end{align}
 %the authors in \cite{LYJ} derived  $L^2-L^2$ decay estimates with additional $L^1$ regularity on initial data for the solution and its higher-order space derivatives.  We expect our approach can also be applied to investigate the following  non-linear structurally damped wave  equations with power type non-linearity $|u|^p$, namely,  
%\begin{align}
%	\begin{cases}
%		u_{tt}(t, \eta)+\left(-\mathcal{L} _{\mathbb{H}}\right)^{\alpha} u(t, \eta)+ \left(-\mathcal{L} _{\mathbb{H}}\right)^{\alpha}u_{t}(t,\eta) 	=|u|^p, & \eta\in \mathbb{H}^{n},t>0, \\
%		u(0, \eta)=u_{0}(\eta),\quad u_{t}(0, \eta)=u_{1}(\eta),   &\eta \in  \mathbb{H}^{n},
%	\end{cases}
%\end{align}
%for $\alpha\geq 1.$  In a similar manner, nonlinear fractional structurally damped wave  equation on the Heisenberg group will be considered in a forthcoming paper.  
\section{Acknowledgement}
SSM is supported by the DST-INSPIRE Faculty Fellowship DST/INSPIRE/04/2023/ 002038. AT is supported by Institute Post-Doctoral fellowship from Tata Institute of Fundamental Research, Centre For Applicable Mathematics, Bangalore, India.  
\section{Declarations}
%\noindent\textbf{Ethical Approval:} Not applicable.\\
\noindent\textbf{Competing interests:} No potential competing of interest was reported by the author.\\
%\textbf{Author's contributions:} 		 All the authors contributed equally.\\
\textbf{Availability of data and materials:} 	The authors confirm that the data supporting the findings of this study are available within the article.
\bibliographystyle{alphaabbr}
\bibliography{time-fractional.bib}

\end{document}